\documentclass[11pt, a4paper]{article}
\usepackage{amsfonts}
\usepackage{amsmath}
\usepackage{amssymb}

\newtheorem{theorem}{Theorem}[subsection]
\newtheorem{definition}{Definition}[subsection]
\newtheorem{proposition}{Proposition}[subsection]
\newtheorem{remark}{Remark}[subsection]
\newtheorem{corollary}{Corollary}[subsection]
\newtheorem{example}{Example}[subsection]

\begin{document}
\title{The geometry of Lie algebroids and its applications to optimal control}
\author{LIVIU POPESCU}
\date{Department of Applied Mathematics, University of
Craiova, Romania, \\e-mail: liviupopescu@central.ucv.ro}
\maketitle
\begin{abstract}
The paper presents the geometry of Lie algebroids and its
applications to optimal control. The first part deals with the
theory of Lie algebroids, connections on Lie algebroids and
dynamical systems defined on Lie algebroids (mainly Lagrangian and
Hamiltonian systems). In the second part we use the framework of
Lie algebroids in the study of distributional systems (drift less
control affine systems) with holonomic or nonholonomic
distributions.
\end{abstract}

\newpage
\tableofcontents
\newpage\
\begin{center}
\textbf{PREFACE}
\quad
\end{center}

The present paper is devoted to Lie algebroids geometry and its
applications to optimal control and variational calculus. The
framework of the differential geometry is very useful in modelling
and understanding of a large class of natural phenomena. The Lie
geometric methods are applied successfully in differential
equations, optimal control theory or theoretical physics. In the
most of cases the study is starting with a variational problem
formulated for a regular Lagrangian (see \cite{Ab}), on the
tangent bundle $TM$ over the manifold $M$ and very often the whole
set of problems is transferred on the dual space $T^{*}M$, endowed
with a Hamiltonian function, via Legendre transformation. The case
of a non-regular Lagrangians is also studied. The problem in this
case is that the proposed Lagrangian formalism yields a singular
Lagrangian description, which makes the Legendre transform
ill-defined and thus no straightforward Hamiltonian formulation
can be related. One of the motivations for the present work is the
study of Lagrangian systems subjected to external constraints
(holonomic or nonholonomic). These systems have a wide application
in many different areas as optimal control theory, mathematical
economics or sub-Riemannian geometry.

In the last years the investigations have led to a geometric framework which
is covering these phenomena. It is precisely the underlying structure of a
Lie algebroid on the phase space which allows a unified treatment. This idea
was first introduced by A. Weinstein \cite{We2, Ca} in order to define a
Lagrangian formalism which is very useful for the various types of such
systems.

The concept of Lie algebroids have been introduced into differential
geometry since the early 1950, and also can be found in physics and algebra,
under a wide variety of names. However, the fundamental concept has been
introduced in sixties by J. Pradines \cite{Pr} in relation with Lie
groupoids. For every Lie groupoid there exists an associated Lie algebroid,
like as for every Lie group there exists an associated Lie algebra. A Lie
algebroid \cite{Mk1, Mk3} over a smooth manifold $M$ is a real vector bundle
$(E,\pi ,M)$ with a Lie algebra structure on its space of sections, and an
application $\sigma ,$ named anchor, which induces a Lie algebra
homomorphism from sections of $E$ to vector fields on $M$. It is convenient
to think a Lie algebroid as a substituent for the tangent bundle of $M$, an
element $e$ of $E$ as a generalized velocity, and the actual velocity $v$ on
$TM$ is obtained when applying the anchor to $e$, i.e., $\sigma (e)=v$.

The basic example of Lie algebroid over the manifold $M$ is the tangent
bundle $TM$ itself, with the identity mapping as anchor. Every integrable
distribution of $TM$ is a Lie algebroid with the inclusion as anchor and
induced Lie bracket, and every Lie algebra is a Lie algebroid over one
point. An important Lie algebroid is the cotangent bundle of a Poisson
manifold \cite{Kos}. Being related to many areas of geometry, as connections
theory \cite{Mk3, Fe2, Cl, Co1, Hr3, Me1, Pp3, Po11} cohomology \cite{Mk1,
Mk3} foliations and pseudogroups, symplectic and Poisson geometry \cite{Lic,
We1, Va1, Fe1, Va2, Kos1, Cra2, Crs, Po13, Po15, Po23} the Lie algebroids
are today the object of extensive studies. More precisely, Lie algebroids
have applications in mechanical systems and optimal control theory \cite{Co,
Me2, Me3, Hr3, Cle, Gra, An2, Po8, Po12, Po14, Bla} (distributional systems)
and are a natural framework in which one can be developed the theory of
differential operators (exterior derivative and Lie derivative) and
differential equations.

In his papers \cite{Mk1, Mk3} K. Mackenzie has been achieved a unitary study
of Lie groupoids and algebroids and together with P. Higgins \cite{Hi} have
introduced the notion of prolongation of a Lie algebroid over a smooth map,
useful in the study of induced vector bundle by the Lie algebroid structure.
Using the geometry of Lie algebroids, A. Weinstein \cite{We2} shows that is
possible to give a common description of the most interesting classical
mechanical systems. He developed a generalized theory of Lagrangian
mechanics and obtained the equations of motions, using the Poisson structure
on the dual of a Lie algebroid and Legendre transformation associated with a
regular Lagrangian. In the last years the problems raised by A. Weinstein
have been investigated by many authors. Thus, E. Martinez \cite{Ma1, Ma2,
Ma4} obtained the same Euler-Lagrange equations using the symplectic
formalism for Lagrangian and Hamiltonian, similarly with the J. Klein
formalism \cite{Kl} for the classical Lagrangian mechanics.

In the classical version of the tangent bundle ($E=TM)$ the
Klein's method is based on the vector bundle structure of $TM$ and
the existence of a vector-valued $1$-form. Such a form does not
exist for a general Lie algebroid \cite{Li} because of different
dimensions of the horizontal and vertical distributions, and so
Klein's approach is not applied directly. To overcome this
difficulty, E. Martinez, M de Leon, J. C. Marero \cite{Ma1, Le}
have proposed a modified version, in which the bundles tangent to $E$ and $%
E^{*}$ are replaced by the prolongations $\mathcal{T}E$ and $\mathcal{T}%
E^{*} $ (in sense Higgins and Mackenzie \cite{Hi}). The nonholonomic
Lagrangian systems and Hamiltonian mechanics on Lie algebroids are studied
by a group of E. Martinez \cite{Le}. The first step in studying the
mechanical control systems on Lie algebroids seems to be done by J. Cortes
and E. Martinez \cite{Co}, which also approached the problem of
accessibility and controllability. A framework for nonholonomic systems,
using a subbundle of a Lie algebroids is proposed by T. Mestdag and B.
Langerock \cite{Me2}. A start in the study of some problems of control
affine systems and sub-Riemannian geometry, using the framework of Lie
algebroids is due to D. Hrimiuc and L. Popescu \cite{Hr3, Po12, Po14}.

Control theory is splitting in two major branches: the first is
the control theory of problems described by partial differential
equations where the objective functionals are mostly quadratic
forms, and the second is the control theory of problems described
by the parameter dependent ordinary differential equations. In
this last case it is more frequent to deal with non-linear systems
and non-quadratic objective functional. The mathematical models
from the optimal control theory cover also the economic growth in
both open and closed economies, exploitation of (non-) renewable
resources, pollution control, behavior of firms or differential
games \cite{Fei, Se1, Se2}.

The geometric methods in the control theory have been applied by many
authors (see \cite{Br, Ju, Bl, Ma3}). One of the most important issues in
the geometric approach is the analysis of the solution to the optimal
control problem as provided by Pontryagin's Maximum Principle; that is, the
curve $c(t)=(x(t),u(t))$ is an optimal trajectory if there exists a lifting
of $x(t)$ to the dual space $(x(t),p(t))$ satisfying the Hamilton equations,
together with a maximization condition for the Hamiltonian with respect to
the control variables $u(t)$.

In the paper \cite{Ma3} E. Martinez presents the Pontryagin
Maximum Principle on Lie algebroids using the prolongation (in
sense of Higgins and Mackenzie \cite{Hi}) of the Lie algebroid
over the vector bundle projection of a dual bundle. In this paper
we study some distributional systems with positive homogeneous
cost, using the Pontryagin Maximum Principle at the level of a Lie
algebroid.

''In spite of that, the control theory can be considered part of the general
theory of differential equations, the problems that inspires it and some of
the results obtained so far, have configured a theory with a strong and
definite personality, that is already offering interesting returns to its
ancestors. For instance, the geometrization of non-linear affine-input
control theory problems by introducing Lie-geometrical methods into its
analysis, started already by R, Brocket \cite{Br}, is inspiring classical
Riemannian geometry and creating what is called today sub-Riemannian
geometry'' \cite{Str1, Mo, Bel, Ag, Cal, Be1, Be2, Be3, Be4, Be5}.

If $M$ is a smooth $n$-dimensional manifold then a sub-Riemannian structure
on $M$ is a pair ($D,g$) where $D$ is a distribution of rank $m$ and $g$ is
a Riemannian metric on $D$. A piecewise smooth curve on $M$ is called
horizontal if its tangent vectors are in $D$. The length of a horizontal
curve $c$ is defined by
\begin{equation}
L(c)=\int_I\sqrt{g(\dot c(t))}dt,  \tag{1}
\end{equation}
where $g$ is a Riemannian metric on $D$. The distance between two points $a$
and $b$ is $d(a,b)=infL(c)$, where the infimum is taken over all horizontal
curves connecting $a$ to $b$. The distance is assumed to be infinite if
there is no horizontal curve that connects these two points. If locally, the
distribution $D$ of rank $m$ is generated by $X_i$, $i=\overline{1,m}$ a
sub-Riemannian structure on $M$ is locally given by a control system
\begin{equation}
\dot x=\sum_{i=1}^mu_i(t)X_i(x),  \tag{2}
\end{equation}
of constant rank $m$, with the controls $u(.)$. The controlled paths are
obtained by integrating the system (2) and are the geodesics in the
framework of sub-Riemannian geometry. If $D$ is assumed to be bracket
generating, \textit{i.e}. sections of $D$ and iterated brackets span the
entire tangent space $TM$, by a well-known theorem of Chow \cite{Ch} the
system (2) is controllable, that is for any two points $a$ and $b,$ there
exists a horizontal curve which connects these points ($M$ is assumed to be
connected).

The concept of sub-Riemannian geometry can be extended to a more
general setting, \cite{Hr3, Cl1, Cl2} by replacing the Riemannian
metric with a posi-tive homogeneous one. For the theory of optimal
control this extension is equivalent to the change of the
quadratic cost of a control affine system with a positive
homogeneous cost. Also, the case of distribution $D$ with
non-constant rank is generating interesting examples (Grushin case
\cite{Fa, Hr3}).

The case when the distribution $D$ generated by vector fields $X_i$, $i=%
\overline{1,m}$ is integrable is also studied. In this case the distribution
determines a foliation on $M$ and two points can be joined if and only if
belongs to the same leaf. In order to find the optimal trajectory of the
system one uses the Pontryagin Maximum Principle at the level of Lie
algebroids, built different in the case of holonomic or nonholonomic
distributions.

\begin{center}
* * *
\end{center}

\newpage\

The paper is organized in two parts. The first part entitled
\textit{The geometry of Lie algebroids} contains seven chapters.
In the first chapter some preliminaries concerning geometrical
structures on the total space of a vector bundle are presented
\cite{Mi3}. We focus on the notions of nonlinear connection and
covariant derivative. In the next chapter we present the notion of
\textit{Lie algebroid} including the cohomology and structure
equations \cite{Mk1}. The notion of prolongation of a Lie
algebroid over the vector bundle projection is studied in the
chapter three. The Ehresmann nonlinear connection
$\mathcal{N}=-\mathcal{L}_{\mathcal{S}}J$ with the coefficients
given by
\[
\mathcal{N}_\alpha ^\beta =\frac 12\left( -\frac{\partial \mathcal{S}^\beta
}{\partial y^\alpha }+y^\varepsilon L_{\alpha \varepsilon }^\beta \right) ,
\]
is investigated and the relations with the Ehresmann connections on tangent
bundles $TE$ and $TM$ are pointed out. In the chapter four we introduce the
notion of dynamical covariant derivative and metric nonlinear connection at
the level of the Lie algebroid $\mathcal{T}E$. The Lagrangian formalism on
Lie algebroids yields a canonical semispray \cite{Ma2}
\[
\mathcal{S}^\varepsilon =g^{\varepsilon \beta }\left( \sigma
_\beta ^i\frac{
\partial L}{\partial x^i}-\sigma _\alpha ^i\frac{\partial ^2L}{\partial
x^i\partial y^\beta }y^\alpha -L_{\beta \alpha }^\theta y^\alpha
\frac{
\partial L}{\partial y^\theta }\right) ,
\]
and a canonical Ehresman connection, which is a metric nonlinear connection.
We also have the Lagrange equations on Lie algebroids given by \cite{We2}
\[
\frac{dx^i}{dt}=\sigma _\alpha ^iy^\alpha ,\quad \frac d{dt}\left(
\frac{
\partial L}{\partial y^\alpha }\right) =\sigma _\alpha ^i\frac{\partial L}{
\partial x^i}-L_{\alpha \beta }^\theta y^\beta \frac{\partial L}{\partial
y^\theta }.
\]
In the case of positive homogeneous Lagrangian (Finsler function) we find a
canonical Ehresmann connection which depends only on Finsler function and
the structure functions of the Lie algebroid.

In the chapter five we deal with the prolongation of a Lie
algebroid over the vector bundle projections of a dual bundle. We
introduce the notions of dual adapted tangent structure
$\mathcal{J}$ and $\mathcal{J}$-regular sections. These structures
induce a canonical nonlinear connection $\mathcal{
N}=-\mathcal{L}_\rho \mathcal{J}$ with the coefficients given by
\cite{Hr4}
\[
\mathcal{N}_{\alpha \beta }=\frac 12\left( t_{\alpha \gamma
}\frac{\partial \rho _\beta }{\partial \mu _\gamma }-\sigma
_\alpha ^it_{\gamma \beta }\frac{
\partial \xi ^\gamma }{\partial q^i}-\rho (t_{\alpha \beta })+\xi ^\gamma
t_{\lambda \beta }L_{\gamma \alpha }^\lambda \right) .
\]
In the case of Hamiltonian formalism these coefficients become \cite{Po19}\
\begin{eqnarray}
\mathcal{N}_{\alpha \beta } &=&\frac 12(\sigma _\gamma
^i\{g_{\alpha \beta }, \mathcal{H}\}-\frac{\partial
^2\mathcal{H}}{\partial q^i\partial \mu _\varepsilon }(\sigma
_\beta ^ig_{\alpha \varepsilon }+\sigma _\alpha
^ig_{\beta \varepsilon })+  \nonumber \\
&&\ \ \ \ +\ \mu _\gamma L_{\varepsilon \kappa }^\gamma
\frac{\partial \mathcal{H}}{\partial \mu _\varepsilon
}\frac{\partial g_{\alpha \beta }}{
\partial \mu _\kappa }+\mu _\gamma L_{\alpha \beta }^\gamma +\frac{\partial
\mathcal{H}}{\partial \mu _\delta }(g_{\alpha \varepsilon }L_{\delta \beta
}^\varepsilon +g_{\beta \varepsilon }L_{\delta \alpha }^\varepsilon )),
\nonumber
\end{eqnarray}
where $\{\cdot ,\cdot \}$ is the Poisson bracket. The corresponding Hamilton
equations on Lie algebroid are given by \cite{We2, Le}
\[
\frac{dq^i}{dt}=\sigma _\alpha ^i\frac{\partial
\mathcal{H}}{\partial \mu _\alpha },\quad \frac{d\mu _\alpha
}{dt}=-\sigma _\alpha ^i\frac{\partial \mathcal{H}}{\partial
q^i}-\mu _\gamma L_{\alpha \beta }^\gamma \frac{
\partial \mathcal{H}}{\partial \mu _\beta }.
\]
In the chapter six we introduce the notion of dynamical covariant
derivative and metric nonlinear connection at the level of a Lie
algebroid $\mathcal{T} E^{*}$. We prove that the canonical
nonlinear connection induces by a regular Hamiltonian is a unique
metric and symmetric nonlinear connection. In the chapter seven we
investigate some aspects of the Lie algebroids geometry endowed
with a Poisson structures, the so-called Poisson-Lie algebroids.

Author's papers \cite{Hr3, Hr4, Po1, Po8, Po10, Po11, Po12, Po13, Po15,
Po16, Po19, Po20, Po21, Po22, Po23, Po25} are used in writting this part.

The purpose of the second part entitled \textit{Optimal Control} is to study
the drift less control affine systems (distributional systems) with positive
homogeneous cost, using the Pontryagin Maximum Principle at the level of a
Lie algebroid in the case of constant rank of distribution.

We prove that the framework of Lie algebroids is better than cotangent
bundle in order to solve some problems of drift less control affine systems.
In the first chapter the known results on the optimal control systems are
recalled by geometric viewpoint. In the next chapter the distributional
systems are presented and the relation between the Hamiltonians on $E^{*}$
and $T^{*}M$ is given by
\[
H(p)=\mathcal{H}(\mu ),\quad \mu =\sigma ^{\star }(p),\quad p\in
T_x^{*}M,\quad \mu \in E_x^{*}.
\]
We investigate the cases of holonomic and nonholonomic
distributions with constant rank. In the holonomic case, we will
consider the Lie algebroid being just the distribution whereas in
the nonholonomic case (i.e., strong bracket generating
distribution) the Lie algebroid is the tangent bundle with the
basis given by vectors of distribution completed by the first Lie
brackets. Also, the case of distribution $D$ with non-constant
rank is studied in the last two sections of the chapter and some
interesting examples are given. In the last chapter we present the
intrinsic relation between the distributional systems and
sub-Riemannian geometry. Thus, the optimal trajectory of our
distributional systems are the geodesics in the framework of
sub-Riemannian geometry. We investigate two classical cases:
Grusin plan and Heisenberg group, but equipped with positive
homogeneous costs (Randers metric). We are using the Pontryagin
Maximum Principle at the level of Lie algebroids, in the case of
Heisenberg group and show that this idea is very useful in order
to solve a large class of distributional systems. Author's papers
\cite{Hr4, Po12, Po14, Po17, Po18, Po22, Po24} are used in writing
this part.

In my opinion, the paper is useful to a large class of readers:
graduate students, mathematicians and to everybody else interested
in the subject of differential geometry, differential equations or
optimal control. I want to address my thanks to all authors
mentioned in this paper and to everybody else I forgot to mention,
without any intention, in the Bibliography.

Finally, I wish to address my thanks to the referees for many
useful remarks and suggestions concerning this paper. I should
like to express the deep gratitude to professor D. Hrimiuc for the
collaboration during the postdoctoral fellowship at the University
of Alberta, Edmonton, Canada, where many ideas presented in this
paper have been started. Also, I want to address my thanks to
Professor P. Stavre for support and guidance given me in life and
in mathematics.

\textit{Acknowledgments}: This work was supported by the strategic grant
POSDRU/89/1.5/S/61968, Project ID61968 (2009), co-financed by the European
Social Fund within the Sectorial Operational Program Human Resources
Development 2007-2013.

$
\begin{array}{c}
\\
\end{array}
$

Craiova,\qquad \qquad \qquad \qquad \qquad \qquad \qquad \qquad \qquad
\qquad  December 2012

$
\begin{array}{c}
\end{array}
$\qquad \qquad \qquad

Assoc. Prof. Liviu Popescu

University of Craiova

Department of Applied Mathematics

13 ''Al. I. Cuza'' st., Craiova 200585

e-mail: liviupopescu@central.ucv.ro

\qquad \quad \ \newpage\ \

\section{\textbf{\ THE GEOMETRY OF LIE ALGEBROIDS}}

\quad \\

The purpose of this first part is to study the geometry of a Lie algebroid
and its prolongations over the vector bundles projections. A Lie algebroid
\cite{Mk1, Mk3} over a smooth manifold $M$ is a real vector bundle $(E,\pi
,M)$ with a Lie algebra structure on its space of sections, and an
application $\sigma ,$ named the anchor, which induces a Lie algebra
homomorphism from the sections of $E$ to vector fields on $M$. For this
reason, in the first chapter we present some results on the geometry of the
total space of a vector bundle, including nonlinear connections and
covariant derivatives. In the next chapter we give only the relevant
formulas for Lie algebroid cohomology we shall need later, and refer the
reader to the monograph \cite{Mk1} for further details.

The chapter three deals with the prolongation $\mathcal{T}E$ of a Lie
algebroid over the vector bundle projection. We introduce the Ehresmann
nonlinear connection on the Lie algebroid $\mathcal{T}E$ and study its
properties \cite{Po11, Po8}. We show that the vertical part of the Lie
brackets of horizontal sections from the basis represents the components of
the curvature tensor of the nonlinear connection. We study the related
connections and show that a connection on the tangent bundle $TE$ induces a
connection on the Lie algebroid $\mathcal{T}E$. We introduce an almost
complex structure on Lie algebroids and prove that its integrability is
characterized by zero torsion and curvature property of the connection. We
present the notion of dynamical covariant derivative at the level of a Lie
algebroid and show that the metric compatibility of the semispray and
associated nonlinear connection gives the one of the so called Helmholtz
conditions of the inverse problem of Lagrangian Mechanics. In the
homogeneous case a canonical nonlinear connection associated to a Finsler
function is determined. We study the linear connections on $\mathcal{T}E$
and determine the torsion and curvature.

In the chapter four we study the dynamical covariant derivative and metric
nonlinear connection on $\mathcal{T}E$ \cite{Po20}. We introduce the
dynamical covariant derivative as a tensor derivation and study the
compatibility conditions with a pseudo-Riemannian metric. In the case of
SODE connection we find the expression of Jacobi endomorphism and its
relation with curvature tensor. We prove that the canonical nonlinear
connection induced by a regular Lagrangian is a unique connection which is
metric and compatible with symplectic structure. Also the invariant form of
the Helmholtz conditions on Lie algebroids are given (see also \cite{Cr1}).

The chapter five deals with the prolongation $\mathcal{T}E^{*}$ of a Lie
algebroid over the vector bundle projection of a dual bundle. We study the
properties of the connections on $\mathcal{T}E^{*}$ \cite{Hr4, Po10, Po12}
and introduce the notions of adapted almost tangent structure, almost
complex structure and characterize the integrability conditions in terms of
torsion and curvature of the connection. We prove that every $\mathcal{J}$%
-regular section (in particular, any regular Hamiltonian on $E^{*}$)
determines a canonical Ehresmann connection on the Lie algebroid $\mathcal{T}%
E^{*}.$ We introduce some generalizations of the Hamilton sections, as a
mechanical structures and semi-Hamiltonian sections and study their
properties. In the last part of this chapter, using the diffeomorphism from $%
\mathcal{T}E^{*}$ and $\mathcal{T}E$ \cite{Le} induced by a regular
Hamiltonian, we can transfer many geometrical results between these spaces.
Thus, a semispray on $\mathcal{T}E$ is transformed into a semi-Hamiltonian
section on $\mathcal{T}E^{*}$ if and only if the nonlinear connection on $%
\mathcal{T}E$ determined by semispray is just the canonical nonlinear
connection induced by regular Lagrangian, via Legendre transformation.

In the chapter six we study the dynamical covariant derivative and metric
nonlinear connection on $\mathcal{T}E^{*}$ \cite{Po25}. Using the notion of $%
\mathcal{J}$-regular section we introduce the dynamical covariant derivative
as a tensor derivation. In the case of nonlinear connection induce by a $%
\mathcal{J}$-regular section we find the expression of Jacobi endomorphism
and its relation with curvature tensor. Finally, we prove that the canonical
nonlinear connection induced by a regular Hamiltonian is the unique metric
and symmetric nonlinear connection.

In the chapter seven we investigate some aspects of the Lie algebroids
geometry endowed with a Poisson structures \cite{Po13, Po15, Po23}, which
generalize the Poisson manifolds. We recall the Cartan calculus and the
Schouten-Nijenhuis bracket at the level of Lie algebroids and introduce the
Poisson structure on Lie algebroids. We study the properties of linear
contravariant connection and its tensors of torsion and curvature. In the
last part of this section we find a Poisson connection which depends only on
the Poisson bivector and structural functions of Lie algebroid, which
generalize some results of Fernandes from \cite{Fe1}. Also the geodesic
equations are given. We study the properties of the complete lift of a
Poisson bivector on $\mathcal{T}E$ and introduce the notion of horizontal
lift. The compatibility conditions of these bivectors are investigated.
Finally, the compatibility conditions between the canonical Poisson
structure and the horizontal lift on $\mathcal{T}E^{*}$ are given.

\newpage\

\subsection{\textbf{The geometry of the total space of a vector bundle}}

\subsubsection{\textbf{Connections on vector bundles}}

The connections theory is an important topic of the differential
geometry with important applications in Differential Equations or
Optimal Control. In this section we shall present only the notion
of Ehresmann nonlinear connection and induced geometrical
structures. Also, the covariant derivative is described. For more
details and complete proofs we refer to the monographs \cite{Ko,
LeR, Mi, Mi3}.

Let us consider the $n$-dimensional differentiable manifold $M$
and a vector bundle $(E,\pi ,M)$ over $M$, with type fibre
$F=\Bbb{R}^m$. The structure of the vector bundle is given by a
vectorial atlas $\{(U_i,\psi _i,\Bbb{R} ^m)\}_{i\in I}$ such that

1) $(U_i)_{i\in I}$ is an open covering of the manifold $M$.

2) The mappings $\psi _i:\pi ^{-1}(U_i)\rightarrow U_i\times \Bbb{R}^m$ are
bijective and satisfy the relation
\[
\pi \left( \psi _i^{-1}(x,f)\right) =x,\quad x\in M,\ f\in \Bbb{R}^m.
\]

3) For every pair $(i,j)\in I\times I,$ such that $U_i\cap U_j\not
=\emptyset $ there exists a smooth mapping $g_{ij}:U_i\cap
U_j\rightarrow GL(m,\Bbb{R})$ with $\Psi _{i,x}^{-1}=\Psi
_{j,x}^{-1}\circ g_{ji}(x)$ for every $x\in U_i\cap U_j$ where
$\Psi _{i,x}^{-1}$ is the restriction of $ \psi _i^{-1}$ to
$\{x\}\times \Bbb{R}^m.$

Let $\{(U_i,\varphi _i)\}_{i\in I}$ be an atlas on the manifold
$M$ such that $U_i$ belongs to the maps domain into vectorial
atlas $\{(U_i,\psi _i, \Bbb{R}^m)\}_{i\in I}$. We obtain that
$\{(\pi ^{-1}(U_i),h_i)\}_{i\in I}$ with
\[
h_i:\pi ^{-1}(U_i)\rightarrow \Bbb{R}^n\times \Bbb{R}^m,\quad
h_i(u)=(\varphi _i(\pi (u)),\psi _{i,\pi (u)}(u)),
\]
is a differentiable atlas on the manifold $E$. If $(U_j,\varphi _j)$ is
another local chart on $M$ with $U_i\cap U_j\not =\emptyset $ and $(U_j,\psi
_j,\Bbb{R}^m)$ is a bundle chart, then
\[
(h_j\circ h_i^{-1})(x,y)=\left( (\varphi _j\circ \varphi
_i)(x),g_{ji}(\varphi _i(x)y\right) ,\quad (x,y)\in
\Bbb{R}^n\times \Bbb{R} ^m.
\]
We will denote by $M_a^{a^{\prime }}(x)$ the entries of the matrix
associated with the linear application $g_{ji}(\varphi _i(x))$.
For $x\in U\subset M$ we take $\varphi _i(x)=(x^i)\in
\Bbb{R}\mathbf{^n},$ $i= \overline{1,n}$ and considering $\varphi
_j(x)=(x^{i^{\prime }})$ then $ \varphi _j\circ \varphi _i^{-1}$
has the form
\[
x^{i^{\prime }}=x^{i^{\prime }}(x^1,...,x^n),\text{ }rank\left(
\frac{
\partial x^{i^{\prime }}}{\partial x^i}\right) =n.
\]
A tangent vector $X_x$ at the point $x\in M$ will be locally represented by
the pair $(x^i,y^i)$, where $(y^i)$ is given by $X_x=y^i\frac \partial
{\partial x^i}$. Therefore, the transformations of coordinates on
differentiable manifold $TM$ have the form
\begin{equation}
x^{i^{\prime }}=x^{i^{\prime }}(x^1,...,x^n),\text{ }rank\left(
\frac{
\partial x^{i^{\prime }}}{\partial x^i}\right) =n,\quad y^{i^{\prime }}=
\frac{\partial x^{i^{\prime }}}{\partial x^i}y^i.  \tag{1.1.1}
\end{equation}
A point $u$ $\in E$ on the total space of a vector bundle $(E,\pi ,M)$ is
locally represented by the pair $(x^i,y^a)\in \Bbb{R}^n\times \Bbb{R}^m$ .
The transformations of coordinates $(x^i,y^a)\rightarrow (x^{i^{\prime
}},y^{a^{\prime }})$ on $E$ are of the form \cite{Mi3}
\begin{equation}
\begin{array}{c}
x^{i^{\prime }}=x^{i^{\prime }}(x^1,...,x^n),\text{ }rank\left(
\frac{
\partial x^{i^{\prime }}}{\partial x^i}\right) =n, \\
\\
y^{a^{\prime }}=M_a^{a^{\prime }}(x)y^a,\quad \text{ }rank(M_a^{a^{\prime
}}(x))=m.
\end{array}
\tag{1.1.2}
\end{equation}
For the vector bundle $(E,\pi ,M)$ we consider $\pi
_{*}:TE\rightarrow TM$ the tangent application of $\pi $. The
application $\pi _{*}$ is a $\pi $ -morphism of vector bundles
between tangent bundles ($TE,\pi _E,E$) and ($ TM,\pi _M,M)$. The
kernel of this $\pi $-morphism is a subbundle of vector bundle
($TE,\pi _E,E$), denoted $(VE,\pi _V,E)=Ker\ \pi _{*}$, which will
be called the \textit{vertical subbundle}. The total space is
$VE={\cup }V_u$ , where $V_u=Ker\ \pi _{*},$ $u\in E.$

A tangent vector $X_u$ at the point $u\in E$ has the local
representation $ (x^i,y^a,X^i,A^a),$ where the coefficients
$\left( X^i\right) \in \Bbb{R}^n$ and $(A^a)\in \Bbb{R}^m$ are
defined by the relation $X_u=X^i\frac \partial {\partial
x^i}+A^a\frac \partial {\partial y^a}.$ The tangent application $
\pi _{*}$ is locally represented by $\pi _{*}(x,y,X,A)=(x,X).$
Consequently, the local fibres of vector bundle $(TE,\pi _{*},TM)$
are isomorph with $ \{x\}\times \Bbb{R}^m\times \{X\}\times
\Bbb{R}^m\simeq \Bbb{R}^{2m}$. Because $\pi _{*}(\frac \partial
{\partial y^a})=0$, it results that the functions $\{\frac
\partial {\partial y^a}\},a=\overline{1,m}$ determine a local
basis of the vertical distribution $\{u\rightarrow V_u\mid u\in
E\}$, which means that $V_u$ is integrable. The elements of the
vertical subbundle have de form $(x,y,0,A)$, that is the fibres of
vertical subbundle $VE$ are locally isomorph with $\Bbb{R}^m$. Let
$\pi ^{*}TM$ be the induced vector bundle of tangent bundle over
the application $\pi :E\rightarrow M$ and $\pi !:TE\rightarrow \pi
^{*}TM$ given by
\[
\pi !(X_u)=(u,\pi _{*}(X_u)).
\]
This application is a morphism of vector bundles, and follows that the
application $\pi !$ is a surjection and
\[
Ker\ \pi !=Ker\ \pi _{*}=VE.
\]
Therefore, it results that the following sequence of vector bundles over $E$
is exact
\begin{equation}
0\rightarrow VE\stackrel{i}{\rightarrow }TE\stackrel{\pi !}{\rightarrow }\pi
^{*}TM\rightarrow 0.  \tag{1.1.3}
\end{equation}
where $i:VE\rightarrow TE$ is the inclusion map.

We can present now a definition of the Ehresmann connection, called usually
nonlinear connection.

\begin{definition}
The nonlinear connection in the vector bundle $(E,\pi ,M)$ is a splitting on
the left of the exact sequences (1.1.3).
\end{definition}

It results that a nonlinear connection in $(E,\pi ,M)$ is a morphism of
vector bundle $N:TE\rightarrow VE$ such that $N\circ i=Id\left|
_{VE}.\right. $ The kernel of the morphism $N$ is a vector subbundle of the
bundle ($TE,\pi _E,E$), which will be called \textit{horizontal subbundle}
and will be denoted $(HE,\pi _H,E)$. Therefore, the vector bundle $(TE,\pi
_E,E)$ is the Whitney sum of the horizontal and vertical subbundles. Thus we
have the following characterization of the nonlinear connection.

\begin{proposition}
A nonlinear connection in $(E,\pi ,M)$ is determined by the existence of the
vector subbundle $(HE,\pi _H,E)$ of the tangent bundle over $E,$ ($TE,\pi
_E,E$) such that $TE=HE\oplus VE.$
\end{proposition}

The restriction $\pi !\left| _{HE}\right. $ of the application$\ \pi !$ to
the horizontal subbundle $HE$ is an isomorphism of vector bundles. The
component $\pi _{*}:HE\rightarrow TM$ of the application $\pi !\mid _{HE}$
is a $\pi $-morphism and its restriction to the fibres is an isomorphism.
Therefore, for any vector field $X$ on $M$, there exists a horizontal vector
field on $E$, such that $\pi _{*}(X^h)=X.$ The vector field $X^h$ is called
the \textit{horizontal lift }of the vector field $X.$ The horizontal lift
has a local representation
\[
\left( \frac \partial {\partial x^i}\right) ^h=\frac \delta {\delta
x^i},\quad i=\overline{1,n},
\]
thus, we determine a local basis \{$\frac \delta {\delta x^i}$\} of $H_uE$.
This vector field can be represented in the form $\frac \delta {\delta
x^i}=A_i^j\frac \partial {\partial x^j}+B_i^a\frac \partial {\partial y^a},$
but the condition $\pi _{*}\left( \frac \delta {\delta x^i}\right) =\frac
\partial {\partial x^i}$ implies $A_j^i=\delta _j^i$ and the fact that $%
\frac \delta {\delta x^i}$ are the kernel of the mapping $N$ gives $%
B_i^a=-N_i^a$. It results
\begin{equation}
\frac \delta {\delta x^i}=\frac \partial {\partial x^i}-N_i^a\frac \partial
{\partial y^a},  \tag{1.1.4}
\end{equation}
We obtain a new basis $\left( \frac \delta {\delta x^i},\frac \partial
{\partial y^a}\right) $ of the tangent bundle $T_uE$ which is called the
\textit{Berwald basis }associated to the nonlinear connection $N$.

If we denote $N_i^a$, $i=\overline{1,n},$ $a=\overline{1,m}$ the
coefficients of a nonlinear connection, then it results \cite{Mi3}:

\begin{proposition}
By a change of the local coordinates (1.1.2) on the vector bundle $(E,\pi
,M),$ the local coefficients $N_i^a$ of a nonlinear connection $N$ change as
follows
\end{proposition}

\begin{equation}
N_{i^{\prime }}^{a^{\prime }}\frac{\partial x^{i^{\prime
}}}{\partial x^i} =M_a^{a^{\prime }}N_i^a-\frac{\partial
M_{a^{\prime }}^a}{\partial x^i}y^a. \tag{1.1.5}
\end{equation}

A nonlinear connection in the vector bundle $(E\backslash \{0\},\pi ,M)$ is
said to be homogeneous (respective linear) if the coefficients $N_i^a(x,y)$
are homogeneous (respective linear) with respect to the second argument.

We define a morphism $v:\mathcal{X}(E)\rightarrow \mathcal{X}(E)$
such that $ v(X)=-N(X)$ if $X\in \Gamma (VE)$ and $v(X)=0$, for
$X\in \Gamma (HE)$. It follows that the morphism
$v:\mathcal{X}(E)\rightarrow \mathcal{X}(E)$ has the properties:

1$^{\circ }$ $v(\mathcal{X}(E))\subset \Gamma (VE),$

2$^{\circ }$ $\{v(X)=X\}\Leftrightarrow X\in $ $\Gamma (VE),$\\and it
determines a nonlinear connection in $(E,\pi ,M)$ since $Kerv$ $=\Gamma
(HE). $

It results that $v^2=v,$ \textit{i.e}. $v$ is a projector, which is called
the \textit{vertical projector} of the nonlinear connection $N$.
Analogously, a nonlinear connection in the vector bundle $(E,\pi ,M)$ is
characterized by the morphism $h:\mathcal{X}(E)\rightarrow \mathcal{X}(E)$
with the properties $h^2=h,$ $Ker\ h=\Gamma (VE),$ and $h$ is called the
\textit{horizontal projector} of the nonlinear connection $N.$ It results
that $h+v=Id.$ In \cite{Mi3} one proves:

\begin{proposition}
A nonlinear connection in a vector bundle $(E,\pi ,M)$ is
characterized by an almost product structure $P$ on $E$ whose
distribution of eigensubspaces which correspond to the eigenvalue
$-1$ coincides to the vertical distribution.
\end{proposition}

From the previous consideration it follows that
\[
P=2h-Id=Id-2v=h-v,P^2=Id.
\]
The existence of the nonlinear connection in $(E,\pi ,M)$ leads to the
decomposition
\[
\mathcal{X}(E)=\Gamma (HE)\oplus \Gamma (VE)\Rightarrow X=hX+vX,\ \forall
X\in \mathcal{X}(E).
\]
\newpage\

\subsubsection{\textbf{Covariant derivative}}

We consider a local basis $\{s_a\}$ of the sections of the bundle
$\pi ^{-1}(U)\rightarrow U$, $U\subset M$ and it results $\rho
=\rho ^a(x)s_a$, $ \rho \in \Gamma (E)$. The nonlinear connection
$N$ induces a covariant derivative for the sections in $E$ defined
as follows
\[
D_X\rho =X^i\left( \frac{\partial \rho ^a}{\partial x^i}+N_i^a(x,\rho
(x))X^i\right) s_a,
\]
where
\[
X=X^i\frac \partial {\partial x^i}\in \Gamma (TM).
\]
This covariant derivative has the properties:

1$^{\circ }$ $D_{X+Y}=D_X+D_Y,$

2$^{\circ }$ $D_{fX}\rho =fD_X\rho ,\quad X,Y\in \Gamma (TM).$ \\If $N$ is
homogeneous one also has

3$^{\circ }$ $D_X(f\rho )=X(f)\rho +fD_X\rho ,\quad \rho \in \Gamma (E),\
f\in \mathcal{F}(M)$. \\If $N$ is linear, then the covariant derivative
satisfies 1$^{\circ }$, 2$^{\circ }$, 3$^{\circ }$ and

4$^{\circ }$ $D_X(\rho _1+\rho _2)=D_X\rho _1+D_X\rho _2,\quad \rho _1,\rho
_2\in \Gamma (E)$.

Let $c:[a,b]\rightarrow M$, $t\rightarrow c(t)$ be a smooth curve
on the manifold $M$ and $\stackrel{\cdot }{c}:[a,b]\rightarrow TM$
the tangent field along the curve $c$. Setting $D_{\stackrel{\cdot
}{c}}\rho =\frac{ D\rho }{dt}$ we say that the section $\rho $ in
$(E,\pi ,M)$ is parallel along $c$ if $\frac{D\rho }{dt}=0$ and it
results:

\begin{proposition}
A local section $\rho $ of the bundle $\pi ^{-1}(U)\rightarrow U$, $U\subset
M$ is parallel along the curve $c$ with $c([a,b])\subset U$ if and only if
\[
\frac{d\rho ^a}{dt}+N_i^a(x,\rho (x))\frac{dx^i}{dt}=0,\quad \rho =\rho
^a(x)s_a,\quad x\in U.
\]
\end{proposition}

The decomposition $T_uE=H_uE\oplus V_uE$ in the vector bundle
$(E,\pi ,M)$ leads to the decomposition $(T_uE)^{*}=(V_uE)^{\perp
}\oplus (H_uE)^{\perp }$ , for $u\in E,$ where $(V_uE)^{\perp }$
denotes the subspace of $1$-forms on $T_uE$ which vanish on
horizontal vectors, and $(H_uE)^{\perp }$ denotes the subspace of
$1$-forms on $T_uE$ which vanish on vertical vectors. Therefore,
we have
\[
\omega =h\omega +v\omega ,\quad \forall \omega \in \mathcal{\wedge }^1(E),
\]
where $h\omega (X)=\omega (hX)$ and $v\omega (X)=\omega (vX)$. The
\textit{ dual Berwald basis }is $(dx^i,\delta y^a)$, where
\begin{equation}
\delta y^a=dy^a+N_i^adx^i.  \tag{1.1.6}
\end{equation}

\begin{proposition}
The Berwald basis $\left( \frac \delta {\delta x^i},\frac \partial {\partial
y^a}\right) $ and its dual $(dx^i,\delta y^a)$ transform under a change of
coordinates (1.1.2) as follows
\[
\frac \delta {\delta x^i}=\frac{\partial x^{i^{\prime }}}{\partial x^i}\frac
\delta {\delta x^{i^{\prime }}},\quad \ \frac \partial {\partial
y^a}=M_a^{a^{\prime }}\frac \partial {\partial y^{a^{\prime }}},
\]
\[
dx^i=\frac{\partial x^i}{\partial x^{i^{\prime }}}dx^{i^{\prime }},\ \quad
\delta y^a=M_{a^{\prime }}^a\delta y^{a^{\prime }}.
\]
\end{proposition}

The Lie brackets of the vector fields of Berwald basis are given by
\begin{equation}
\left[ \frac \delta {\delta x^i},\frac \delta {\delta x^j}\right]
=R_{ji}^a\frac \partial {\partial y^a},\quad \left[ \frac \delta
{\delta x^i},\frac \partial {\partial y^a}\right] =\frac{\partial
N_i^b}{\partial y^a }\frac \partial {\partial y^b},\quad \left[
\frac \partial {\partial y^a},\frac \partial {\partial y^b}\right]
=0,  \tag{1.1.7}
\end{equation}
where
\begin{equation}
R_{jk}^a=\frac{\delta N_j^a}{\delta x^k}-\frac{\delta N_k^a}{\delta x^j}.
\tag{1.1.8}
\end{equation}

In the following we deal with the curvature of a nonlinear connection $N$.

\begin{definition}
The curvature of a nonlinear connection $N$ is given by
\begin{equation}
\Omega =-\mathbf{N}_v,  \tag{1.1.9}
\end{equation}
where $v$ is the vertical projector induced by $N$ and $\mathbf{N}_v$ is the
Nijenhuis tensor associated to $v,$%
\[
\mathbf{N}_v(X,Y)=[vX,vY]-v[vX,Y]-v[X,vY]+v^2[X,Y],\quad X,Y\in
\mathcal{X} (E).
\]
\end{definition}

In local coordinates, we set
\[
\Omega =\frac 12\Omega _{ij}^a(dx^i\wedge dx^j)\otimes \frac \partial
{\partial y^a},
\]
and using (1.1.9), we obtain
\begin{equation}
\Omega (hX,hY)=-v[hX,hY],\quad \Omega (hX,vY)=\Omega (vX,vY)=0,  \tag{1.1.10}
\end{equation}
and then
\begin{equation}
\Omega _{ij}^a=-R_{ij}^a=\frac{\delta N_j^a}{\delta
x^i}-\frac{\delta N_i^a}{ \delta x^j}.  \tag{1.1.11}
\end{equation}
The equality $v=Id-h$ leads to $\mathbf{N}_v=\mathbf{N}_h,$ and we
obtain $ \Omega =-\mathbf{N}_h.$

\begin{proposition}
The following equation holds
\[
[X^h,Y^h]=[X,Y]^h-\Omega (X^h,Y^h).
\]
\end{proposition}

Using (1.1.10) follows a characterization of the integrability of the
horizontal distribution.

\begin{theorem}
The horizontal distribution of a nonlinear connection $N$ is integrable if
and only if the curvature vanishes.
\end{theorem}

See \cite{Bu2, Mi3} for the particular case of the tangent bundle.

\newpage\

\subsection{\textbf{\ Lie algebroids}}

\subsubsection{\textbf{Cohomology}}

Let $M$ be a real, $C^\infty $-differentiable, $n$-dimensional
manifold and $ (TM,\pi _M,M)$ its tangent bundle.

\begin{definition}
A Lie algebroid over a manifold $M$ is a triple $(E,[\cdot ,\cdot
]_E,\sigma )$, where ($E,\pi ,M)$ is a vector bundle of rank $m$
over $M,$ which satisfies the following conditions: \\a) $C^\infty
(M)$-module of sections $
\Gamma (E)$ is equipped with a Lie algebra structure $[\cdot ,\cdot ]_E$. \\
b) $\sigma :E\rightarrow TM$ is a bundle map (called the anchor) which
induces a Lie algebra homomorphism (also denoted $\sigma $) from the Lie
algebra of sections $(\Gamma (E),[\cdot ,\cdot ]_E)$ to the Lie algebra of
vector fields $(\mathcal{\chi }(M),[\cdot ,\cdot ])$ satisfying the Leibniz
rule
\begin{equation}
[s_1,fs_2]_E=f[s_1,s_2]_E+(\sigma (s_1)f)s_2,\ \forall s_1,s_2\in \Gamma
(E),\ f\in C^\infty (M).  \tag{1.2.1}
\end{equation}
\end{definition}

From the above definition it results: \\$1^{\circ }$ $[\cdot ,\cdot ]_E$ is
a $\Bbb{R}$-bilinear operation, \\$2^{\circ }$ $[\cdot ,\cdot ]_E$ is
skew-symmetric, i.e.
\[
\lbrack s_1,s_2]_E=-[s_2,s_1]_E,\quad \forall s_1,s_2\in \Gamma (E),
\]
$3^{\circ }$ $[\cdot ,\cdot ]_E$ verifies the Jacobi identity
\begin{equation}
\lbrack s_1,[s_2,s_3]_E]_E+[s_2,[s_3,s_1]_E]_E+[s_3,[s_1,s_2]_E]_E=0,\
\tag{1.2.2}
\end{equation}
and $\sigma $ being a Lie algebra homomorphism, means that
\begin{equation}
\sigma [s_1,s_2]_E=[\sigma (s_1),\sigma (s_2)].  \tag{1.2.3}
\end{equation}

The existence of a Lie bracket on the space of sections of a Lie algebroid
leads to a calculus on its sections analogous to the usual Cartan calculus
on differential forms. In this paragraph we give only the relevant formulas
for Lie algebroid cohomology we shall need later, and refer the reader to
the monograph \cite{Mk1} for further details.

If $f$ is a function on $M$, then $df(x)\in E_x^{*}$ is given by $
\left\langle df(x),a\right\rangle =\sigma (a)f$, for $\forall a\in
E_x$. For $\omega $ $\in $ $\bigwedge^k(E^{*})$ the
\textit{exterior derivative} $ d^E\omega \in
\bigwedge^{k+1}(E^{*})$ is given by the formula
\newpage
\begin{equation*}
d^E\omega(s_1,...,s_{k+1})=\stackrel{k+1}{\sum_{i=1}}(-1)^{i+1}\sigma
(s_i)\omega (s_1,...,\stackrel{\symbol{94}}{s}_i,...,s_{k+1})+
\end{equation*}
\begin{equation}
+\sum_{1\leq i<j\leq k+1}(-1)^{i+j}\omega ([s_{i,}s_j]_E,s_1,...,
\stackrel{\symbol{94}}{s_i},...,\stackrel{\symbol{94}}{s_j},...s_{k+1}).
\tag{1.2.4}
\end{equation}
where $s_i\in \Gamma (E)$, $i=\overline{1,k+1}$, and the hat over
an argument means the absence of the argument. It results that
\begin{equation}
(d^E)^2=0,  \tag{1.2.5}
\end{equation}
\begin{equation}
d^E(\omega _1\wedge \omega _2)=d^E\omega _1\wedge \omega _2+(-1)^{\deg
\omega _1}\omega _1\wedge d^E\omega _2,  \tag{1.2.6}
\end{equation}
The cohomology associated with $d^E$ is called the \textit{Lie
algebroid cohomology} of $E$ and is denoted by $H^{\bullet }(E).$
Also, for $\xi $ $ \in \Gamma (E)$ one can define the \textit{Lie
derivative} with respect to $ \xi $ by
\begin{equation}
\mathcal{L}_\xi =i_\xi \circ d^E+d^E\circ i_\xi ,  \tag{1.2.7}
\end{equation}
where $i_\xi $ is the contraction with $\xi $.

\subsubsection{\textbf{Structure equations on Lie algebroids}}

If we take the local coordinates $(x^i)$ on an open $U\subset $ $M$, a local
basis $\{s_\alpha \}$ of the sections of the bundle $\pi ^{-1}(U)\rightarrow
U$ generates local coordinates $(x^i,y^\alpha )$ on $E$. The local functions
$\sigma _\alpha ^i(x)$, $L_{\alpha \beta }^\gamma (x)$ on $M$ given by
\begin{equation}
\sigma (s_\alpha )=\sigma _\alpha ^i\frac \partial {\partial
x^i},\quad [s_\alpha ,s_\beta ]_E=L_{\alpha \beta }^\gamma
s_\gamma ,\quad i=\overline{ 1,n},\quad \alpha ,\beta ,\gamma
=\overline{1,m},  \tag{1.2.8}
\end{equation}
are called the \textit{structure functions of the Lie algebroid}, and
satisfy the \textit{structure equations }on Lie algebroids
\begin{equation}
\sum_{(\alpha ,\beta ,\gamma )}\left( \sigma _\alpha ^i\frac{\partial
L_{\beta \gamma }^\delta }{\partial x^i}+L_{\alpha \eta }^\delta L_{\beta
\gamma }^\eta \right) =0,  \tag{1.2.9}
\end{equation}
\begin{equation}
\sigma _\alpha ^j\frac{\partial \sigma _\beta ^i}{\partial x^j}-\sigma
_\beta ^j\frac{\partial \sigma _\alpha ^i}{\partial x^j}=\sigma _\gamma
^iL_{\alpha \beta }^\gamma .  \tag{1.2.10}
\end{equation}
Locally, if $f\in C^\infty (M)$ then $d^Ef=\frac{\partial
f}{\partial x^i} \sigma _\alpha ^is^\alpha ,$ where $\{s^\alpha
\}$ is the dual basis of $ \{s_\alpha \}$ and if $\theta \in
\Gamma (E^{*}),$ $\theta =\theta _\alpha s^\alpha $ then
\begin{equation}
d^E\theta =\left( \sigma _\alpha ^i\frac{\partial \theta _\beta }{\partial
x^i}-\frac 12\theta _\gamma L_{\alpha \beta }^\gamma \right) s^\alpha \wedge
s^\beta .  \tag{1.2.11}
\end{equation}
Particularly, we get
\begin{equation}
d^Ex^i=\sigma _\alpha ^is^\alpha ,\quad d^Es^\alpha =-\frac 12L_{\beta
\gamma }^\alpha s^\beta \wedge s^\gamma .  \tag{1.2.12}
\end{equation}
Under a change of coordinates
\begin{equation}
\left\{
\begin{array}{l}
x^{i^{\prime }}=x^{i^{\prime }}(x^i),\quad i,i^{\prime }=\overline{1,n}\ on\
M, \\
y^{\alpha ^{\prime }}=A_\alpha ^{\alpha ^{\prime }}(x^i)y^\alpha ,\quad
\alpha ,\alpha ^{\prime }=\overline{1,m}\ on\ E,
\end{array}
\right.  \tag{1.2.13}
\end{equation}
corresponding to a new base $s_\alpha =A_\alpha ^{\alpha ^{\prime
}}s_{\alpha ^{\prime }}$, then the transformation rules of the structure
functions are
\begin{equation}
\sigma _{\alpha ^{\prime }}^{i^{\prime }}A_\alpha ^{\alpha
^{\prime }}=\frac{
\partial x^{i^{\prime }}}{\partial x^i}\sigma _\alpha ^i,  \tag{1.2.14}
\end{equation}
\begin{equation}
L_{\alpha \beta }^\gamma A_\gamma ^{\gamma ^{\prime }}=L_{\alpha ^{\prime
}\beta ^{\prime }}^{\gamma ^{\prime }}A_\alpha ^{\alpha ^{\prime }}A_\beta
^{\beta ^{\prime }}+\sigma _\alpha ^i\frac{\partial A_\beta ^{\gamma
^{\prime }}}{\partial x^i}-\sigma _\beta ^i\frac{\partial A_\alpha ^{\gamma
^{\prime }}}{\partial x^i}.  \tag{1.2.15}
\end{equation}

Some examples of Lie algebroids which will be used in this paper
(see \cite {Mk3} for more examples and details)

\begin{example}
The tangent bundle $E=TM$ itself, with identity mapping as anchor.
Then $ H^{\bullet }(E)=H_{de\ Rham}^{\bullet }(M)$ is the
\textit{de Rham cohomology. }With respect to the usual coordinates
$(x,\stackrel{\cdot }{x})$ , the structure functions are
$L_{jk}^i=0$, $\sigma _j^i=\delta _j^i$, but if we were to change
to another basis for the vector fields, the structure functions
would become nonzero.
\end{example}

\begin{example}
Any integrable subbundle of $TM$ is a Lie algebroid with the
inclusion as anchor and the induced bracket. So,
$E=T\mathcal{F}\subset TM,$ an involutive distribution associated
with a regular foliation $\mathcal{F}$, where one gets the
\textit{tangential cohomology }denoted $H_{\mathcal{F}+ }^{\bullet
}(M).$
\end{example}

\begin{example}
The cotangent bundle of a Poisson manifold $E=T^{*}M$ with $(M,\Pi
)$ Poisson manifold. These carry a bracket characterized by the
rule $ \{df,dg\}=d\{f,g\}$ and the anchor is the map $\Pi
^{\#}:T^{*}M\rightarrow TM $ associated to the Poisson bivector
field $\Pi .$ We obtain the \textit{ Poisson cohomology }denoted
$H_\Pi ^{\bullet }(M)$ \cite{Va1}.
\end{example}

\newpage\

\subsection{\textbf{The prolongation of a Lie algebroid over the vector
bundle projection}}

Let $(E,\pi ,M)$ be a vector bundle. For the projection $\pi :E\rightarrow M$
we can construct the prolongation of $E$ (see \cite{Hi, Ma1, Ma2, Le}). The
associated vector bundle is ($\mathcal{T}E,\pi _2,E$) where
\[
\mathcal{T}E={\cup}\mathcal{T}_wE, \quad w\in E,
\]
\[
\mathcal{T}_wE=\{(u_x,v_w)\in E_x\times T_wE\mid \sigma (u_x)=T_w\pi
(v_w),\quad \pi (w)=x\in M\},
\]
and the projection $\pi _2(u_x,v_w)=\pi _E(v_w)=w$, where $\pi
_E:TE\rightarrow E$ is the tangent projection.

We also have the canonical projection $\pi
_1:\mathcal{T}E\rightarrow E$ given by $\pi _1(u,v)=u$. The
projection onto the second factor $\sigma ^1:
\mathcal{T}E\rightarrow TE$, $\sigma ^1(u,v)=v$ will be the anchor
of a new Lie algebroid over the manifold $E$. An element of
$\mathcal{T}E$ is said to be vertical if it is in the kernel of
the projection $\pi _1$.

We will denote $(V\mathcal{T}E,\pi _{2\mid _{V\mathcal{T}E}},E)$ the
vertical bundle of $(\mathcal{T}E,\pi _2,E)$ and
\[
\sigma ^1\left| _{V\mathcal{T}E}\right. :V\mathcal{T}E\rightarrow VTE,
\]
is an isomorphism. If $f\in C^\infty (M)$ we will denote by $f^c$ and $f^v$
the \textit{complete and vertical lift} to $E$ of $f$ defined by
\[
f^c(u)=\sigma (u)(f),\quad f^v(u)=f(\pi (u)),\quad u\in E.
\]
For $s\in \Gamma (E)$ we can consider the \textit{vertical lift} of $s$
given by
\[
s^v(u)=s(\pi (u))_u^v,
\]
for $u\in E,$ where
\[
_u^v:E_{\pi (u)}\rightarrow T_u(E_{\pi (u)})
\]
is the canonical isomorphism.

There exists an unique vector field $s^c$ on $E$, the \textit{complete lift}
of $s$ satisfying the two following conditions:

i) $s^c$ is $\pi $-projectable on $\sigma (s),$

ii) $s^c(\stackrel{\wedge }{\alpha })=\widehat{\mathcal{L}_s\alpha
},$\\for all $\alpha \in \Gamma (E^{*}),$ where $\stackrel{\wedge
}{\alpha } (u)=\alpha (\pi (u))(u)$, $u\in E$ (see \cite{Gu1,
Gu2}). \\Considering the prolongation $\mathcal{T}E$ of $E$ over
the projection $\pi $, we may introduce the \textit{vertical lift
}$s^{\mathrm{v}}$ and the \textit{ complete lift} $s^{\mathrm{c}}$
of a section $s\in \Gamma (E)$ as the sections of
$\mathcal{T}E\rightarrow E$ given by (see \cite{Ma2})
\[
s^{\mathrm{v}}(u)=(0,s^v(u)),\quad s^{\mathrm{c}}(u)=(s(\pi
(u)),s^c(u)),\quad u\in E.
\]
Another two canonical objects on $\mathcal{T}E$ are the
\textit{Euler section } $C$ and the \textit{almost tangent
structure} (\textit{vertical endomorphism)} $J$.

\begin{definition}
The Euler section $C$ is the section of $\mathcal{T}E\rightarrow E$ defined
by
\[
C(u)=(0,u_u^v),\ \forall u\in E.
\]
\end{definition}

\begin{definition}
The vertical endomorphism $J$ is the section of the bundle
$(\mathcal{T} E)\oplus (\mathcal{T}E)^{*}\rightarrow E$
characterized by
\[
J(s^{\mathrm{v}})=0,\quad J\left( s^{\mathrm{c}}\right)
=s^{\mathrm{v} },\quad s\in \Gamma (E).
\]
\end{definition}

The vertical endomorphism satisfies
\[
J^2=0,\ ImJ=\ker J=V\mathcal{T}E,\quad \ [C,J]_{\mathcal{T}E}=-J.
\]

\begin{definition}
A section $\mathcal{S}$ of $\mathcal{T}E\rightarrow E$ is called
\textit{ semispray} (or \textit{second order differential equation
-SODE}) on $E$ if
\begin{equation}
J(\mathcal{S})=C.  \tag{1.3.1}
\end{equation}
\end{definition}

The local basis of $\Gamma (\mathcal{T}E)$ is given by $\{\mathcal{X}_\alpha
,\mathcal{V}_\alpha \}$, where \cite{Ma2}
\begin{equation}
\mathcal{X}_\alpha (u)=\left( s_\alpha (\pi (u)),\left. \sigma _\alpha
^i\frac \partial {\partial x^i}\right| _u\right) ,\quad \mathcal{V}_\alpha
(u)=\left( 0,\left. \frac \partial {\partial y^\alpha }\right| _u\right) ,
\tag{1.3.2}
\end{equation}
and $(\partial /\partial x^i,\partial /\partial y^\alpha )$ is the local
basis on $TE.$ The structure functions of $\mathcal{T}E$ are given by the
following formulas
\begin{equation}
\sigma ^1(\mathcal{X}_\alpha )=\sigma _\alpha ^i\frac \partial {\partial
x^i},\quad \sigma ^1(\mathcal{V}_\alpha )=\frac \partial {\partial y^\alpha
},  \tag{1.3.3}
\end{equation}
\begin{equation}
\lbrack \mathcal{X}_\alpha ,\mathcal{X}_\beta
]_{\mathcal{T}E}=L_{\alpha \beta }^\gamma \mathcal{X}_\gamma
,\quad [\mathcal{X}_\alpha ,\mathcal{V} _\beta
]_{\mathcal{T}E}=0,\quad [\mathcal{V}_\alpha ,\mathcal{V}_\beta
]_{ \mathcal{T}E}=0.  \tag{1.3.4}
\end{equation}
If $V$ is a section of $\mathcal{T}E,$ then in terms of basis
$\{\mathcal{X} _\alpha ,\mathcal{V}_\alpha \}$ it is
\[
V=Z^\alpha \mathcal{X}_\alpha +V^\alpha \mathcal{V}_\alpha ,
\]
and the vector field $\sigma ^1(V)$ $\in \chi (E)$ has the expression
\[
\sigma ^1(V)=\sigma _\alpha ^iZ^\alpha \frac \partial {\partial
x^i}+V^\alpha \frac \partial {\partial y^\alpha }.
\]
The vertical lift of a section $\rho =\rho ^\alpha s_\alpha $ and the
corresponding vector field are
\[
\rho ^{\mathrm{v}}=\rho ^\alpha \mathcal{V}_\alpha ,
\]
respectively
\[
\sigma ^1(\rho ^{\mathrm{v}})=\rho ^\alpha \frac \partial {\partial y^\alpha
}.
\]
$\,$ The coordinate expressions of Euler section $C$ and $\sigma ^1(C)$ are
\[
C=y^\alpha \mathcal{V}_\alpha ,\quad \sigma ^1(C)=y^\alpha \frac \partial
{\partial y^\alpha },
\]
and the local expression of $J$ is given by
\begin{equation}
J=\mathcal{X}^\alpha \otimes \mathcal{V}_\alpha ,  \tag{1.3.5}
\end{equation}
where $\{\mathcal{X}^\alpha ,\mathcal{V}^\alpha \}$ denotes the
corresponding dual basis of $\{\mathcal{X}_\alpha ,\mathcal{V}_\alpha \}$.

The Nijenhuis tensor of the vertical endomorphism
\[
\mathbf{N}_J((z,w)=[Jz,Jw]_E-J[Jz,w]_E-J[z,Jw]_E+J^2[z,w]_E
\]
vanishes, and it results that $J$ is a integrable structure.

The expression of the complete lift of a section $\rho =\rho ^\alpha
s_\alpha $ is
\[
\rho ^{\mathrm{c}}=\rho ^\alpha \mathcal{X}_\alpha +(\sigma
_\varepsilon ^i \frac{\partial \rho ^\alpha }{\partial
x^i}-L_{\beta \varepsilon }^\alpha \rho ^\beta )y^\varepsilon
\mathcal{V}_\alpha ,
\]
and therefore
\[
\sigma ^1(\rho ^{\mathrm{c}})=\rho ^\alpha \sigma _\alpha ^i\frac
\partial {\partial x^i}+\left( \sigma _\varepsilon
^i\frac{\partial \rho ^\alpha }{
\partial x^i}-L_{\beta \varepsilon }^\alpha \rho ^\beta \right)
y^\varepsilon \frac \partial {\partial y^\alpha }.
\]
In particular
\[
s_\alpha ^{\mathrm{v}}=\mathcal{V}_\alpha ,\quad s_\alpha
^{\mathrm{c}}= \mathcal{X}_\alpha -L_{\alpha \varepsilon }^\beta
y^\varepsilon \mathcal{V} _\beta .
\]
The local expression of the differential of a function $L$ on $\mathcal{T}E$
is
\[
d^EL=\sigma _\alpha ^i\frac{\partial L}{\partial
x^i}\mathcal{X}^\alpha + \frac{\partial L}{\partial y^\alpha
}\mathcal{V}^\alpha ,
\]
and therefore, we have
\[
d^Ex^i=\sigma _\alpha ^i\mathcal{X}^\alpha ,\ d^Ey^\alpha
=\mathcal{V} ^\alpha .
\]
The differential of sections of $(\mathcal{T}E)^{*}$ is determined by
\begin{equation}
d^E\mathcal{X}^\alpha =-\frac 12L_{\beta \gamma }^\alpha \mathcal{X}^\beta
\wedge \mathcal{X}^\gamma ,\quad d^E\mathcal{V}^\alpha =0.  \tag{1.3.6}
\end{equation}
In local coordinates a semispray has the expression
\begin{equation}
\mathcal{S}(x,y)=y^\alpha \mathcal{X}_\alpha +\mathcal{S}^\alpha
(x,y) \mathcal{V}_\alpha .  \tag{1.3.7}
\end{equation}
The integral curves of $\sigma ^1(\mathcal{S})$ satisfy the differential
equations
\begin{equation}
\frac{dx^i}{dt}=\sigma _\alpha ^i(x)y^\alpha ,\quad
\frac{dy^\alpha }{dt}= \mathcal{S}^\alpha (x,y).  \tag{1.3.8}
\end{equation}
If we have the relation
\[
\lbrack C,\mathcal{S}]_{\mathcal{T}E}=\mathcal{S},
\]
then $\mathcal{S}$ is called a spray and the functions $\mathcal{S}^\alpha $
are homogeneous functions of degree $2$ in $y^\alpha .$

Under a change of coordinates $x^{i^{\prime }}=x^{i^{\prime
}}(x^i),$ $ y^{\alpha ^{\prime }}=A_\alpha ^{\alpha ^{\prime
}}(x^i)y^\alpha $ on $E$ the transformation rule of the
coordinates on $\mathcal{T}E$ is given by
\[
\begin{array}{l}
x^{i^{\prime }}=x^{i^{\prime }}(x^i),\quad \\
y^{\alpha ^{\prime }}=A_\alpha ^{\alpha ^{\prime }}(x^i)y^\alpha , \\
u^{\alpha ^{\prime }}=A_\alpha ^{\alpha ^{\prime }}(x^i)u^\alpha , \\
\nu ^{\alpha ^{\prime }}=A_\alpha ^{\alpha ^{\prime }}\nu ^\alpha
+\sigma _\alpha ^iu^\alpha \frac{\partial A_\beta ^{\alpha
^{\prime }}}{\partial x^i} y^\beta .\quad
\end{array}
\]
and the rule of change of base is
\begin{eqnarray*}
\mathcal{X}_\alpha &=&A_\alpha ^{\alpha ^{\prime
}}\mathcal{X}_{\alpha ^{\prime }}+\sigma _\alpha ^i\frac{\partial
A_\beta ^{\alpha ^{\prime }}}{
\partial x^i}y^\beta \mathcal{V}_{\alpha ^{\prime }}, \\
\mathcal{V}_\alpha &=&A_\alpha ^{\alpha ^{\prime }}\mathcal{V}_{\alpha
^{\prime }}.
\end{eqnarray*}
The rule of change of dual base has the form
\begin{eqnarray*}
\mathcal{X}^\alpha &=&A_{\alpha ^{\prime }}^\alpha \mathcal{X}^{\alpha
^{\prime }}, \\
\mathcal{V}^\alpha &=&A_{\alpha ^{\prime }}^\alpha
\mathcal{V}^{\alpha ^{\prime }}+\sigma _{\alpha ^{\prime
}}^i\frac{\partial A_\beta ^\alpha }{
\partial x^i}y^\beta \mathcal{X}^{\alpha ^{\prime }}.
\end{eqnarray*}

\newpage\

\subsubsection{\textbf{Ehresmann nonlinear connections on Lie algebroid }$
\mathcal{T}E$}

As in the case $E=TM$ \cite{Cr, Gr} we can define the nonlinear connection.

\begin{definition}
An Ehresmann nonlinear connection (or connection) on $\mathcal{T}E$ is an
almost product structure $\mathcal{N}$ on $\pi _2:\mathcal{T}E\rightarrow E$
(i.e. a bundle morphism $\mathcal{N}:\mathcal{T}E\rightarrow \mathcal{T}E$,
such that $\mathcal{N}^2=Id$) smooth on $\mathcal{T}E\backslash \{0\}$ such
that
\[
V\mathcal{T}E=\ker (Id+\mathcal{N}).
\]
\end{definition}

If $\mathcal{N}$ is a connection on $\mathcal{T}E$ then $H\mathcal{T}E=\ker
(Id-\mathcal{N})$ is the horizontal subbundle associated to $\mathcal{N}$
and
\[
\mathcal{T}E=V\mathcal{T}E\oplus H\mathcal{T}E.
\]
Each $\rho \in \Gamma (\mathcal{T}E)$ can be written as $\rho
=\rho ^{ \mathrm{h}}+\rho ^{\mathrm{v}}$ where $\rho
^{\mathrm{h}}$, $\rho ^{\mathrm{v }}$ are sections in the
horizontal and respective vertical subbundles. If $ \rho
^{\mathrm{h}}=0,$ then $\rho $ is called\textit{\ vertical }and if
$ \rho ^{\mathrm{v}}=0,$ then $\rho $ is called
\textit{horizontal}. A connection $\mathcal{N}$ on $\mathcal{T}E$
induces two projectors $\mathrm{h}
,\mathrm{v}:\mathcal{T}E\rightarrow \mathcal{T}E$ such that
$\mathrm{h}(\rho )=\rho ^{\mathrm{h}}$ and $\mathrm{v}(\rho )=\rho
^{\mathrm{v}}$ for every $ \rho \in \Gamma (\mathcal{T}E)$. We
have
\begin{equation}
\mathrm{h}=\frac 12(Id+\mathcal{N}),\quad \mathrm{v}=\frac
12(Id-\mathcal{N} ),  \tag{1.3.9}
\end{equation}
\[
\ker \mathrm{h}=Im\mathrm{v}=V\mathcal{T}E,\quad Im\mathrm{h}=\ker
\mathrm{v} =H\mathcal{T}E.
\]
Locally, a connection can be expressed as
\begin{equation}
\mathcal{N}(\mathcal{X}_\alpha )=\mathcal{X}_\alpha
-2\mathcal{N}_\alpha ^\beta \mathcal{V}_\beta ,\quad
\mathcal{N}(\mathcal{V}_\beta )=-\mathcal{V} _\beta ,
\tag{1.3.10}
\end{equation}
where $\mathcal{N}_\alpha ^\beta =\mathcal{N}_\alpha ^\beta (x,y)$ are \ the
local coefficients of $\mathcal{N}$. The sections
\begin{equation}
\delta _\alpha =\mathrm{h}(\mathcal{X}_\alpha )=\mathcal{X}_\alpha
-\mathcal{ N}_\alpha ^\beta \mathcal{V}_\beta ,  \tag{1.3.11}
\end{equation}
generate a basis of $H\mathcal{T}E$. The frame $\{\delta _\alpha
,\mathcal{V} _\alpha \}$ is a local basis of $\mathcal{T}E$ called
\textit{adapted}. The dual adapted basis is $\{\mathcal{X}^\alpha
,\delta \mathcal{V}^\alpha \}$ where
\[
\ \delta \mathcal{V}^\alpha =\mathcal{V}^\alpha -\mathcal{N}_\beta ^\alpha
\mathcal{X}^\beta .
\]

\begin{proposition}
The Lie brackets of the adapted basis $\{\delta _\alpha ,\mathcal{V}_\alpha
\}$ are
\begin{equation}
[\delta _\alpha ,\delta _\beta ]_E=L_{\alpha \beta }^\gamma \delta
_\gamma + \mathcal{R}_{\alpha \beta }^\gamma \mathcal{V}_\gamma
,\quad [\delta _\alpha ,\mathcal{V}_\beta ]_E=\frac{\partial
\mathcal{N}_\alpha ^\gamma }{\partial y^\beta }\mathcal{V}_\gamma
,\quad [\mathcal{V}_\alpha ,\mathcal{V}_\beta ]_E=0,  \tag{1.3.12}
\end{equation}
where
\begin{equation}
\mathcal{R}_{\alpha \beta }^\gamma =\sigma _\beta ^i\frac{\partial
\mathcal{N }_\alpha ^\gamma }{\partial x^i}-\sigma _\alpha
^i\frac{\partial \mathcal{N} _\beta ^\gamma }{\partial
x^i}-\mathcal{N}_\beta ^\varepsilon \frac{\partial
\mathcal{N}_\alpha ^\gamma }{\partial y^\varepsilon
}+\mathcal{N}_\alpha ^\varepsilon \frac{\partial \mathcal{N}_\beta
^\gamma }{\partial y^\varepsilon }+L_{\alpha \beta }^\varepsilon
\mathcal{N}_\varepsilon ^\gamma .  \tag{1.3.13}
\end{equation}
\end{proposition}

\textbf{Proof}. Using (1.2.1) and (1.3.11) we get
\[
\lbrack \delta _\alpha ,\delta _\beta ]_{\mathcal{T}E}=\left(
\sigma _\beta ^i\frac{\partial \mathcal{N}_\alpha ^\varepsilon
}{\partial x^i}-\mathcal{N} _\beta ^\gamma \frac{\partial
\mathcal{N}_\alpha ^\varepsilon }{\partial y^\gamma }-\sigma
_\alpha ^i\frac{\partial \mathcal{N}_\beta ^\varepsilon }{
\partial x^i}+\mathcal{N}_\alpha ^\gamma \frac{\partial \mathcal{N}_\beta
^\varepsilon }{\partial y^\gamma }\right) \mathcal{V}_\varepsilon +L_{\alpha
\beta }^\gamma \mathcal{X}_\gamma .
\]
If we insert $\mathcal{X}_\gamma =\delta _\gamma +\mathcal{N}_\gamma
^\varepsilon \mathcal{V}_\varepsilon $ then the first relation from (1.3.12)
is obtained. By direct computation the second relation is verified.\hfill
\hbox{\rlap{$\sqcap$}$\sqcup$}\\We recall that the Nijenhuis tensor of an
endomorphism $A$ is given by
\[
\mathbf{N}_A(z,w)=[Az,Aw]_{\mathcal{T}E}-A[Az,w]_{\mathcal{T}E}-A[z,Aw]_{
\mathcal{T}E}+A^2[z,w]_{\mathcal{T}E}.
\]

\begin{definition}
The curvature of the connection $\mathcal{N}$ on $\mathcal{T}E$ is
given by $ \Omega =-\mathbf{N}_{\mathrm{h}}$ where $\mathrm{h}$ is
the horizontal projector and $\mathbf{N}_{\mathrm{h}}$ is the
Nijenhuis tensor of $\mathrm{h }$.
\end{definition}

\begin{proposition}
In local coordinates we have
\[
\Omega =-\frac 12\mathcal{R}_{\alpha \beta }^\gamma \mathcal{X}^\alpha
\wedge \mathcal{X}^\beta \otimes \mathcal{V}_\gamma ,
\]
where $\mathcal{R}_{\alpha \beta }^\gamma $ are given by (1.3.13)
and represent the local coordinate functions of the curvature
tensor $\Omega $ in the frame $\bigwedge^2\mathcal{T}E^{*}\otimes
\mathcal{T}E$ induced by $\{ \mathcal{X}_\alpha
,\mathcal{V}_\alpha \}$.
\end{proposition}

\textbf{Proof}. Since $\mathrm{h}^2=\mathrm{h}$ we obtain
\[
\Omega
(z,w)=-[\mathrm{h}z,\mathrm{h}w]_{\mathcal{T}E}+\mathrm{h}[\mathrm{h}
z,w]_{\mathcal{T}E}+\mathrm{h}[z,\mathrm{h}w]_{\mathcal{T}E}-\mathrm{h}
[z,w]_{\mathcal{T}E},
\]
\[
\Omega
(\mathrm{h}z,\mathrm{h}w)=-\mathrm{v}[\mathrm{h}z,\mathrm{h}w]_{
\mathcal{T}E},\quad \Omega (\mathrm{h}z,\mathrm{v}w)=\Omega
(\mathrm{v}z, \mathrm{v}w)=0,
\]
and in local coordinates we get
\[
\Omega (\delta _\alpha ,\delta _\beta )=-\mathrm{v}[\delta _\alpha
,\delta _\beta ]_{\mathcal{T}E}=-\mathcal{R}_{\alpha \beta
}^\gamma \mathcal{V} _\gamma ,
\]
which ends the proof.\hfill
\hbox{\rlap{$\sqcap$}$\sqcup$}

The curvature of the nonlinear connection is an obstruction to the
integrability of $H\mathcal{T}E$, understanding that a vanishing curvature
entails that horizontal sections are closed under the Lie algebroid bracket
of $\mathcal{T}E$.

\begin{remark}
{}$H\mathcal{T}E$ is integrable if and only if the curvature
$\Omega =- \mathbf{N}_{\mathrm{h}}$ of the nonlinear connection
vanishes.
\end{remark}

Let $\Psi $ a morphism of vector bundles $E$ and $\overline{E}$.
We recall that the connections $N$ on $E$ and $\overline{N}$ on
$\overline{E}$ are $ \Psi $-related if
\[
\Psi \circ N=\overline{N}\circ \Psi .
\]
We consider the connections $\mathcal{N}$ on $\mathcal{T}E$ and \
$\Bbb{N}$ on $TE$ which are $\sigma ^1$-related and a connection
$N$ on $T^2M$ which is $\sigma _{*}$-related with $\Bbb{N}$ on
$TE$ and $\widetilde{\sigma }$ -related with $\mathcal{N}$ on
$\mathcal{T}E$, where $\widetilde{\sigma }:
\mathcal{T}E\rightarrow T^2M$ is given by $\widetilde{\sigma
}=\sigma _{*}\circ \sigma ^1$ and $\sigma _{*}:TE\rightarrow T^2M$
is the tangent application of $\sigma $. It follows
\begin{equation}
\Bbb{N}\circ \sigma ^1=\sigma ^1\circ \mathcal{N},\quad N\circ
\sigma _{*}=\sigma _{*}\circ \Bbb{N},,\quad N\circ
\widetilde{\sigma }=\widetilde{ \sigma }\circ \mathcal{N},
\tag{1.3.14}
\end{equation}
Let us consider the adapted basis $(\stackrel{E}{\delta }_i,\frac \partial
{\partial y^\beta })$ of \ $\Bbb{N}$ and $(\stackrel{TM}{\delta }_i,\frac
\partial {\partial y^j})$ of $N$ given by
\[
\stackrel{E}{\delta }_i=\frac \partial {\partial x^i}-\Bbb{N}_i^\beta \frac
\partial {\partial y^\beta },
\]
and
\[
\stackrel{TM}{\delta }_i=\frac \partial {\partial x^i}-N_i^j\frac \partial
{\partial y^j}.
\]
Therefore, we get
\[
\sigma _{*}\left( \frac \partial {\partial x^i}\right) =\frac \partial
{\partial x^i}+\frac{\partial \sigma ^k}{\partial x^i}\frac \partial
{\partial y^k},\quad \sigma _{*}\left( \frac \partial {\partial y^\alpha
}\right) =\sigma _\alpha ^i\frac \partial {\partial y^i}.
\]

\begin{theorem}
The following relations hold
\[
\sigma ^1(\delta _\alpha )=\sigma _\alpha ^i\stackrel{E}{\delta }_i,\quad
\mathcal{N}_\alpha ^\beta =\sigma _\alpha ^i\Bbb{N}_i^\beta ,
\]
\begin{equation}
\sigma _{*}(\stackrel{E}{\delta }_i)=\stackrel{TM}{\delta
}_i,\quad \frac{
\partial \sigma ^j}{\partial x^i}+N_i^j=\Bbb{N}_i^\beta \sigma _\beta ^j,
\tag{1.3.15}
\end{equation}
\[
\widetilde{\sigma }(\delta _\alpha )=\sigma _\alpha
^i\stackrel{TM}{\delta } _i,\quad \sigma _\alpha ^i\frac{\partial
\sigma ^j}{\partial x^i}+\sigma _\alpha ^iN_i^j=\mathcal{N}_\alpha
^\beta \sigma _\beta ^j.
\]
\end{theorem}

\textbf{Proof.} The first relation from (1.3.14) leads to the
relation $\Bbb{ N}(\sigma ^1(\delta _\alpha ))=\sigma ^1(\delta
_\alpha )$ from which we get $\sigma ^1(\delta _\alpha )=\sigma
_\alpha ^i\stackrel{E}{\delta }_i$ and $ \mathcal{N}_\alpha ^\beta
=\sigma _\alpha ^i\Bbb{N}_i^\beta $. In the similar way the others
relations are obtained.\hfill \hbox{\rlap{$\sqcap$}$\sqcup$}

\begin{proposition}
For the curvature tensors of $\sigma ^1$-related connections $\mathcal{N}$
and $\Bbb{N}$ we have the relation
\begin{equation}
\mathcal{R}_{\alpha \beta }^\gamma =\sigma _\alpha ^i\sigma _\beta
^j\Bbb{R} _{ij}^\gamma ,  \tag{1.3.16}
\end{equation}
where
\[
\Bbb{R}_{ij}^\gamma =\stackrel{E}{\delta }_i(\Bbb{N}_j^\gamma
)-\stackrel{E}{ \delta }_j(\Bbb{N}_i^\gamma ),
\]
is the curvature tensor of the nonlinear connection on $TE.$
\end{proposition}

\textbf{Proof}. Using the relation $\mathcal{N}_\alpha ^\varepsilon =\sigma
_\alpha ^i\Bbb{N}_i^\varepsilon $ we obtain
\[
\mathcal{R}_{\alpha \beta }^\gamma =\sigma _\beta ^j\sigma _\alpha
^i\left( \stackrel{E}{\delta }_j\left( \Bbb{N}_i^\gamma \right)
-\stackrel{E}{\delta } _i\left( \Bbb{N}_j^\gamma \right) \right)
+\Bbb{N}_j^\gamma \left( \sigma _\beta ^i\frac{\partial \sigma
_\alpha ^j}{\partial x^i}-\sigma _\alpha ^i \frac{\partial \sigma
_\beta ^j}{\partial x^i}\right) +L_{\alpha \beta }^\varepsilon
\mathcal{N}_\varepsilon ^\gamma ,
\]
and from structure equations of the Lie algebroid (1.2.10), the
second term is $\Bbb{N}_j^\gamma \sigma _\varepsilon ^jL_{\beta
\alpha }^\varepsilon =- \mathcal{N}_\varepsilon ^\gamma L_{\alpha
\beta }^\varepsilon ,$ which concludes the proof. \hfill
\hbox{\rlap{$\sqcap$}$\sqcup$}

\begin{remark}
A $\sigma ^1$--related connection $\Bbb{N}$ on $TE$ determines a
connection $ \mathcal{N}$ on $\mathcal{T}E$ with the coefficients
$\mathcal{N}_\alpha ^\beta =\sigma _\alpha ^i\Bbb{N}_i^\beta $ and
the curvature $\mathcal{R}
_{\alpha \beta }^\gamma =\sigma _\alpha ^i\sigma _\beta ^j\Bbb{R}%
_{ij}^\gamma .$ The converse is not true because $\sigma $ is only injective.
\end{remark}

Let $J$ be the vertical endomorphism.

\begin{remark}
{}Let $\mathcal{N}$ be a bundle morphism of $\pi _2:\mathcal{T}E\rightarrow
E $, smooth on $\mathcal{T}E\backslash \{0\}$. Then $\mathcal{N}$ is a
connection on $\mathcal{T}E$ if and only if
\[
J\mathcal{N}=J,\quad \mathcal{N}J=-J.
\]
\end{remark}

The proof proceeds as in the case $E=TM$ and will be omitted.

\begin{definition}
The torsion of a nonlinear connection $\mathcal{N}$ is the vector valued two
form $t=[J,\mathrm{h}]$ where $\ \mathrm{h}$ is the horizontal projector and
$[\cdot ,\cdot ]$ is the Fr\"olicher-Nijenhuis bracket.
\end{definition}

\begin{definition}
$t$ is a semibasic vector-valued form. Its local expression is
\begin{equation}
t=\frac 12t_{\alpha \beta }^\gamma \mathcal{X}^\alpha \wedge
\mathcal{X} ^\beta \otimes \mathcal{V}_\gamma ,  \tag{1.3.17}
\end{equation}
where
\begin{equation}
t_{\alpha \beta }^\gamma =\frac{\partial \mathcal{N}_\alpha
^\gamma }{
\partial y^\beta }-\frac{\partial \mathcal{N}_\beta ^\gamma }{\partial
y^\alpha }-L_{\alpha \beta }^\gamma .  \tag{1.3.18}
\end{equation}
\end{definition}

\textbf{Proof}. We have
\begin{eqnarray*}
\lbrack J,\mathrm{h}](z,w)
&=&[Jz,\mathrm{h}w]_{\mathcal{T}E}+[\mathrm{h}
z,Jw]_{\mathcal{T}E}+J[z,w]_{\mathcal{T}E}-J[z,\mathrm{h}w]_{\mathcal{T}E}-
\\
&&\ \ \ \ \ \ \ -\
J[\mathrm{h}z,w]_{\mathcal{T}E}-\mathrm{h}[z,Jw]_{
\mathcal{T}E}-\mathrm{h}[Jz,w]_{\mathcal{T}E},
\end{eqnarray*}
and in local coordinates we get
\[
t(\mathcal{X}_\alpha ,\mathcal{X}_\beta )=\left( \frac{\partial
\mathcal{N} _\alpha ^\gamma }{\partial y^\beta }-\frac{\partial
\mathcal{N}_\beta ^\gamma }{\partial y^\alpha }-L_{\alpha \beta
}^\gamma \right) \mathcal{V} _\gamma ,\quad t(\mathcal{X}_\alpha
,\mathcal{V}_\beta )=t(\mathcal{V} _\alpha ,\mathcal{V}_\beta )=0.
\]
\hfill \hbox{\rlap{$\sqcap$}$\sqcup$}\\Now, let us consider the
linear mapping $ \Bbb{F}:\mathcal{T}E\rightarrow \mathcal{T}E$,
defined by
\begin{equation}
\Bbb{F}(\mathrm{h}z)=-\mathrm{v}z,\quad
\Bbb{F}(\mathrm{v}z)=\mathrm{h}z, \text{ }  \tag{1.3.19}
\end{equation}
for $z\in \Gamma (\mathcal{T}E)$ and $\mathrm{h}$, \textrm{v} the horizontal
and vertical projectors of the nonlinear connection on $\mathcal{T}E$.

\begin{proposition}
The mapping $\Bbb{F}$ has the properties:

i) $\Bbb{F}$ is globally defined on $\mathcal{T}E,$

ii) Locally, it is given by
\begin{equation}
\Bbb{F}=-\mathcal{V}_\alpha \otimes \mathcal{X}^\alpha +\delta _\alpha
\otimes \delta \mathcal{V}^\alpha ,  \tag{1.3.20}
\end{equation}

iii) $\Bbb{F}$ is an almost complex structure $\Bbb{F}\circ \Bbb{F}=-Id.$
\end{proposition}

\textbf{Proof}. It results by definition that $\Bbb{F}$ is globally defined
and
\[
\left( \Bbb{F}\circ \Bbb{F}\right) (\mathrm{h}z)=\Bbb{F}\left(
-\mathrm{v} z\right) =-\mathrm{h}z,\quad \left( \Bbb{F}\circ
\Bbb{F}\right) (\mathrm{v} z)=\Bbb{F}\left( \mathrm{h}z\right)
=-\mathrm{v}z.
\]
In local coordinates we get $\Bbb{F}\left( \delta _\alpha \right)
=-\mathcal{ V}_\alpha $ and $\Bbb{F}(\mathcal{V}_\alpha )=\delta
_\alpha $ which ends the proof.\hfill
\hbox{\rlap{$\sqcap$}$\sqcup$}

\begin{proposition}
The almost complex structure is integrable if and only if the nonlinear
connection is locally flat, that is the curvature and torsion vanish.
\end{proposition}

\textbf{Proof}. Let $\mathbf{N}_{\Bbb{F}}$ be the Nijenhuis tensor of the
almost complex structure. We find
\begin{equation}
\begin{array}{l}
\mathbf{N}_{\Bbb{F}}(\delta _\alpha ,\delta _\beta )=t_{\alpha
\beta }^\gamma \delta _\gamma -\mathcal{R}_{\alpha \beta }^\gamma
\mathcal{V}
_\gamma , \\
\\
\mathbf{N}_{\Bbb{F}}(\delta _\alpha ,\mathcal{V}_\beta
)=-\mathcal{R} _{\alpha \beta }^\gamma \delta _\gamma -t_{\alpha
\beta }^\gamma \mathcal{V}
_\gamma , \\
\\
\mathbf{N}_{\Bbb{F}}(\mathcal{V}_\alpha ,\mathcal{V}_\beta
)=-\mathbf{N}_{ \Bbb{F}}(\delta _\alpha ,\delta _\beta ).
\end{array}
\tag{1.3.21}
\end{equation}
From (1.3.21) one reads immediately that $\mathbf{N}_{\Bbb{F}}=0$ if and
only if $t=0$ and $\Omega =0.$\hfill
\hbox{\rlap{$\sqcap$}$\sqcup$}

A curve $u:[t_0,t_1]\rightarrow E$ is called \textit{admissible} if $\sigma
(u(t))=\dot c(t)$ where $c(t)=\pi (u(t))$ is the base curve. A nonlinear
connection on $\mathcal{T}E$ induces a covariant derivative of the sections
defined locally as follows
\[
\mathcal{D}_\rho \eta =\rho ^\alpha \left( \sigma _\alpha ^i\frac{\partial
\eta ^\beta }{\partial x^i}+\mathcal{N}_\alpha ^\beta \right) s_\beta ,
\]
where $\rho =\rho ^\alpha s_\alpha $ and $\eta =\eta ^\alpha s_\alpha $. The
derivative is linear in the first argument and it respects multiplication of
second argument by real numbers, but not necessarily sum, except the case
when the coefficients $\mathcal{N}_\alpha ^\beta $ are the local
coefficients of a linear connection. The linearity in the first argument
permits us to define the derivative of a section $\eta \in \Gamma (E)$ with
respect to $a\in E_u$ by setting
\[
\mathcal{D}_a\eta =(\mathcal{D}_\rho \eta )(u),
\]
where $\rho \in \Gamma (E)$ is satisfying $\rho (u)=a$. Also, the
covariant derivative allows us to take the derivative of sections
along curves. If we have a morphism of Lie algebroids $\Phi
:F\rightarrow E$ over the map $ \varphi :N\rightarrow M$ and a
section $\eta :N\rightarrow E$ along $\varphi ,$ i.e $\eta (n)\in
E_{\varphi (n)}$, $n\in N$, then $\eta $ can be written in the
form
\[
\eta =\sum_{l=1}^pF_l(\xi _l\circ \varphi ),
\]
for some sections $\{\xi _1,...,\xi _p\}$ of $E$ and some
functions $ \{F_1,...,F_p)\in C^\infty (N)$ and the derivative of
$\eta $ along $\varphi $ is given by
\[
\mathcal{D}_b\eta =\sum_{l=1}^p[(\sigma _F(b)F_l)\xi _l(\varphi
(n))+F_l(n) \mathcal{D}_{\Phi (b)}\xi _l],\quad b\in F_n,
\]
where $\sigma _F$ is the anchor map of the Lie algebroid $F\rightarrow N$
(see \cite{Co}).

Let $a:I\rightarrow E$ be an admissible curve and let
$b:I\rightarrow E$ be a curve in $E$, both of them projecting by
$\pi $ onto the same curve $ \gamma $ in $M$, $\pi (a(t))=\pi
(b(t))=\gamma (t).$ Take the particular case of Lie algebroid
structure $TI\rightarrow I$ and the morphism $\Phi :TI\rightarrow
E$, $\Phi (t,\stackrel{\cdot }{t})=\stackrel{\cdot }{t}\gamma (t)$
over $\gamma :I\rightarrow M$. Then one can define the derivative
of $ b(t)$ along $a(t)$ as $\mathcal{D}_{d/dt}b(t)$. In local
coordinates, we obtain
\[
\mathcal{D}_{a(t)}b(t)=\left( \frac{db^\beta }{dt}+\mathcal{N}_\alpha ^\beta
a^\alpha \right) s_\beta (\gamma (t)).
\]

\begin{definition}
An admissible curve $c(t)$ is a path (autoparallel) for nonlinear connection
$\mathcal{N}$ if and only if
\[
\mathcal{D}_{c(t)}c(t)=0.
\]
\end{definition}

In local coordinates we get
\[
\frac{dc^\beta }{dt}+\mathcal{N}_\alpha ^\beta (x,y)c^\alpha =0.
\]

From the previous considerations we have:

\begin{proposition}
An admissible curve $c(t)$ in $E$ is autoparallel for the nonlinear
connection if and only if
\begin{equation}
\frac{dx^i}{dt}=\sigma _\alpha ^iy^\alpha ,\quad \frac{dy^\beta
}{dt}+ \mathcal{N}_\alpha ^\beta y^\alpha =0,  \tag{1.3.22}
\end{equation}
where $x^i=x^i(t)=x^i(c(t))$, $y^\alpha =y^\alpha (t)=y^\alpha
(c(t))$, $ \sigma _\alpha ^i=\sigma _\alpha ^i(t)=\sigma _\alpha
^i(c(t)).$
\end{proposition}

Let $\mathcal{N}$ be a nonlinear connection on $\mathcal{T}E$,
$\mathcal{S} ^{\prime }$ an arbitrary semispray on $\mathcal{T}E$
and $\mathrm{h}$ the horizontal projector of $\mathcal{N}$. We
consider $\mathcal{S}=\mathrm{h} \mathcal{S}^{\prime }$ and for
any other semispray $\mathcal{S}^{\prime \prime }$ on
$\mathcal{T}E$ we have $\mathrm{h}(\mathcal{S}^{\prime }-
\mathcal{S}^{\prime \prime })=\mathrm{h}((\mathcal{S}^{\prime
\alpha }- \mathcal{S}^{\prime \prime \alpha })\mathcal{V}_\alpha
)=0$ and it results that $\mathcal{S}$ does not depend on the
choose of $\mathcal{S}^{\prime }$. We have
\[
J\mathcal{S}=J\mathrm{h}\mathcal{S}^{\prime }=J\mathcal{S}^{\prime }=C,
\]
so $\mathcal{S}$ is a semispray, which is called the \textit{associated
semispray} to $\mathcal{N}$.

\begin{proposition}
A nonlinear connection $\mathcal{N}$ and its associated semispray have the
same paths.
\end{proposition}

\textbf{Proof}. For the arbitrary semispray $\mathcal{S}^{\prime }=y^\alpha
\mathcal{X}_\alpha +\mathcal{S}^{\prime \alpha }\mathcal{V}_\alpha ,$ the
associated semispray of $\mathcal{N}$ is
\[
\mathcal{S}=\mathrm{h}\mathcal{S}^{\prime }=y^\alpha
\mathcal{X}_\alpha - \mathcal{N}_\alpha ^\beta y^\alpha
\mathcal{V}_\beta ,
\]
so
\[
\mathcal{S}^\beta =-\mathcal{N}_\alpha ^\beta y^\alpha .
\]
From (1.3.8) and (1.3.22) it results the conclusion.\hfill
\hbox{\rlap{$\sqcap$}$\sqcup$}

\begin{remark}
If $\mathcal{S}$ is a semispray on $\mathcal{T}E$, then we have
\begin{equation}
J[\mathcal{S},Jz]_{\mathcal{T}E}=-Jz,\quad z\in \Gamma (\mathcal{T}E).
\tag{1.3.23}
\end{equation}
\end{remark}

\begin{theorem}
{}Let $J$ be the vertical endomorphism on $\mathcal{T}E.$ If $\xi $ is a
semispray then
\begin{equation}
\mathcal{N}=-\mathcal{L}_{\mathcal{S}}J,  \tag{1.3.24}
\end{equation}
is a connection on $\mathcal{T}E.$
\end{theorem}

\textbf{Proof}. Since
\[
\mathcal{N(}\upsilon
)=-\mathcal{L}_{\mathcal{S}}J\mathcal{(}\upsilon )=-[
\mathcal{S},J\upsilon ]_{\mathcal{T}E}+J[\mathcal{S},\upsilon
]_{\mathcal{T} E}
\]
using (1.3.23) we get
\[
J\mathcal{N(}\upsilon )=-J[\mathcal{S},J\upsilon
]_{\mathcal{T}E}+J^2[ \mathcal{S},\upsilon
]_{\mathcal{T}E}=J\upsilon ,
\]
and
\[
\mathcal{N}J(\upsilon )=-[\mathcal{S},J^2\upsilon
]_{\mathcal{T}E}+J[ \mathcal{S},J\upsilon
]_{\mathcal{T}E}=-J\upsilon .
\]
By using the Remark 1.3.3 we get the proof of the theorem.\hfill
\hbox{\rlap{$\sqcap$}$\sqcup$}

\begin{remark}
The connection $\mathcal{N}=-\mathcal{L}_{\mathcal{S}}J$ is induced by the
semispray $\mathcal{S}$. Its local \ coefficients are given by
\begin{equation}
\mathcal{N}_\alpha ^\beta =\frac 12\left( -\frac{\partial \mathcal{S}^\beta
}{\partial y^\alpha }+y^\varepsilon L_{\alpha \varepsilon }^\beta \right) .
\tag{1.3.25}
\end{equation}
\end{remark}

\textbf{Proof}. By direct computation it results
\[
\mathcal{N}(\mathcal{X}_\alpha
)=-[\mathcal{S},J(\mathcal{X}_\alpha )]_{
\mathcal{T}E}+J[\mathcal{S},\mathcal{X}_\alpha ]_{\mathcal{T}E}
\]
\[
=\mathcal{X}_\alpha +\frac{\partial \mathcal{S}^\beta }{\partial
y^\alpha } \mathcal{V}_\beta +J(y^\beta L_{\beta \alpha }^\gamma
\mathcal{X}_\gamma -\sigma _\alpha ^i\frac{\partial
\mathcal{S}^\beta }{\partial x^i}\mathcal{V} _\beta )
\]
\[
=\mathcal{X}_\alpha +\left( \frac{\partial \mathcal{S}^\beta }{\partial
y^\alpha }+y^\gamma L_{\gamma \alpha }^\beta \right) \mathcal{V}_\beta ,
\]
and using (1.3.10) we obtain (1.3.25).\hfill
\hbox{\rlap{$\sqcap$}$\sqcup$}

\begin{proposition}
The torsion of the connection $\mathcal{N}=-\mathcal{L}_{\mathcal{S}}J$
vanishes.
\end{proposition}

\textbf{Proof}. We have
\[
t=[J,h]_{\mathcal{T}E}=\frac 12\left(
[J,Id]_{\mathcal{T}E}+[J,-[\xi ,J]_{
\mathcal{T}E}]_{\mathcal{T}E}\right) =\frac 12[J,[J,\xi
]_{\mathcal{T}E}]_{ \mathcal{T}E}.
\]
Using Jacobi identity we obtain that $t=0.$ Also, if we use (1.3.25) into
(1.3.18), by direct computation, the same result is obtained. \hfill
\hbox{\rlap{$\sqcap$}$\sqcup$}

\begin{proposition}
The associated semispray of $\mathcal{N}=-\mathcal{L}_{\mathcal{S}}J$ is
given by
\[
\frac 12(\mathcal{S}-[\mathcal{S},\mathcal{C}]_{\mathcal{T}E}).
\]
\end{proposition}

\textbf{Proof}. The associated semispray is
\[
\mathrm{h}\mathcal{S}=\frac 12\mathcal{S}+\frac
12\mathcal{N}(\mathcal{S} )=\frac
12(\mathcal{S}-[\mathcal{S},J\mathcal{S}]_{\mathcal{T}E}+J[\mathcal{S
},\mathcal{S}]_{\mathcal{T}E})=\frac
12(\mathcal{S}-[\mathcal{S},\mathcal{C} ]_{\mathcal{T}E}).
\]
\hfill
\hbox{\rlap{$\sqcap$}$\sqcup$}\newpage\

\subsubsection{\textbf{Lagrangian formalism on Lie algebroids}}

We consider the Cartan $1$-section
\[
\theta _L=J(d^EL),
\]
which, in local coordinates is
\begin{equation}
\theta _L=\frac{\partial L}{\partial y^\alpha }\mathcal{X}^\alpha ,
\tag{1.3.26}
\end{equation}
The differential of $\theta _L$ is the Cartan $2$-section
\[
\omega _L=d^E\theta _L,
\]
and from the local coordinate expression of $\theta _L$ we get
\[
\omega _L=d^E\left( \frac{\partial L}{\partial y^\alpha }\right)
\wedge \mathcal{X}^\alpha +\frac{\partial L}{\partial y^\alpha
}\wedge d^E\mathcal{X }^\alpha .
\]
But
\[
d^E\mathcal{X}^\alpha =-\frac 12L_{\beta \gamma }^\alpha \mathcal{X}^\beta
\wedge \mathcal{X}^\gamma ,
\]
and it results \cite{Ma2}
\begin{equation}
\omega _L=\frac{\partial L}{\partial y^\alpha \partial y^\beta
}\mathcal{V} ^\beta \wedge \mathcal{X}^\alpha +\frac 12\left(
\frac{\partial ^2L}{
\partial x^i\partial y^\beta }\sigma _\alpha ^i-\frac{\partial ^2L}{\partial
x^i\partial y^\alpha }\sigma _\beta ^i-\frac{\partial L}{\partial
y^\gamma } L_{\alpha \beta }^\gamma \right) \mathcal{X}^\alpha
\wedge \mathcal{X}^\beta .  \tag{1.3.27}
\end{equation}
The function $L$ is said to be a regular Lagrangian if $\omega _L$ is
regular at every point as a bilinear form. Let us consider the energy
function given by

\[
E_L\stackrel{def}{=}y^\alpha \frac{\partial L}{\partial y^\alpha }-L,
\]
and the symplectic equation
\begin{equation}
i_S\omega _L=-d^EE_L,\quad S\in \Gamma (\mathcal{T}E).  \tag{1.3.28}
\end{equation}
In local coordinates, considering the section
\[
\mathcal{S}=f^\alpha \mathcal{X}_\alpha +\mathcal{S}^\alpha
\mathcal{V} _\alpha ,
\]
we obtain the equations
\[
i_S\omega _L=\left( \mathcal{S}^\beta g_{\alpha \beta }+f^\beta
\left( \sigma _\beta ^i\frac{\partial ^2L}{\partial x^i\partial
y^\alpha }-\sigma _\alpha ^i\frac{\partial ^2L}{\partial
x^i\partial y^\beta }+\frac{\partial L }{\partial y^\gamma
}L_{\alpha \beta }^\gamma \right) \right) \mathcal{X} ^\alpha
-f^\beta g_{\alpha \beta }\mathcal{V}^\alpha ,
\]
\[
-dE_L=\left( \sigma _\alpha ^i\frac{\partial L}{\partial
x^i}-\sigma _\alpha ^i\frac{\partial ^2L}{\partial x^i\partial
y^\beta }y^\beta \right) \mathcal{ X}^\alpha -y^\beta g_{\alpha
\beta }\mathcal{V}^\alpha .
\]
The equality of the $\mathcal{V}^\alpha $ components yields
\[
g_{\alpha \beta }(y^\beta -f^\beta )=0,
\]
and the regularity of the Lagrangian implies $y^\beta =f^\beta $,
which means that $\mathcal{S}=y^\alpha \mathcal{X}_\alpha
+\mathcal{S}^\alpha \mathcal{V}_\alpha $ is a semispray.
Analogously, the equality of the $ \mathcal{X}^\alpha $ components
leads to the equation
\[
\mathcal{S}^\beta g_{\alpha \beta }+\sigma _\beta ^i\frac{\partial
^2L}{
\partial x^i\partial y^\alpha }y^\beta +\frac{\partial L}{\partial y^\gamma }
y^\beta L_{\alpha \beta }^\gamma =\sigma _\alpha ^i\frac{\partial
L}{
\partial x^i},
\]
and the regularity condition of the Lagrangian determines the components of
the semispray
\begin{equation}
\mathcal{S}^\varepsilon =g^{\varepsilon \beta }\left( \sigma
_\beta ^i\frac{
\partial L}{\partial x^i}-\sigma _\alpha ^i\frac{\partial ^2L}{\partial
x^i\partial y^\beta }y^\alpha -L_{\beta \alpha }^\theta y^\alpha
\frac{
\partial L}{\partial y^\theta }\right) ,  \tag{1.3.29}
\end{equation}
where $g_{\alpha \beta }g^{\beta \gamma }=\delta _\alpha ^\gamma $. From
(1.3.29) and (1.3.25) it results:

\begin{corollary}
For a regular Lagrangian $L$, there exists a nonlinear connection
$\mathcal{N }$ with the coefficients given by
\begin{equation}
\mathcal{N}_\alpha ^\varepsilon =\frac 12\left( -\frac{\partial
\mathcal{S} ^\varepsilon }{\partial y^\alpha }+y^\beta L_{\alpha
\beta }^\varepsilon \right) ,  \tag{1.3.30}
\end{equation}
where
\[
\mathcal{S}^\varepsilon =g^{\varepsilon \beta }\left( \sigma
_\beta ^i\frac{
\partial L}{\partial x^i}-\sigma _\alpha ^i\frac{\partial ^2L}{\partial
x^i\partial y^\beta }y^\alpha -L_{\beta \alpha }^\theta y^\alpha
\frac{
\partial L}{\partial y^\theta }\right) ,
\]
\end{corollary}

and will be called the \textit{canonical nonlinear connection} induced by a
regular Lagrangian $L$.

If $(x^i)$ are coordinates on $M$, $\{s_\alpha \}$ is a local
basis of $ \Gamma (E)$, $(x^i,y^\alpha )$ are the corresponding
coordinates on $E$ and $ \gamma (t)=(x^i(t),y^\alpha (t))$ then,
$\gamma $ is a solution of the Euler-Lagrange equations if and
only if \cite{We2}
\begin{equation}
\frac{dx^i}{dt}=\sigma _\alpha ^iy^\alpha ,\quad \frac d{dt}\left(
\frac{
\partial L}{\partial y^\alpha }\right) =\sigma _\alpha ^i\frac{\partial L}{
\partial x^i}-L_{\alpha \beta }^\theta y^\beta \frac{\partial L}{\partial
y^\theta }.  \tag{1.3.31}
\end{equation}

\newpage\

\subsubsection{\textbf{Homogeneous connections}}

\begin{definition}
The morphism
\[
\mathcal{T}=\frac 12\mathcal{L}_C\mathcal{N},
\]
is called the \textit{tension} of the nonlinear connection.
\end{definition}

In local coordinates we get
\[
\mathcal{T}(\mathcal{X}_\alpha )=\left( \mathcal{N}_\alpha ^\beta
-\frac{
\partial \mathcal{N}_\alpha ^\beta }{\partial y^\gamma }y^\gamma \right)
\mathcal{V}_\beta ,\quad \mathcal{T}(\mathcal{V}_\alpha )=0,
\]
and it results
\begin{equation}
\mathcal{T}=\left( \mathcal{N}_\alpha ^\beta -\frac{\partial
\mathcal{N} _\alpha ^\beta }{\partial y^\gamma }y^\gamma \right)
\mathcal{X}^\alpha \otimes \mathcal{V}_\beta .  \tag{1.3.32}
\end{equation}
It is obvious that $\mathcal{T}$ is vanishing, if and only if the nonlinear
connection is homogeneous of degree $1$ with respect to $y^\alpha $.

\begin{proposition}
If $\mathcal{S}$ is a spray then $\mathcal{N}=-\mathcal{L}_{\mathcal{S}}J$
is a homogeneous nonlinear connection.
\end{proposition}

\textbf{Proof}. Using (1.3.25) we get
\[
\mathcal{T}=\left( -\frac{\partial \mathcal{S}^\gamma }{\partial
y^\alpha } +y^\beta \frac{\partial ^2\mathcal{S}^\gamma }{\partial
y^\alpha \partial y^\beta }\right) \mathcal{X}^\alpha \otimes
\mathcal{V}_\gamma .
\]
But $\mathcal{S}$ is a spray and it results that $\mathcal{S}^\gamma $ is
homogeneous of degree $2$, there is
\[
2\mathcal{S}^\gamma =y^\beta \frac{\partial \mathcal{S}^\gamma }{\partial
y^\beta },
\]
and
\[
\frac{\partial \mathcal{S}^\gamma }{\partial y^\alpha }=y^\beta
\frac{
\partial ^2\mathcal{S}^\gamma }{\partial y^\alpha \partial y^\beta },
\]
therefore, the tension vanishes.\hfill
\hbox{\rlap{$\sqcap$}$\sqcup$}

\begin{definition}
The strong torsion $T$ of $\mathcal{N}$ is given by
\[
T=i_\xi t-\mathcal{T},
\]
where $\mathcal{T}$ is the tension, $t$ is the torsion of $\mathcal{N}$, and
$i_\xi $ is the contraction with $\xi $.
\end{definition}

Locally, we obtain
\[
T(\delta _\alpha )=\left( \frac{\partial \mathcal{N}_\beta ^\gamma
}{
\partial y^\alpha }y^\beta -\mathcal{N}_\alpha ^\gamma +y^\beta L_{\alpha
\beta }^\gamma \right) \mathcal{V}_\gamma ,\quad T(\mathcal{V}_\alpha )=0.
\]

\begin{proposition}
The strong torsion $T$ of a nonlinear $\mathcal{N}$ connection vanishes if
and only if the torsion and the tension of $\mathcal{N}$ vanish.
\end{proposition}

\textbf{Proof}. If $T=0$ then we have
\[
N_\alpha ^\gamma =\frac{\partial \mathcal{N}_\beta ^\gamma }{\partial
y^\alpha }y^\beta +y^\beta L_{\alpha \beta }^\gamma
\]
and it results
\begin{eqnarray*}
t(\delta _\alpha ,\delta _\beta ) &=&\frac{\partial
\mathcal{N}_\beta ^\gamma }{\partial y^\alpha }+\frac{\partial
^2\mathcal{N}_\varepsilon ^\gamma }{\partial y^\alpha \partial
y^\beta }y^\varepsilon -\frac{\partial \mathcal{N}_\alpha ^\gamma
}{\partial y^\beta }-\frac{\partial ^2\mathcal{N} _\varepsilon
^\gamma }{\partial y^\alpha \partial y^\beta }y^\varepsilon
+L_{\alpha \beta }^\gamma -L_{\beta \alpha }^\gamma -L_{\alpha
\beta }^\gamma
\\
&=&\frac{\partial \mathcal{N}_\beta ^\gamma }{\partial y^\alpha
}-\frac{
\partial \mathcal{N}_\alpha ^\gamma }{\partial y^\beta }-L_{\alpha \beta
}^\gamma =-t(\delta _\alpha ,\delta _\beta ).
\end{eqnarray*}
which yields $t(\delta _\alpha ,\delta _\beta )=0$ and $\mathcal{T}=0$.
\hfill
\hbox{\rlap{$\sqcap$}$\sqcup$}

\begin{definition}
A function $\mathcal{F}:E\rightarrow [0,\infty ]$ which satisfies
the following properties\\1) $\mathcal{F}$ is $C^\infty $ on
$E\backslash \{0\}$
\\2) $\mathcal{F}(\lambda u)=\lambda \mathcal{F}(u)$ for $\lambda >0$ and $
u\in E_x.$ $x\in M.$\\3) For each $y\in E_x\backslash \{0\}$ the quadratic
form
\[
g_{\alpha \beta }(x,y)=\frac 12\frac{\partial ^2\mathcal{F}^2}{\partial
y^\alpha \partial y^\beta },
\]
is positive definite, will be called the Finsler function on a Lie algebroid.
\end{definition}

If we insert $\mathcal{L}=\frac 12\mathcal{F}^2=\frac 12g_{\alpha \beta
}y^\alpha y^\beta $ into the expression of semispray (1.3.29) we obtain \cite
{Po11}

\begin{corollary}
A homogeneous nonlinear connections has the coefficients given by
\[
\mathcal{N}_\alpha ^\varepsilon =\frac 12\left( -\frac{\partial
\mathcal{S} ^\varepsilon }{\partial y^\alpha }+y^\beta L_{\alpha
\beta }^\varepsilon \right) ,
\]
with
\begin{equation}
\mathcal{S}^\delta =\frac 12g^{\delta \beta }\left( \sigma _\alpha
^i\frac{
\partial g_{\beta \gamma }}{\partial x^i}+\sigma _\gamma ^i\frac{\partial
g_{\alpha \beta }}{\partial x^i}-\sigma _\beta ^i\frac{\partial g_{\alpha
\gamma }}{\partial x^i}+g_{\varepsilon \alpha }L_{\beta \gamma }^\varepsilon
+g_{\varepsilon \gamma }L_{\beta \alpha }^\varepsilon -g_{\varepsilon \beta
}L_{\gamma \alpha }^\varepsilon \right) y^\alpha y^\gamma .  \tag{1.3.33}
\end{equation}
and is called the canonical nonlinear connection associated to a Finsler
function.
\end{corollary}

\begin{remark}
In the particular case of the standard Lie algebroid $E=TM$ and $\sigma =Id$
the Cartan nonlinear connection is obtained.
\end{remark}

We consider the canonical nonlinear connection and $\left\| y\right\|
^2=g_{\alpha \beta }y^\alpha y^\beta =\mathcal{F}^2$ is the square of the
norm of the \textit{Euler} section. The almost complex structure
characterized by (1.3.20) does not preserve the property of homogeneity of
the sections. Indeed, it applies the $1$-homogeneous section $\delta _\alpha
$ onto the $0$-homogeneous section $\mathcal{V}_\alpha $, $\alpha \in
\overline{1,m}$. We define a new almost complex structure $\Bbb{F}_0:%
\mathcal{T}E\rightarrow \mathcal{T}E$ given by
\[
\Bbb{F}_0(\delta _\alpha )=-\frac{\mathcal{F}}a\mathcal{V}_\alpha ,\quad
\Bbb{F}_0(\mathcal{V}_\alpha )=\frac a{\mathcal{F}}\delta _\alpha ,\quad
a>0.
\]
It is not difficult to prove that $\Bbb{F}_0^2=-Id$ and $\Bbb{F}_0$
preserves the property of the homogeneity of the sections.

\begin{theorem}
The almost complex structure $\Bbb{F}_0$ is integrable if and only if the
following relations hold
\begin{equation}
\begin{array}{l}
\mathcal{R}_{\alpha \beta }^\gamma =\frac 1{a^2}\left( y_\alpha \delta
_\beta ^\gamma -y_\beta \delta _\alpha ^\gamma \right) , \\
\\
\delta _\alpha (\mathcal{F}^2)\delta _\beta ^\gamma =\delta _\beta
(\mathcal{ F}^2)\delta _\alpha ^\gamma ,
\end{array}
\tag{1.3.34}
\end{equation}
where $y_\alpha =g_{\alpha \beta }y^\beta ,$ $\alpha ,\beta
,\gamma = \overline{1,m}.$
\end{theorem}

\textbf{Proof.} For the Nijenjuis tensor $\mathbf{N}_{\Bbb{F}_0}$ we have
\[
\begin{array}{l}
\mathbf{N}_{\Bbb{F}_0}(\delta _\alpha ,\delta _\beta )=\left(
t_{\alpha \beta }^\gamma +\frac 1{2\mathcal{F}^2}\left( \delta
_\beta (\mathcal{F} ^2)\delta _\alpha ^\gamma -\delta _\alpha
(\mathcal{F}^2)\delta _\beta
^\gamma \right) \right) \delta _\gamma + \\
\\
\qquad \qquad \qquad \quad +\left( \frac 1{a^2}(y_\alpha \delta _\beta
^\gamma -y_\beta \delta _\alpha ^\gamma )-\mathcal{R}_{\alpha \beta }^\gamma
\right) \mathcal{V}_\gamma , \\
\\
\mathbf{N}_{\Bbb{F}_0}(\delta _\alpha ,\mathcal{V}_\beta )=\left(
\frac 1{ \mathcal{F}^2}(y_\alpha \delta _\beta ^\gamma -y_\beta
\delta _\alpha ^\gamma
)-\frac{a^2}{\mathcal{F}^2}\mathcal{R}_{\alpha \beta }^\gamma
\right) \delta _\gamma - \\
\\
\qquad \qquad \qquad \quad -\left( t_{\alpha \beta }^\gamma +\frac
1{2 \mathcal{F}^2}\left( \delta _\beta (\mathcal{F}^2)\delta
_\alpha ^\gamma -\delta _\alpha (\mathcal{F}^2)\delta _\beta
^\gamma \right) \right)
\mathcal{V}_\gamma , \\
\\
\mathbf{N}_{\Bbb{F}_0}(\delta _\alpha ,\delta _\beta
)=-\frac{\mathcal{F}^2}{
a^2}\mathbf{N}_{\Bbb{F}_0}(\mathcal{V}_\alpha ,\mathcal{V}_\beta
).
\end{array}
\]
It follows that $\mathbf{N}_{\Bbb{F}_0}=0$ if and only if the relations
(1.3.34) are satisfied.\hfill
\hbox{\rlap{$\sqcap$}$\sqcup$}\newpage\

\subsubsection{\textbf{Linear connections on Lie algebroids}}

A linear connection on a Lie algebroid $(E,[,]_E,\sigma )$ is a map
\[
\mathcal{D}:\Gamma (E)\times \Gamma (E)\rightarrow \Gamma (E),
\]
which satisfies the rules\\i) $\mathcal{D}_{\rho +\omega }\eta
=\mathcal{D} _\rho \eta +\mathcal{D}_\omega \eta ,$\\ii)
$\mathcal{D}_\rho (\eta +\omega )=\mathcal{D}_\rho \eta
+\mathcal{D}_\rho \omega $\\iii) $\mathcal{D}_{f\rho }\eta
=f\mathcal{D}_{_\rho \eta },$\\iv) $\mathcal{D}_\rho (f\eta
)=(\sigma (\rho )f)\eta +f\mathcal{D}_{_\rho \eta }$\\for any
function $f\in C^\infty
(M)$ and $\rho ,\eta ,\omega \in \Gamma (E).$\\

For $\rho ,\eta \in \Gamma (E)$ the section $\mathcal{D}_\rho \eta \in
\Gamma (E)$ is called the \textit{covariant derivative} of the section $\eta
$ with respect to the section $\rho $. Let $\mathcal{N}$ be a nonlinear
connection. We have:

\begin{definition}
A linear connection $\mathcal{D}$ on Lie algebroids is called
$\mathcal{N}-$ linear connection if

i) $D$ preserves by parallelism the horizontal distribution $H\mathcal{T}E$.

ii) The tangent structure $J$ is absolute parallel with $\mathcal{D}$, that
is $\mathcal{D}J=0.$
\end{definition}

Consequently, the following properties hold:
\[
(\mathcal{D}_\rho \eta ^{\mathrm{h}})^{\mathrm{v}}=0,\quad
(\mathcal{D}_\rho \eta ^{\mathrm{v}})^{\mathrm{h}}=0,\quad
\mathcal{D}_\rho h=0,\quad \mathcal{ D}_\rho v=0,
\]
\[
\mathcal{D}_\rho (J\eta ^{\mathrm{h}})=J(\mathcal{D}_\rho \eta
^{\mathrm{h} }),\quad \mathcal{D}_\rho (J\eta
^{\mathrm{v}})=J(\mathcal{D}_\rho \eta ^{ \mathrm{v}}).
\]
If we denote
\[
\mathcal{D}_\rho ^{\mathrm{h}}\eta =\mathcal{D}_{\rho
^{\mathrm{h}}}\eta ,\quad \mathcal{D}_\rho ^{\mathrm{v}}\eta
=\mathcal{D}_{\rho ^{\mathrm{v} }}\eta ,
\]
then the following decomposition is obtained
\[
\mathcal{D}_\rho =\mathcal{D}_\rho ^{\mathrm{h}}+\mathcal{D}_\rho
^{\mathrm{v }},\quad \rho \in \Gamma (E).
\]
We remark that $\mathcal{D}^{\mathrm{h}}$ and
$\mathcal{D}^{\mathrm{v}}$ are not covariant derivative, because
$\mathcal{D}_\rho ^{\mathrm{h}}f=\sigma (\rho ^{\mathrm{h}})f\not
=\sigma (\rho )f,$ $\mathcal{D}_\rho ^{\mathrm{v} }f=\sigma (\rho
^{\mathrm{v}})f\not =\sigma (\rho )f,$ but, it still preserves
many properties of $\mathcal{D}$. Indeed, $\mathcal{D}^{\mathrm{h}
} $ and $\mathcal{D}^{\mathrm{v}}$ satisfy the Leibniz rule, and
$\mathcal{D} ^{\mathrm{h}}$ and $\mathcal{D}^{\mathrm{v}}$ will be
called the $\mathrm{h} - $\textit{covariant derivation} and
$\mathrm{v}-$\textit{covariant derivation}, respectively. Using
the fact that $\mathcal{D}_{\delta _\alpha }=\mathcal{D}_{\delta
_\alpha }^{\mathrm{h}}$, $\mathcal{D}_{\mathcal{V} _\alpha
}=\mathcal{D}_{\mathcal{V}_\alpha }^{\mathrm{v}}$ we get

\begin{proposition}
In the adapted basis $\{\delta _\alpha ,\mathcal{V}_\alpha \}$ a
$\mathcal{N} -$linear connection can be uniquely represented in
the form
\[
\mathcal{D}_{\delta _\beta }^{\mathrm{h}}\delta _\alpha =F_{\alpha
\beta }^\gamma \delta _\gamma ,\quad \mathcal{D}_{\delta _\beta
}^{\mathrm{v}} \mathcal{V}_\alpha =F_{\alpha \beta }^\gamma
\mathcal{V}_\gamma ,
\]
\[
\mathcal{D}_{\mathcal{V}_\beta }^{\mathrm{h}}\delta _\alpha
=C_{\alpha \beta }^\gamma \delta _\gamma ,\quad
\mathcal{D}_{\mathcal{V}_\beta }^{\mathrm{v}} \mathcal{V}_\alpha
=C_{\alpha \beta }^\gamma \mathcal{V}_\gamma .
\]
\end{proposition}

The system of functions ($F_{\alpha \beta }^\gamma (x,y)$, $C_{\alpha \beta
}^\gamma (x,y)$) represent the local coefficients of $\mathrm{h}-$covariant
derivation and of $\mathrm{v}-$covariant derivation, respectively. Under a
change of coordinates (1.2.13) the coefficients satisfy the rules
\[
F_{\alpha \beta }^\gamma =F_{\alpha ^{\prime }\beta ^{\prime }}^{\gamma
^{\prime }}A_\alpha ^{\alpha ^{\prime }}A_\beta ^{\beta ^{\prime }}A_{\gamma
^{\prime }}^\gamma +\sigma _\alpha ^i\frac{\partial A_\beta ^{\gamma
^{\prime }}}{\partial x^i}A_{\gamma ^{\prime }}^\gamma ,
\]
\[
C_{\alpha \beta }^\gamma =C_{\alpha ^{\prime }\beta ^{\prime }}^{\gamma
^{\prime }}A_\alpha ^{\alpha ^{\prime }}A_\beta ^{\beta ^{\prime }}A_{\gamma
^{\prime }}^\gamma .
\]
Let us consider a $d-$tensor $T$ in the local adapted basis, given by
\[
T=T_{\gamma \varepsilon }^{\alpha \beta }\delta _\alpha \otimes
\mathcal{V} _\beta \otimes \mathcal{X}^\gamma \otimes \delta
\mathcal{V}^\varepsilon
\]
A $d-$tensor means a tensor on $\mathcal{T}E$, whose components, under a
change of coordinates on $\mathcal{T}E$ behave like the components of a
tensor field on the base manifold $E$. Its covariant derivative with respect
to $\xi =\xi ^{\mathrm{h}}+\xi ^{\mathrm{v}}=\xi ^\alpha \delta _\alpha
+\zeta ^\beta \mathcal{V}_\beta $ is given by
\[
\mathcal{D}_\xi T=\left( \xi ^\kappa T_{\gamma \varepsilon /\kappa }^{\alpha
\beta }+\zeta ^\nu T_{\gamma \varepsilon }^{\alpha \beta }/_\nu \right)
\delta _\alpha \otimes \mathcal{V}_\beta \otimes \mathcal{X}^\gamma \otimes
\delta \mathcal{V}^\varepsilon
\]
where we have the $\mathrm{h}-$covariant derivative
$\mathcal{D}_\xi ^{ \mathrm{h}}T=\xi ^\kappa T_{\gamma \varepsilon
/\kappa }^{\alpha \beta }\delta _\alpha \otimes \mathcal{V}_\beta
\otimes \mathcal{X}^\gamma \otimes \delta \mathcal{V}^\varepsilon
$ with
\[
T_{\gamma \varepsilon /\kappa }^{\alpha \beta }=\sigma _\kappa
^i\frac{
\partial T_{\gamma \varepsilon }^{\alpha \beta }}{\partial x^i}-\mathcal{N}
_\kappa ^\tau \frac{\partial T_{\gamma \varepsilon }^{\alpha \beta
}}{
\partial y^\tau }+F_{\tau \kappa }^\alpha T_{\gamma \varepsilon }^{\tau
\beta }+F_{\tau \kappa }^\beta T_{\gamma \varepsilon }^{\alpha \tau
}-F_{\varepsilon \kappa }^\tau T_{\gamma \tau }^{\alpha \beta }-F_{\gamma
\kappa }^\tau T_{\tau \varepsilon }^{\alpha \beta }
\]
and $^{\prime \prime }\prime ^{\prime \prime }$ is the operator of
$\mathrm{h }-$ covariant derivative.\\The $\mathrm{v}-$covariant
derivative of $T$ is $ \mathcal{D}_\xi ^{\mathrm{v}}T=\xi ^\kappa
T_{\gamma \varepsilon }^{\alpha \beta }/_\kappa \delta _\alpha
\otimes \mathcal{V}_\beta \otimes \mathcal{X} ^\gamma \otimes
\delta \mathcal{V}^\varepsilon $ with
\[
T_{\gamma \varepsilon }^{\alpha \beta }/_\kappa =\frac{\partial T_{\gamma
\varepsilon }^{\alpha \beta }}{\partial y^\kappa }+C_{\tau \kappa }^\alpha
T_{\gamma \varepsilon }^{\tau \beta }+C_{\tau \kappa }^\beta T_{\gamma
\varepsilon }^{\alpha \tau }-C_{\varepsilon \kappa }^\tau T_{\gamma \tau
}^{\alpha \beta }-C_{\gamma \kappa }^\tau T_{\tau \varepsilon }^{\alpha
\beta }
\]
where $^{\prime \prime }/^{\prime \prime }$ is the operator of
$\mathrm{v}-$ covariant derivative.\newpage\

\subsubsection{\textbf{Torsion and curvature of a }$\mathcal{N}-$\textbf{
linear connection}}

The torsion tensor of a $\mathcal{N}-$linear connection is defined as usual
\[
T(\xi ,\omega )=\mathcal{D}_\xi \omega -\mathcal{D}_\omega \xi -[\xi ,\omega
]_{\mathcal{T}E}.
\]
As in the case of tangent bundle we have:

\begin{proposition}
The torsion of a $\mathcal{N}-$linear connection is completely determined by
the following five components
\begin{equation}
\left\{
\begin{array}{l}
\mathrm{h}T(\mathrm{h}\xi ,\mathrm{h}\omega )=\mathcal{D}_\xi
^{\mathrm{h}} \mathrm{h}\omega -\mathcal{D}_\omega
^{\mathrm{h}}\mathrm{h}\xi -\mathrm{h}[
\mathrm{h}\xi ,\mathrm{h}\omega ]_{\mathcal{T}E}, \\
\mathrm{v}T(\mathrm{h}\xi ,\mathrm{h}\omega
)=-\mathrm{v}[\mathrm{h}\xi ,
\mathrm{h}\omega ]_{\mathcal{T}E}, \\
\mathrm{h}T(\mathrm{h}\xi ,\mathrm{v}\omega )=-\mathcal{D}_\omega
^{\mathrm{v }}\mathrm{h}\xi -\mathrm{h}[\mathrm{h}\xi
,\mathrm{v}\omega ]_{\mathcal{T}E},
\\
\mathrm{v}T(\mathrm{h}\xi ,\mathrm{v}\omega )=\mathcal{D}_\xi
^{\mathrm{h}} \mathrm{v}\omega -\mathrm{v}[\mathrm{h}\xi
,\mathrm{v}\omega ]_{\mathcal{T}
E}, \\
\mathrm{v}T(\mathrm{v}\xi ,\mathrm{v}\omega )=\mathcal{D}_\xi
^{\mathrm{v}} \mathrm{v}\omega -\mathcal{D}_\omega
^{\mathrm{v}}\mathrm{v}\xi -\mathrm{v}[ \mathrm{v}\xi
,\mathrm{v}\omega ]_{\mathcal{T}E}.
\end{array}
\right.  \tag{1.3.35}
\end{equation}
\end{proposition}

With respect to the adapted basis the components of torsion are given by
\[
\left\{
\begin{array}{l}
hT(\delta _\beta ,\delta _\alpha )=T_{\alpha \beta }^\gamma \delta _\gamma
=\left( F_{\alpha \beta }^\gamma -F_{\beta \alpha }^\gamma -L_{\alpha \beta
}^\gamma \right) \delta _\gamma , \\
vT(\delta _\beta ,\delta _\alpha )=\mathcal{R}_{\alpha \beta }^\gamma
\mathcal{V}_\gamma =\left( \sigma _\beta ^i\frac{\partial \mathcal{N}_\alpha
^\gamma }{\partial x^i}-\sigma _\alpha ^i\frac{\partial \mathcal{N}_\beta
^\gamma }{\partial x^i}-\mathcal{N}_\beta ^\varepsilon \frac{\partial
\mathcal{N}_\alpha ^\gamma }{\partial y^\varepsilon }+\mathcal{N}_\alpha
^\varepsilon \frac{\partial \mathcal{N}_\beta ^\gamma }{\partial
y^\varepsilon }+L_{\alpha \beta }^\varepsilon \mathcal{N}_\varepsilon
^\gamma \right) \mathcal{V}_\gamma , \\
hT(\mathcal{V}_\beta ,\delta _\alpha )=C_{\alpha \beta }^\gamma \delta
_\gamma , \\
vT(\mathcal{V}_\beta ,\delta _\alpha )=P_{\alpha \beta }^\gamma
\mathcal{V} _\gamma =\left( \frac{\partial N_\beta ^\gamma
}{\partial y^\alpha }
-F_{\alpha \beta }^\gamma \right) \mathcal{V}_\gamma , \\
vT(\mathcal{V}_\beta ,\mathcal{V}_\alpha )=S_{\alpha \beta
}^\gamma \mathcal{ V}_\gamma =\left( C_{\alpha \beta }^\gamma
-C_{\beta \alpha }^\gamma \right) \mathcal{V}_\gamma .
\end{array}
\right.
\]

The curvature of a $\mathcal{N}-$linear connection is defined by
\[
R(\xi ,\omega )\varphi =\mathcal{D}_\xi \mathcal{D}_\omega \varphi
-\mathcal{ D}_\omega \mathcal{D}_\xi \varphi -\mathcal{D}_{[\xi
,\omega ]_{\mathcal{T} E}}\varphi .
\]

\begin{proposition}
The tensor of curvature has three essential components
\begin{equation}
\left\{
\begin{array}{c}
R(\delta _\gamma ,\delta _\beta )\delta _\alpha =R_{\alpha \beta \gamma
}^\varepsilon \delta _\varepsilon , \\
R(\mathcal{V}_\gamma ,\delta _\beta )\delta _\alpha =P_{\alpha \beta \gamma
}^\varepsilon \delta _\varepsilon , \\
R(\mathcal{V}_\gamma ,\mathcal{V}_\beta )\delta _\alpha =S_{\alpha \beta
\gamma }^\varepsilon \delta _\varepsilon ,
\end{array}
\right.  \tag{1.3.36}
\end{equation}
given by
\begin{eqnarray*}
R_{\gamma \beta \alpha }^\delta =\sigma _\alpha ^i\frac{\partial F_{\gamma
\beta }^\delta }{\partial x^i}-\mathcal{N}_\alpha ^\varepsilon \frac{%
\partial F_{\gamma \beta }^\delta }{\partial y^\varepsilon }-\sigma _\beta ^i
\frac{\partial F_{\gamma \alpha }^\delta }{\partial x^i}+N_\beta
^\varepsilon \frac{\partial F_{\gamma \alpha }^\delta }{\partial
y^\varepsilon }+ \\
\ \ +F_{\rho \alpha }^\delta F_{\gamma \beta }^\rho -F_{\rho \beta }^\delta
F_{\gamma \alpha }^\rho -L_{\beta \alpha }^\rho F_{\gamma \rho }^\delta
+C_{\gamma \varepsilon }^\delta \mathcal{R}_{\beta \alpha }^\varepsilon .
\end{eqnarray*}
\[
P_{\alpha \beta \gamma }^\gamma =\frac{\partial F_{\alpha \beta
}^\gamma }{
\partial y^\delta }-\sigma _\beta ^i\frac{\partial C_{\alpha \delta }^\gamma
}{\partial x^i}+N_\beta ^\varepsilon \frac{\partial C_{\alpha
\delta }^\gamma }{\partial y^\varepsilon }+F_{\alpha \beta
}^\varepsilon C_{\varepsilon \delta }^\gamma -C_{\alpha \delta
}^\varepsilon F_{\varepsilon \beta }^\gamma -\frac{\partial
N_\beta ^\varepsilon }{
\partial y^\delta }C_{\alpha \varepsilon }^\delta .
\]
\[
S_{\gamma \beta \alpha }^\delta =\frac{\partial C_{\gamma \beta
}^\delta }{
\partial y^\alpha }-\frac{\partial C_{\gamma \alpha }^\delta }{\partial
y^\beta }+\ C_{\rho \alpha }^\delta C_{\gamma \beta }^\rho -C_{\rho \beta
}^\delta C_{\gamma \alpha }^\rho .
\]
\end{proposition}

\begin{proposition}
The Ricci identities have the following form
\begin{eqnarray*}
X_{/\beta /\gamma }^\alpha -X_{/\gamma /\beta }^\alpha =R_{\rho
\beta \gamma }^\alpha X^\rho -T_{\beta \gamma }^\rho X_{/\rho
}^\alpha -\mathcal{R}
_{\beta \alpha }^\rho X^\alpha /_\rho - \\
\ -L_{\beta \gamma }^\rho \sigma _\rho ^i\frac{\partial X^\alpha
}{\partial x^i}+\left( \sigma _\gamma ^i\frac{\partial \sigma
_\beta ^j}{\partial x^i} -\sigma _\beta ^i\frac{\partial \sigma
_\gamma ^j}{\partial x^i}\right) \frac{\partial X^\alpha
}{\partial x^j},
\end{eqnarray*}
\[
X_{/\beta }^\alpha /_\gamma -X^\alpha /\gamma _{/\beta }=P_{\rho \beta
\gamma }^\alpha X^\rho -C_{\beta \gamma }^\rho X_{/\rho }^\alpha -P_{\beta
\gamma }^\rho X^\alpha /_\rho ,
\]
\[
X^\alpha /_\beta /_\gamma -X^\alpha /_\gamma /_\beta =S_{\rho \beta \gamma
}^\alpha X^\rho -S_{\beta \gamma }^\rho X^\alpha /_\rho .
\]
\end{proposition}

\begin{definition}
A $\mathcal{N}-$linear connection is called of Cartan type if
\begin{equation}
\mathcal{D}_\xi ^{\mathrm{h}}C=0,\quad \mathcal{D}_\xi
^{\mathrm{v}}C= \mathrm{v}\xi ,  \tag{1.3.37}
\end{equation}
where $C=y^\alpha \mathcal{V}_\alpha $ is the \textit{\ }Euler section.
\end{definition}

By direct computation, it results that a $\mathcal{N}-$linear
connection $ \mathcal{D}$ is of Cartan type if and only if
\begin{equation}
\mathcal{N}_\beta ^\alpha =F_{\varepsilon \beta }^\alpha y^\varepsilon
,\quad y^\varepsilon C_{\varepsilon \beta }^\alpha =0.  \tag{1.3.38}
\end{equation}
Introducing these relations into the coefficients expression of the
curvature, we obtain the following result:

\begin{proposition}
A $\mathcal{N}-$linear connection of Cartan type has the properties
\begin{equation}
\mathcal{R}_{\beta \gamma }^\alpha =R_{\varepsilon \beta \gamma }^\alpha
y^\varepsilon -L_{\beta \gamma }^\alpha ,\quad P_{\beta \gamma }^\alpha
=P_{\varepsilon \beta \gamma }^\alpha y^\varepsilon ,\quad S_{\beta \gamma
}^\alpha =S_{\varepsilon \beta \gamma }^\alpha y^\varepsilon .  \tag{1.3.39}
\end{equation}
\end{proposition}

\newpage\

\subsection{\textbf{Dynamical covariant derivative and metric non-linear
connection on} \textbf{Lie algebroid }$\mathcal{T}E$}

In this section we will introduce the notion of dynamical covariant
derivative on Lie algebroids as a tensor derivation and study the
compatibility between nonlinear connection and a pseudo-Riemannian metric.

\begin{definition}
A map $\nabla :\frak{T}(\mathcal{T}E\backslash \{0\})\rightarrow
\frak{T}( \mathcal{T}E\backslash \{0\})$ is said to be a tensor
derivation on $ \mathcal{T}E\backslash \{0\}$ if the following
conditions are satisfied:\\i) $\nabla $ is $\Bbb{R}$-linear\\ii)
$\nabla $ is type preserving, i.e. $ \nabla
(\frak{T}_s^r(\mathcal{T}E\backslash \{0\})\subset \frak{T}_s^r(
\mathcal{T}E\backslash \{0\})$, for each $(r,s)\in \Bbb{N}\times \Bbb{N.}$\\
iii) $\nabla $ obeys the Leibnitz rule $\nabla (P\otimes S)=\nabla
P\otimes S+P\otimes \nabla S$, for any tensors $P,S$ on
$\mathcal{T}E\backslash \{0\}. $\\iv) $\nabla \,$commutes with any
contractions, where $\frak{T} _{\bullet }^{\bullet
}(\mathcal{T}E\backslash \{0\})$ is the space of tensors on
$\mathcal{T}E\backslash \{0\}.$
\end{definition}

For a semispray $\mathcal{S}$ we consider the $\Bbb{R}$-linear map
\[
\nabla _0:\Gamma (\mathcal{T}E\backslash \{0\})\rightarrow \Gamma
(\mathcal{T }E\backslash \{0\}),
\]
given by
\begin{equation}
\nabla _0\rho =\mathrm{h}[\mathcal{S},\mathrm{h}\rho
]_{\mathcal{T}E}+ \mathrm{v}[\mathcal{S},\mathrm{v}\rho
]_{\mathcal{T}E},\quad \forall \rho \in \Gamma
(\mathcal{T}E\backslash \{0\}).  \tag{1.4.1}
\end{equation}
It results that
\begin{equation}
\nabla _0(f\rho )=\mathcal{S}(f)\rho +f\nabla _0\rho ,\quad \forall f\in
C^\infty (E),\ \rho \in \Gamma (\mathcal{T}E\backslash \{0\}).  \tag{1.4.2}
\end{equation}
Any tensor derivation on $\mathcal{T}E\backslash \{0\}$ is
completely determined by its actions on smooth functions and
sections on $\mathcal{T} E\backslash \{0\}$ (see \cite{Si}
generalized Willmore's theorem, p. 1217). Therefore there exists a
unique tensor derivation $\nabla $ on $\mathcal{T} E\backslash
\{0\}$ such that
\[
\nabla \mid _{C^\infty (E)}=\mathcal{S},\quad \nabla \mid _{\Gamma
(\mathcal{ T}E\backslash \{0\})}=\nabla _0.
\]
We will call the tensor derivation $\nabla $, the \textit{dynamical
covariant derivative} induced by the semispray $\mathcal{S}$ and a nonlinear
connection $N$.

\begin{proposition}
The following formulas hold
\begin{equation}
[\mathcal{S},\mathcal{V}_\beta ]_{\mathcal{T}E}=-\delta _\beta -\left(
\mathcal{N}_\beta ^\alpha +\frac{\partial \mathcal{S}^\alpha }{\partial
y^\beta }\right) \mathcal{V}_\alpha ,  \tag{1.4.3}
\end{equation}
\begin{equation}
[\mathcal{S},\delta _\beta ]_{\mathcal{T}E}=\left(
\mathcal{N}_\beta ^\alpha -L_{\beta \varepsilon }^\alpha
y^\varepsilon \right) \delta _\alpha + \mathcal{R}_\beta ^\gamma
\mathcal{V}_\gamma ,  \tag{1.4.4}
\end{equation}
where
\begin{equation}
\mathcal{R}_\beta ^\gamma =-\sigma _\beta ^i\frac{\partial
\mathcal{S} ^\gamma }{\partial x^i}-\mathcal{S}(\mathcal{N}_\beta
^\gamma )+\mathcal{N} _\beta ^\alpha \mathcal{N}_\alpha ^\gamma
+\mathcal{N}_\beta ^\alpha \frac{
\partial \mathcal{S}^\gamma }{\partial y^\alpha }+\mathcal{N}_\varepsilon
^\gamma L_{\alpha \beta }^\varepsilon y^\alpha .  \tag{1.4.5}
\end{equation}
\end{proposition}

The action of the dynamical covariant derivative on the Berwald basis is
given by
\[
\nabla \mathcal{V}_\beta =\mathrm{v}[\mathcal{S},\mathcal{V}_\beta
]_{ \mathcal{T}E}=-\left( \mathcal{N}_\beta ^\alpha
+\frac{\partial \mathcal{S} ^\alpha }{\partial y^\beta }\right)
\mathcal{V}_\alpha
\]
\[
\nabla \delta _\beta =\mathrm{h}[\mathcal{S},\delta _\beta
]_{\mathcal{T} E}=\left( \mathcal{N}_\beta ^\alpha -L_{\beta
\varepsilon }^\alpha y^\varepsilon \right) \delta _\alpha .
\]
It is not difficult to extend the action of $\nabla $ to the algebra of
tensors by requiring for $\nabla $ to preserve the tensor product. For a
pseudo-Riemannian metric $g$ on $V\mathcal{T}E$ (i.e. a $(2,0)$-type
symmetric tensor $g=g_{\alpha \beta }(x,y)\mathcal{V}^\alpha \otimes
\mathcal{V}^\beta $ of rank $m$ on $V\mathcal{T}E$) we have
\begin{equation}
(\nabla g)(\rho _1,\rho _2)=\mathcal{S}(g(\rho _1,\rho _2))-g(\nabla \rho
_1,\rho _2)-g(\rho _1,\nabla \rho _2),  \tag{1.4.6}
\end{equation}
and in local coordinates we get
\begin{equation}
g_{\alpha \beta /}:=(\nabla g)(\mathcal{V}_\alpha
,\mathcal{V}_\beta )= \mathcal{S}(g_{\alpha \beta })+g_{\gamma
\beta }\left( \mathcal{N}_\alpha ^\gamma +\frac{\partial
\mathcal{S}^\gamma }{\partial y^\alpha }\right) +g_{\gamma \alpha
}\left( \mathcal{N}_\beta ^\gamma +\frac{\partial \mathcal{
S}^\gamma }{\partial y^\beta }\right) .  \tag{1.4.7}
\end{equation}

\begin{definition}
The nonlinear connection $N$ is called metric or compatible with the metric
tensor $g$ if $\bigtriangledown g=0$, that is
\begin{equation}
\mathcal{S}(g(\rho _1,\rho _2))=g(\nabla \rho _1,\rho _2)+g(\rho _1,\nabla
\rho _2).  \tag{1.4.8}
\end{equation}
\end{definition}

If $\mathcal{S}$ be a semispray, $\mathcal{N}$ a nonlinear
connection and $ \nabla $ the dynamical covariant derivative
induced by $(\mathcal{S}, \mathcal{N})$, then we set:

\begin{proposition}
The nonlinear connection $\widetilde{\mathcal{N}}$ with the coefficients
given by
\begin{equation}
\widetilde{\mathcal{N}}_\beta ^\alpha =\mathcal{N}_\beta ^\alpha -\frac
12g^{\alpha \gamma }g_{\gamma \beta /},  \tag{1.4.9}
\end{equation}
is a metric nonlinear connection.
\end{proposition}

{\textbf{Proof}. Since $N_\beta ^\alpha $ are the coefficients of a
nonlinear connection and $g^{\alpha \gamma }g_{\gamma \beta /},$ are the
components of a tensor of type (1,1) it results that $\widetilde{N}_\beta
^\alpha $ are also the coefficients of a nonlinear connection. We consider
the dynamical covariant derivative induced by $(S,\widetilde{N})$ and we
have
\begin{eqnarray*}
(\nabla g)(\mathcal{V}_\alpha ,\mathcal{V}_\beta )
&=&\mathcal{S}(g_{\alpha \beta })+g_{\gamma \beta }\left(
\widetilde{\mathcal{N}}_\alpha ^\gamma + \frac{\partial
\mathcal{S}^\gamma }{\partial y^\alpha }\right) +g_{\gamma \alpha
}\left( \widetilde{\mathcal{N}}_\beta ^\gamma +\frac{\partial
\mathcal{S}^\gamma }{\partial y^\beta }\right) \\
\ &=&\mathcal{S}(g_{\alpha \beta })+g_{\gamma \beta }\left(
\mathcal{N} _\alpha ^\gamma +\frac{\partial \mathcal{S}^\gamma
}{\partial y^\alpha } \right) +g_{\gamma \alpha }\left(
\mathcal{N}_\beta ^\gamma +\frac{\partial
\mathcal{S}^\gamma }{\partial y^\beta }\right) - \\
\ &&-{g_{\gamma \beta }\frac 12g^{\gamma \varepsilon }g_{\varepsilon \alpha
/}-}g_{\gamma \alpha }\frac 12g^{\gamma \varepsilon }g_{\varepsilon \beta /}=
\\
&=&{g_{\alpha \beta /}-\frac 12g_{\alpha \beta /}-\frac 12g_{\alpha \beta
/}=0,}
\end{eqnarray*}
that is the connection $\widetilde{\mathcal{N}}$ is metric.\hfill
\hbox{\rlap{$\sqcap$}$\sqcup$} }

For the particular case of the tangent bundle see \cite{Bu1, Bu3}\newpage\

\subsubsection{\textbf{The case of SODE connection}}

A semispray (\textit{SODE}) given by $\mathcal{S}=y^\alpha
\mathcal{X} _\alpha +\mathcal{S}^\alpha \mathcal{V}_\alpha $
determines an associated nonlinear connection
\[
\mathcal{N}=-\mathcal{L}_{\mathcal{S}}J,
\]
with local coefficients (1.3.25)
\[
\mathcal{N}_\alpha ^\beta =\frac 12\left( -\frac{\partial \mathcal{S}^\beta
}{\partial y^\alpha }+y^\varepsilon L_{\alpha \varepsilon }^\beta \right) .
\]

\begin{proposition}
The following equations hold
\begin{equation}
[\mathcal{S},\mathcal{V}_\beta ]_{\mathcal{T}E}=-\delta _\beta +\left(
\mathcal{N}_\beta ^\alpha -L_{\beta \varepsilon }^\alpha y^\varepsilon
\right) \mathcal{V}_\alpha ,  \tag{1.4.10}
\end{equation}
\begin{equation}
[\mathcal{S},\delta _\beta ]_{\mathcal{T}E}=\left(
\mathcal{N}_\beta ^\alpha -L_{\beta \varepsilon }^\alpha
y^\varepsilon \right) \delta _\alpha + \mathcal{R}_\beta ^\alpha
\mathcal{V}_\alpha ,  \tag{1.4.11}
\end{equation}
where
\begin{equation}
\mathcal{R}_\beta ^\alpha =-\sigma _\beta ^i\frac{\partial
\mathcal{S} ^\alpha }{\partial x^i}-\mathcal{S}(\mathcal{N}_\beta
^\alpha )-\mathcal{N} _\gamma ^\alpha \mathcal{N}_\beta ^\gamma
+(L_{\varepsilon \beta }^\gamma \mathcal{N}_\gamma ^\alpha
+L_{\gamma \varepsilon }^\alpha \mathcal{N}_\beta ^\gamma
)y^\varepsilon .  \tag{1.4.12}
\end{equation}
\end{proposition}

The dynamical covariant derivative induced by $\mathcal{S}$ and associated
nonlinear connection is characterized by
\begin{equation}
\nabla \mathcal{V}_\beta =\mathrm{v}[\mathcal{S},\mathcal{V}_\beta
]_{ \mathcal{T}E}=\left( \mathcal{N}_\beta ^\alpha -L_{\beta
\varepsilon }^\alpha y^\varepsilon \right) \mathcal{V}_\alpha
=-\frac 12\left( \frac{
\partial \mathcal{S}^\alpha }{\partial y^\beta }+L_{\beta \varepsilon
}^\alpha y^\varepsilon \right) \mathcal{V}_\alpha ,  \tag{1.4.13}
\end{equation}
\[
\nabla \delta _\beta =\mathrm{h}[\mathcal{S},\delta _\beta
]_{\mathcal{T} E}=\left( \mathcal{N}_\beta ^\alpha -L_{\beta
\varepsilon }^\alpha y^\varepsilon \right) \delta _\alpha .
\]
\begin{equation}
g_{\alpha \beta /}:=(\nabla g)(\mathcal{V}_\alpha
,\mathcal{V}_\beta )= \mathcal{S}(g_{\alpha \beta })-g_{\gamma
\beta }\mathcal{N}_\alpha ^\gamma -g_{\gamma \alpha
}\mathcal{N}_\beta ^\gamma +\left( g_{\gamma \beta }L_{\alpha
\varepsilon }^\gamma +g_{\gamma \alpha }L_{\beta \varepsilon
}^\gamma \right) y^\varepsilon ,  \tag{1.4.14}
\end{equation}
which is equivalent to
\begin{equation}
(\nabla g)(\mathcal{V}_\alpha ,\mathcal{V}_\beta
)=\mathcal{S}(g_{\alpha \beta })+\frac 12\frac{\partial
\mathcal{S}^\gamma }{\partial y^\alpha } g_{\gamma \beta }+\frac
12\frac{\partial \mathcal{S}^\gamma }{\partial y^\beta }g_{\gamma
\alpha }+\frac 12\left( g_{\gamma \beta }L_{\alpha \varepsilon
}^\gamma +g_{\gamma \alpha }L_{\beta \varepsilon }^\gamma \right)
y^\varepsilon .  \tag{1.4.15}
\end{equation}

\begin{definition}
The Jacobi endomorphism is given by
\[
\Phi =\mathrm{v}[\mathcal{S},\mathrm{h}\rho ]_{\mathcal{T}E}.
\]
\end{definition}

Locally, from (1.4.11) we obtain that $\Phi =\mathcal{R}_\beta ^\alpha
\mathcal{V}_\alpha \otimes \mathcal{X}^\beta ,$ where $\mathcal{R}_\beta
^\alpha $ is given by (1.4.12) and represent the local coefficients of the
Jacobi endomorphism.

\begin{proposition}
The following result holds
\[
\Phi =i_{\mathcal{S}}\Omega +\mathrm{v}[\mathrm{v}\mathcal{S},\mathrm{h}\rho
]_{\mathcal{T}E}.
\]
\end{proposition}

\textbf{Proof}. Indeed, $\Phi (\rho
)=\mathrm{v}[\mathcal{S},\mathrm{h}\rho
]_{\mathcal{T}E}=\mathrm{v}[\mathrm{h}\mathcal{S},\mathrm{h}\rho
]_{\mathcal{ T}E}+\mathrm{v}[\mathrm{v}\mathcal{S},\mathrm{h}\rho
]_{\mathcal{T}E}$ and $ \Omega (\mathcal{S},\rho
)=\mathrm{v}[\mathrm{h}\mathcal{S},\mathrm{h}\rho
]_{\mathcal{T}E},$ which yields $\Phi (\rho )=\Omega
(\mathcal{S},\rho )+
\mathrm{v}[\mathrm{v}\mathcal{S},\mathrm{h}\rho
]_{\mathcal{T}E}.$\hfill \hbox{\rlap{$\sqcap$}$\sqcup$}

If $\mathcal{S}$ is a spray, then the coefficients $\mathcal{S}^\alpha $ are
2-homogeneous with respect to the variables $y^\beta $ and it results
\[
2\mathcal{S}^\alpha =\frac{\partial \mathcal{S}^\alpha }{\partial
y^\beta } y^\beta =-2\mathcal{N}_\beta ^\alpha y^\beta +L_{\beta
\gamma }^\alpha y^\beta y^\gamma =-2\mathcal{N}_\beta ^\alpha
y^\beta .
\]
\[
\mathcal{S}=\mathrm{h}\mathcal{S}=y^\alpha \delta _\alpha ,\quad
\mathrm{v} \mathcal{S}=0,\quad \mathcal{N}_\beta ^\alpha
=\frac{\partial \mathcal{N} _\varepsilon ^\alpha }{\partial
y^\beta }y^\varepsilon +L_{\beta \varepsilon }^\alpha
y^\varepsilon ,
\]
which yields
\begin{equation}
\Phi =i_{\mathcal{S}}\Omega ,  \tag{1.4.16}
\end{equation}
and locally we get
\begin{equation}
\mathcal{R}_\beta ^\alpha =\mathcal{R}_{\varepsilon \beta }^\alpha
y^\varepsilon ,  \tag{1.4.17}
\end{equation}
which represents the local relation between the Jacobi endomorphism and the
curvature of the nonlinear connection.\newpage\

\subsubsection{\textbf{Lagrangian case}}

Let us consider a regular Lagrangian $L$ on $E$, that is the matrix
\[
g_{\alpha \beta }=\frac{\partial ^2L}{\partial y^\alpha \partial y^\beta },
\]
has constant rank $m$. The symplectic structure induced by the regular
Lagrangian is \cite{Ma2}
\[
\omega _L=\frac{\partial ^2L}{\partial y^\alpha \partial y^\beta
}\mathcal{V} ^\beta \wedge \mathcal{X}^\alpha +\frac 12\left(
\frac{\partial ^2L}{
\partial x^i\partial y^\beta }\sigma _\alpha ^i-\frac{\partial ^2L}{\partial
x^i\partial y^\alpha }\sigma _\beta ^i-\frac{\partial L}{\partial
y^\gamma } L_{\alpha \beta }^\gamma \right) \mathcal{X}^\alpha
\wedge \mathcal{X}^\beta .
\]
Let us consider the energy function given by
\[
E_L=y^\alpha \frac{\partial L}{\partial y^\alpha }-L,
\]
then the symplectic equation
\[
i_{\mathcal{S}}\omega _L=-d^EE_L,\quad \mathcal{S}\in \Gamma (\mathcal{T}E),
\]
and the regularity condition of the Lagrangian determine the components of
the semispray (1.3.29)
\[
\mathcal{S}^\varepsilon =g^{\varepsilon \beta }\left( \sigma
_\beta ^i\frac{
\partial L}{\partial x^i}-\sigma _\alpha ^i\frac{\partial ^2L}{\partial
x^i\partial y^\beta }y^\alpha -L_{\beta \alpha }^\theta y^\alpha
\frac{
\partial L}{\partial y^\theta }\right) ,
\]
where $g_{\alpha \beta }g^{\beta \gamma }=\delta _\alpha ^\gamma
$.\\The connection $\mathcal{N}$ determined by this semispray is
the \textit{ canonical nonlinear connection} induced by a regular
Lagrangian $L$. Its coefficients are given by
\begin{equation}
\mathcal{N}_\beta ^\alpha =\frac 12g^{\alpha \varepsilon }\left[
\mathcal{S} (g_{\varepsilon \beta })+\sigma _\beta
^i\frac{\partial ^2L}{\partial x^i\partial y^\varepsilon }-\sigma
_\varepsilon ^i\frac{\partial ^2L}{
\partial x^i\partial y^\beta }-L_{\beta \varepsilon }^\gamma \frac{\partial L
}{\partial y^\gamma }+\left( g_{\gamma \varepsilon }L_{\beta \theta }^\gamma
+g_{\gamma \beta }L_{\varepsilon \theta }^\gamma \right) y^\theta \right]
\tag{1.4.18}
\end{equation}

\begin{theorem}
The canonical nonlinear connection $\mathcal{N}$ induced by a regular
Lagrangian $L$ is a metric nonlinear connection.
\end{theorem}

\textbf{Proof}. Introducing the expression of the semispray (1.3.29) into
the equation (1.4.15) we obtain
\begin{eqnarray*}
\ (\nabla g)(\mathcal{V}_\alpha ,\mathcal{V}_\beta )
&=&y^\varepsilon \sigma _\varepsilon ^i\frac{\partial g_{\alpha
\beta }}{\partial x^i} +g^{\varepsilon \gamma }\left( \sigma
_\gamma ^i\frac{\partial L}{\partial x^i}-\sigma _\theta
^i\frac{\partial ^2L}{\partial x^i\partial y^\gamma } y^\theta
-L_{\gamma \tau }^\theta y^\tau \frac{\partial L}{\partial
y^\theta
}\right) \frac{\partial g_{\alpha \beta }}{\partial y^\varepsilon } \\
&&+\frac 12\left( g_{\gamma \beta }\frac{\partial g^{\gamma
\varepsilon }}{
\partial y^\alpha }+g_{\gamma \alpha }\frac{\partial g^{\gamma \varepsilon }
}{\partial y^\beta }\right) \left( \sigma _\varepsilon
^i\frac{\partial L}{
\partial x^i}-\sigma _\theta ^i\frac{\partial ^2L}{\partial x^i\partial
y^\varepsilon }y^\theta -L_{\varepsilon \tau }^\theta
\frac{\partial L}{
\partial y^\theta }y^\tau \right) \\
&&+\frac 12\left( \sigma _\beta ^i\frac{\partial ^2L}{\partial x^i\partial
y^\alpha }+\sigma _\alpha ^i\frac{\partial ^2L}{\partial x^i\partial y^\beta
}\right) -\sigma _\varepsilon ^i\frac{\partial g_{\alpha \beta }}{\partial
x^i}y^\varepsilon \\
&&-\frac 12\left( \sigma _\alpha ^i\frac{\partial ^2L}{\partial
x^i\partial y^\beta }+\sigma _\beta ^i\frac{\partial ^2L}{\partial
x^i\partial y^\alpha } \right) -\frac 12\frac{\partial L}{\partial
y^\varepsilon }\left( L_{\beta
\alpha }^\varepsilon +L_{\alpha \beta }^\varepsilon \right) \\
&&-\frac 12\left( g_{\gamma \beta }L_{\alpha \varepsilon }^\gamma +g_{\gamma
\alpha }L_{\beta \varepsilon }^\gamma \right) y^\varepsilon +\frac 12\left(
g_{\gamma \beta }L_{\alpha \varepsilon }^\gamma +g_{\gamma \alpha }L_{\beta
\varepsilon }^\gamma \right) y^\varepsilon .
\end{eqnarray*}
By direct computation, using the equalities
\[
g_{\gamma \beta }\frac{\partial g^{\gamma \varepsilon }}{\partial
y^\alpha } =-g^{\gamma \varepsilon }\frac{\partial g_{\gamma \beta
}}{\partial y^\alpha }=-g^{\gamma \varepsilon }\frac{\partial
g_{\alpha \beta }}{\partial y^\gamma },\quad L_{\alpha \beta
}^\theta =-L_{\beta \alpha }^\theta ,
\]
it results $\ (\nabla g)(\mathcal{V}_\alpha ,\mathcal{V}_\beta )=0,$ which
ends the proof. \hfill
\hbox{\rlap{$\sqcap$}$\sqcup$}

\begin{theorem}
The canonical nonlinear connection induced by a regular Lagrangian is a
unique connection which is metric and compatible with the symplectic
structure $\omega _L$, that is
\begin{equation}
\nabla g=0,  \tag{1.4.19}
\end{equation}
\begin{equation}
\omega _L(\mathrm{h}\rho ,\mathrm{h}\nu )=0,\quad \forall \rho ,\nu \in
\Gamma (\mathcal{T}E\backslash \{0\})  \tag{1.4.20}
\end{equation}
\end{theorem}

\textbf{Proof}. Using the equation $\mathcal{V}^\alpha =\delta
\mathcal{V} ^\alpha -\mathcal{N}_\beta ^\alpha \mathcal{X}^\beta $
it results
\[
\omega _L=g_{\alpha \beta }(\delta \mathcal{V}^\beta
-\mathcal{N}_\gamma ^\beta \mathcal{X}^\gamma )\wedge
\mathcal{X}^\alpha +\frac 12\left( \frac{
\partial ^2L}{\partial x^i\partial y^\beta }\sigma _\alpha ^i-\frac{\partial
^2L}{\partial x^i\partial y^\alpha }\sigma _\beta
^i-\frac{\partial L}{
\partial y^\varepsilon }L_{\alpha \beta }^\varepsilon \right) \mathcal{X}
^\beta \wedge \mathcal{X}^\alpha
\]
\[
=g_{\alpha \beta }\delta \mathcal{V}^\beta \wedge
\mathcal{X}^\alpha +\frac 12\left( g_{\alpha \gamma
}\mathcal{N}_\beta ^\gamma -g_{\beta \gamma } \mathcal{N}_\alpha
^\beta +\frac{\partial ^2L}{\partial x^i\partial y^\beta } \sigma
_\alpha ^i-\frac{\partial ^2L}{\partial x^i\partial y^\alpha
}\sigma _\beta ^i-\frac{\partial L}{\partial y^\varepsilon
}L_{\alpha \beta }^\varepsilon \right) \mathcal{X}^\beta \wedge
\mathcal{X}^\alpha
\]
\[
=g_{\alpha \beta }\delta \mathcal{V}^\beta \wedge
\mathcal{X}^\alpha +\frac 12\left( \mathcal{N}_{\alpha \beta
}-\mathcal{N}_{\beta \alpha }+\frac{
\partial ^2L}{\partial x^i\partial y^\beta }\sigma _\alpha ^i-\frac{\partial
^2L}{\partial x^i\partial y^\alpha }\sigma _\beta
^i-\frac{\partial L}{
\partial y^\varepsilon }L_{\alpha \beta }^\varepsilon \right) \mathcal{X}
^\beta \wedge \mathcal{X}^\alpha ,
\]
where $\mathcal{N}_{\alpha \beta }:=g_{\alpha \gamma }\mathcal{N}_\beta
^\gamma .$ We have that $\omega _L(\mathrm{h}\rho ,\mathrm{h}\nu )=0$ if and
only if the second part of the above relation vanishes, that is
\[
\mathcal{N}_{[\alpha \beta ]}=\frac 12(\mathcal{N}_{\alpha \beta
}-\mathcal{N }_{\beta \alpha })=\frac 12\left( \frac{\partial
^2L}{\partial x^i\partial y^\alpha }\sigma _\beta
^i-\frac{\partial ^2L}{\partial x^i\partial y^\beta } \sigma
_\alpha ^i+\frac{\partial L}{\partial y^\varepsilon }L_{\alpha
\beta }^\varepsilon \right) .
\]
It result that the skew symmetric part of $\mathcal{N}_{\alpha
\beta }$ is uniquely determined by the condition (1.4.20). The
symmetric part of $ \mathcal{N}_{\alpha \beta }$ is completely
determined by the metric condition (1.4.19). Indeed
\begin{eqnarray*}
\mathcal{S}(g_{\alpha \beta }) &=&g_{\gamma \beta }\mathcal{N}_\alpha
^\gamma +g_{\gamma \alpha }\mathcal{N}_\beta ^\gamma -\left( g_{\gamma \beta
}L_{\alpha \varepsilon }^\gamma +g_{\gamma \alpha }L_{\beta \varepsilon
}^\gamma \right) y^\varepsilon \\
\ &=&\mathcal{N}_{\beta \alpha }+\mathcal{N}_{\alpha \beta }-\left(
g_{\gamma \beta }L_{\alpha \varepsilon }^\gamma +g_{\gamma \alpha }L_{\beta
\varepsilon }^\gamma \right) y^\varepsilon \\
\ &=&2\mathcal{N}_{(\alpha \beta )}-\left( g_{\gamma \beta }L_{\alpha
\varepsilon }^\gamma +g_{\gamma \alpha }L_{\beta \varepsilon }^\gamma
\right) y^\varepsilon .
\end{eqnarray*}
that is
\[
2\mathcal{N}_{(\alpha \beta )}=\mathcal{S}(g_{\alpha \beta })+\left(
g_{\gamma \beta }L_{\alpha \varepsilon }^\gamma +g_{\gamma \alpha }L_{\beta
\varepsilon }^\gamma \right) y^\varepsilon .
\]
The equations (1.4.19) and (1.4.20) uniquely determine the coefficients of
the nonlinear connection
\begin{eqnarray*}
\mathcal{N}_\beta ^\gamma &=&g^{\gamma \alpha }\mathcal{N}_{\alpha \beta
}=g^{\gamma \alpha }(\mathcal{N}_{(\alpha \beta )}+\mathcal{N}_{[\alpha
\beta ]}) \\
&=&\frac 12g^{\gamma \alpha }\left[ \mathcal{S}(g_{\alpha \beta
})+\frac{
\partial ^2L}{\partial x^i\partial y^\alpha }\sigma _\beta ^i-\frac{\partial
^2L}{\partial x^i\partial y^\beta }\sigma _\alpha
^i-\frac{\partial L}{
\partial y^\varepsilon }L_{\beta \alpha }^\varepsilon \right] \\
&&+\frac 12g^{\gamma \alpha }\left( g_{\gamma \beta }L_{\alpha \varepsilon
}^\gamma +g_{\gamma \alpha }L_{\beta \varepsilon }^\gamma \right)
y^\varepsilon
\end{eqnarray*}
Conversely, introducing (1.3.29) into (1.3.25) we have (1.4.18) which ends
the proof.

\begin{remark}
The invariant form of Helmholtz conditions on Lie algebroids is given by:
\begin{equation}
\begin{array}{l}
\mathcal{D}_\rho ^{\mathrm{v}}g(\nu ,\theta )=\mathcal{D}_\theta
^{\mathrm{v}
}g(\nu ,\rho ), \\
\nabla g=0, \\
g(\Phi \rho ,\nu )=g(\Phi \nu ,\rho ),
\end{array}
\tag{1.4.21}
\end{equation}
for $\nu ,\rho ,\theta \in \Gamma (E)$, where
$\mathcal{D}^{\mathrm{v}}$ is $ \mathrm{v}$-covariant derivative.
\end{remark}

In local coordinates we obtain the following equations
\[
\frac{\partial g_{\alpha \beta }}{\partial y^\varepsilon }=\frac{\partial
g_{\alpha \varepsilon }}{\partial y^\beta },
\]
\[
\mathcal{S}(g_{\alpha \beta })-g_{\gamma \beta }N_\alpha ^\gamma -g_{\gamma
\alpha }N_\beta ^\gamma =y^\varepsilon \left( g_{\gamma \beta
}L_{\varepsilon \alpha }^\gamma +g_{\gamma \alpha }L_{\varepsilon \beta
}^\gamma \right) ,
\]
\begin{equation}
\ g_{\alpha \gamma }\left( \sigma _\beta ^i\frac{\partial \mathcal{S}^\gamma
}{\partial x^i}+\mathcal{S}N_\beta ^\gamma +N_\beta ^\varepsilon
N_\varepsilon ^\gamma -(L_{\varepsilon \beta }^\delta N_\delta ^\gamma
+L_{\delta \varepsilon }^\gamma N_\beta ^\delta )y^\varepsilon \right) =
\tag{1.4.22}
\end{equation}
\[
\ g_{\beta \gamma }\left( \sigma _\alpha ^i\frac{\partial
\mathcal{S}^\gamma }{\partial x^i}+\mathcal{S}N_\alpha ^\gamma
+N_\alpha ^\varepsilon N_\varepsilon ^\gamma -(L_{\varepsilon
\alpha }^\delta N_\delta ^\gamma +L_{\delta \varepsilon }^\gamma
N_\alpha ^\delta )y^\varepsilon \right).
\]

In the case of standard Lie algebroid $(TM,[\cdot ,\cdot ],id)$ we obtain
the classical Helmholtz conditions \cite{Sa1}.\newpage\

\subsection{\textbf{The prolongation of a Lie algebroid over the vector
bundle projection of the dual bundle}}

Let $\tau :E^{*}\rightarrow M$ be the dual bundle of $\pi :E\rightarrow M$
and $(E,[\cdot ,\cdot ]_E,\sigma )$ a Lie algebroid structure over $M.$ \
One can construct a Lie algebroid structure over $E^{*}$, by taking the
prolongation of $(E,[\cdot ,\cdot ]_E,\sigma )$ over $\tau :E^{*}\rightarrow
M$ (see \cite{Hi}, \cite{Le}, \cite{Hr4}). This structure is given by:

$\bullet $ The associated vector bundle is $(\mathcal{T}E^{*},\tau _1,E^{*})$
where
\[
\mathcal{T}E^{*}={\cup }\mathcal{T}_{u^{*}}E^{*}, \quad {u^{*}\in
E^{*}}
\]
\ with
\[
\mathcal{T}_{u^{*}}E^{*}=\{(u_x,v_{u^{*}})\in E_x\times
T_{u^{*}}E^{*}|\sigma (u_x)=T_{u^{*}}\tau (v_{u^{*}}),\tau \left(
u^{*}\right) =x\in M\},
\]
and the projection $\tau _1:$ $\mathcal{T}E^{*}\rightarrow E^{*}$, $\tau
_1(u_x,v_{u^{*}})=u^{*}.$

$\bullet $ The Lie algebra structure $[\cdot ,\cdot ]_{\mathcal{T}E^{*}}$ on
$\ \Gamma (\mathcal{T}E^{*})$ is defined in the following way: if $\rho
_1,\rho _2\in \Gamma (\mathcal{T}E^{*})$ are such that $\rho
_i(u^{*})=(X_i(\tau \left( u^{*}\right) ),U_i(u^{*}))$ where $X_i\in \Gamma
(E),U_i\in \chi (E^{*})$ and $\sigma (X_i(\tau \left( u^{*}\right)
)=T_{u^{*}}\tau (U_i\left( u^{*}\right) )$, $\ i=1,2$, then
\[
\lbrack \rho _1,\rho
_2]_{\mathcal{T}E^{*}}(u^{*})=([X_1,X_2]_{\mathcal{T} E^{*}}(\tau
\left( u^{*}\right) ),[U_1,U_2]_{\mathcal{T}E^{*}}(u^{*})).
\]

$\bullet $ The anchor is the projection $\sigma ^1:\mathcal{T}%
E^{*}\rightarrow TE^{*}$, $\sigma ^1(u,v)=v$.

Notice that if $\mathcal{T}\tau :\mathcal{T}E^{*}\rightarrow E$,
$\mathcal{T} \tau (u,v)=u$\ then $(V\mathcal{T}E^{*},\tau
_{1|V\mathcal{T}E^{*}},E^{*})\ $ with
$V\mathcal{T}E^{*}=Ker\mathcal{T}\tau $ is a subbundle of
$(\mathcal{T} E^{*},\tau _1,E^{*})$, called the \textit{vertical
subbundle}. If $(q^i,\mu _\alpha )$ are local coordinates on
$E^{*}$ at $u^{*}$ and $\{s_\alpha \}$ is a local basis of
sections of $\pi :E\rightarrow M$ then a local basis of $ \Gamma
(\mathcal{T}E^{*})$ is $\{\mathcal{Q}_\alpha ,\mathcal{P}^\alpha
\}$ where
\begin{equation}
\mathcal{Q}_\alpha (u^{*})=\left( s_\alpha (\tau (u^{*})),\sigma _\alpha
^i\frac \partial {\partial q^i}|_{u^{*}}\right) ,\quad \mathcal{P}^\alpha
(u^{*})=\left( 0,\frac \partial {\partial \mu _\alpha }|_{u^{*}}\right) .
\tag{1.5.1}
\end{equation}
The structure functions on $\mathcal{T}E^{*}$ are given by the following
formulas
\begin{equation}
\sigma ^1(\mathcal{Q}_\alpha )=\sigma _\alpha ^i\frac \partial {\partial
q^i},\quad \sigma ^1(\mathcal{P}^\alpha )=\frac \partial {\partial \mu
_\alpha },  \tag{1.5.2}
\end{equation}
\begin{equation}
\lbrack \mathcal{Q}_\alpha ,\mathcal{Q}_\beta
]_{\mathcal{T}E^{*}}=L_{\alpha \beta }^\gamma \mathcal{Q}_\gamma
,\quad [\mathcal{Q}_\alpha ,\mathcal{P} ^\alpha
]_{\mathcal{T}E^{*}}=0,\quad [\mathcal{P}^\alpha
,\mathcal{P}^\beta ]_{\mathcal{T}E^{*}}=0,  \tag{1.5.3}
\end{equation}
and therefore
\begin{equation}
d^Eq^i=\sigma _\alpha ^i\mathcal{Q}^\alpha ,\ d^E\mu _\alpha
=\mathcal{P} _\alpha ,  \tag{1.5.4}
\end{equation}
\[
d^E\mathcal{Q}^\gamma =-\frac 12L_{\alpha \beta }^\gamma \mathcal{Q}^\alpha
\wedge \mathcal{Q}^\beta ,\quad d^E\mathcal{P}_\alpha =0,
\]
where $\{\mathcal{Q}^\alpha ,\mathcal{P}_\alpha \}$ is the dual
basis of $\{ \mathcal{Q}_\alpha ,\mathcal{P}^\alpha \}.$ Also, if
$\rho =\rho ^\alpha \mathcal{Q}_\alpha +\rho _\alpha
\mathcal{P}^\alpha $ is a section of $ \mathcal{T}E^{*}$, then
\[
\sigma ^1(\rho )=\sigma _\alpha ^i\rho ^\alpha \frac \partial {\partial
q^i}+\rho _\alpha \frac \partial {\partial \mu _\alpha }.
\]
If $u^{*}\in E^{*}$ and $(u_x,v_{u^{*}})\in E_x\times T_{u^{*}}E^{*}$ then
\[
\theta _E(u^{*})(u_x,v_{u^{*}})=u^{*}(u_x),
\]
is called the \textit{Liouville section}. The \textit{canonical symplectic
section} $\omega _E$ is defined by
\[
\omega _E=-d^E\theta _E,
\]
and it results that is a nondegenerate 2-section and $d^E\omega _E=0.$

On $E^{*}$ we have the similar concept of the vertical lift to that in $E$.
If $\alpha \in \Gamma (E^{*})$ we can define the vector field $\alpha ^v$ on
$E^{*}$ as follows
\[
\alpha ^v(u^{*})=\alpha (\tau (u^{*}))_{u^{*}}^v,\ u^{*}\in E^{*},
\]
where
\[
_{u^{*}}^v:E_{\tau \left( u^{*}\right) }^{*}\rightarrow T_{u^{*}}(E_{\tau
\left( u^{*}\right) }^{*}),
\]
is the canonical isomorphism between the vector spaces $E_{\tau
\left( u^{*}\right) }^{*}$ and $T_{u^{*}}(E_{\tau \left(
u^{*}\right) }^{*})$. Also, if $X$ is a section of $\pi
:E\rightarrow M$, there exists a unique vector field $X^{*c}$ on
$E^{*}$, called the complete lift of $X$ to $E^{*}$ , satisfying
the two following conditions:

(i) $X^{*c\text{ }}$is $\tau -$projectable on $\sigma (X),$

(ii) $X^{*c}(\stackrel{\wedge }{Y})=\widehat{\mathcal{L}_XY}$, for
all $Y\in \Gamma (E)$ (see \cite{Gu1}).\\If $X$ is a section of
$E$ then $\stackrel{ \wedge }{X}$ is the linear function
$\stackrel{\wedge }{X}\in C^\infty (E^{*})$ given by
\[
\stackrel{\wedge }{X}(u^{*})=u^{*}(X(\tau (u^{*}))),
\]
for all $u^{*}\in E^{*}$. Now, we may introduce the vertical lift
$\alpha ^{ \mathrm{v}}$ and the complete lift $X^{*\mathrm{c}}$ of
a section $\alpha \in \Gamma (E^{*})$ and a section $X\in \Gamma
(E)$ as the sections of $ \mathcal{T}E^{*}$ given by
\[
\alpha ^{\mathrm{v}}(u^{*})=(0,\alpha ^v(u^{*})),\quad
X^{*\mathrm{c} }=(X(\tau (u^{*}),X^{*c}(u^{*})),\quad u^{*}\in
E^{*}.
\]
The other canonical object on $\mathcal{T}E^{*}$ is \textit{%
Liouville-Hamilton section} $\mathcal{C}$ given by
\[
\mathcal{C}(u^{*})=(0,u_{u^{*}}^{*v}),\quad u^{*}\in E^{*}.
\]
In local coordinates it follows that the Liouville section is given by
\begin{equation}
\theta _E=\mu _\alpha \mathcal{Q}^\alpha ,  \tag{1.5.5}
\end{equation}
and we obtain
\begin{equation}
\omega _E=\mathcal{Q}^\alpha \wedge \mathcal{P}_\alpha +\frac 12\mu _\alpha
L_{\beta \gamma }^\alpha \mathcal{Q}^\beta \wedge \mathcal{Q}^\gamma .
\tag{1.5.6}
\end{equation}
The \textit{Liouville-Hamilton section }$\mathcal{C}$ has local expression
\begin{equation}
\mathcal{C}=\mu _\alpha \mathcal{P}^\alpha .  \tag{1.5.7}
\end{equation}
If $\theta $ is a section of $E^{*}$, $\theta =\theta _\alpha s^\alpha $,
and $X$ is a section of $E$, $X=X^\alpha s_\alpha $ then the vertical and
complete lifts have the expression
\[
\theta ^{\mathrm{v}}=\theta _\alpha \mathcal{P}^\alpha ,
\]
\[
X^{*\mathrm{c}}=X^\alpha \mathcal{Q}_\alpha -\left( \sigma _\alpha
^i\frac{
\partial X^\varepsilon }{\partial q^i}+L_{\alpha \beta }^\varepsilon X^\beta
\right) \mu _\varepsilon \mathcal{P}^\alpha .
\]

The Liouville-Hamilton section on $\mathcal{T}E^{*}$ measures the
homogeneity of the functions and sections. A function $f\in C^\infty (E^{*})$
is said to be homogeneous of degree $r\in \Bbb{Z}$ if
\[
\mathcal{L}_{\mathcal{C}}f=rf,
\]
where $\mathcal{L}_{\mathcal{C}}$ is the Lie derivation with
respect to the Liouville-Hamilton section on the Lie algebroid. A
section $\rho $ of $ \mathcal{T}E^{*}$ is said to be homogeneous
of degree $r\in \Bbb{Z}$ if
\[
\mathcal{L}_{\mathcal{C}}\rho =r\rho .
\]
We remark that $V\mathcal{T}E^{*}$ is Lagrangian for $\omega _E$,
i.e. $ \omega _E(\rho _1,\rho _2)=0$, for every vertical sections
$\rho _1,\rho _2\in \Gamma (V\mathcal{T}E^{*}).$\newpage\

\subsubsection{\textbf{Ehresmann connections on the Lie algebroid }$\mathcal{
T}E^{*}$}

\begin{definition}
The Ehresmann nonlinear connection on $\mathcal{T}E^{*}$ is an
almost product structure $\mathcal{N}$ on $\tau
_1:\mathcal{T}E^{*}\rightarrow E^{*} $(i.e. a bundle morphism
$\mathcal{N}:\mathcal{T}E^{*}\rightarrow \mathcal{T}E^{*}$, such
that $\mathcal{N}^2=Id$) smooth on $\mathcal{T} E^{*}\backslash
\{0\}$ such that $V\mathcal{T}E^{*}=\ker (Id+\mathcal{N})$.
\end{definition}

If $\mathcal{N}$ is a connection on $\mathcal{T}E^{*}$ then
$H\mathcal{T} E^{*}=\ker (Id-\mathcal{N})$ is the horizontal
distribution associated to $ \mathcal{N}$ and
\[
\mathcal{T}E^{*}=V\mathcal{T}E^{*}\oplus H\mathcal{T}E^{*}.
\]
Each $\rho \in \Gamma (\mathcal{T}E^{*})$ can be written as $\rho
=\rho ^h+\rho ^v$ where $\rho ^h$, $\rho ^v$ are sections in the
horizontal and respective, vertical subbundles. If $\rho ^h=0$
then $\rho $\ is called \textit{\ vertical }and if $\rho ^v=0$
then $\rho $ \ is called \textit{ horizontal}. \ A connection
$\mathcal{N}$ on $E^{*}$ induces two projectors $
h,v:\mathcal{T}E^{*}\rightarrow \mathcal{T}E^{*}$ such that
$h(\rho )=\rho ^h $ and $v(\rho )=\rho ^v$ for every $\rho \in
\Gamma (\mathcal{T}E^{*})$. We have
\[
h=\frac 12(Id+\mathcal{N}),\quad v=\frac 12(Id-\mathcal{N}),
\]
\[
\ker h=Imv=V\mathcal{T}E^{*},\quad Imh=\ker v=H\mathcal{T}E^{*},
\]
\[
h^2=h,\quad v^2=v,\quad hv=vh=0,\quad h+v=Id.
\]
Locally, a connection can be expressed as
\[
\mathcal{N}(\mathcal{Q}_\alpha )=\mathcal{Q}_\alpha
+2\mathcal{N}_{\alpha \beta }\mathcal{P}^\beta ,\
\mathcal{N}(\mathcal{P}^\alpha )=-\mathcal{P} ^\alpha ,
\]
where $\mathcal{N}_{\alpha \beta }=\mathcal{N}_{\alpha \beta }(q,\mu )$ are
\ the local coefficients of $\mathcal{N}$. The local coordinate expression
of $\mathcal{N}$ is
\[
\mathcal{N}=\mathcal{Q}_\alpha \otimes \mathcal{Q}^\alpha
-\mathcal{P} ^\alpha \otimes \mathcal{P}_\alpha
+2\mathcal{N}_{\alpha \beta }\mathcal{P} ^\alpha \otimes
\mathcal{Q}^\alpha .
\]
The local sections $\mathcal{P}^\alpha ,$ ($\alpha =\overline{1,m})$ define
a local frame of $V\mathcal{T}E^{*}$, and the sections
\begin{equation}
\delta _\alpha ^{*}=(\mathcal{Q}_\alpha )^h=\mathcal{Q}_\alpha
+\mathcal{N} _{\alpha \beta }\mathcal{P}^\beta ,  \tag{1.5.8}
\end{equation}
generate a local frame of $H\mathcal{T}E^{*}$. The frame $\{\delta _\alpha
^{*},\mathcal{P}^\alpha \}$ is a local basis of $\mathcal{T}E^{*}$ called%
\textit{\ adapted }to the direct sum decomposition. The dual
adapted basis is $\{\mathcal{Q}^\alpha ,\delta \mathcal{P}_\alpha
\}$ where$\ $
\begin{equation}
\delta \mathcal{P}_\alpha =\mathcal{P}_\alpha -\mathcal{N}_{\alpha
\beta } \mathcal{Q}^\beta .  \tag{1.5.9}
\end{equation}
It results that in the adapted basis the expression of Ehresmann connection
becomes
\[
\mathcal{N}=\delta _\alpha ^{*}\otimes \mathcal{Q}^\alpha
+\mathcal{P} ^\alpha \otimes \delta \mathcal{P}_\alpha ,
\]

\begin{definition}
A connection $\mathcal{N}$ is called symmetric if $H\mathcal{T}E^{*}$ is
Lagrangian for $\omega _E.${}
\end{definition}

By a straightforward computation, using (1.5.8) we get
\[
\omega _E(\delta _\alpha ^{*},\delta _\beta ^{*})=\mathcal{N}_{\alpha \beta
}-\mathcal{N}_{\beta \alpha }-\mu _\gamma L_{\alpha \beta }^\gamma ,
\]
and it results that $\mathcal{N}$ is symmetric if and only if
\begin{equation}
\mathcal{N}_{\alpha \beta }-\mathcal{N}_{\beta \alpha }=\mu _\gamma
L_{\alpha \beta }^\gamma .  \tag{1.5.10}
\end{equation}

\begin{proposition}
With respect to a symmetric nonlinear connection, the canonical symplectic
structure $\omega _E$ can be written in the following form
\begin{equation}
\omega _E=\mathcal{Q}^\alpha \wedge \delta \mathcal{P}_\alpha +\mu _\alpha
L_{\beta \gamma }^\alpha \mathcal{Q}^\beta \wedge \mathcal{Q}^\gamma .
\tag{1.5.11}
\end{equation}
\end{proposition}

\textbf{Proof}. Using (1.5.6) and (1.5.9) we get
\[
\omega _E=\mathcal{Q}^\alpha \wedge \delta \mathcal{P}_\alpha
+\frac 12( \mathcal{N}_{\alpha \beta }-\mathcal{N}_{\beta \alpha
})\mathcal{Q}^\alpha \wedge \mathcal{Q}^\beta +\frac 12\mu _\alpha
L_{\beta \gamma }^\alpha \mathcal{Q}^\beta \wedge
\mathcal{Q}^\gamma ,
\]
which ends the proof.\hfill\hbox{\rlap{$\sqcap$}$\sqcup$}

\begin{proposition}
The Lie brackets of the adapted basis $\{\delta _\alpha
^{*},\mathcal{P} ^\alpha \}$ are \cite{Hr4}
\[
[\delta _\alpha ^{*},\delta _\beta
^{*}]_{\mathcal{T}E^{*}}=L_{\alpha \beta }^\gamma \delta _\gamma
^{*}+\mathcal{R}_{\alpha \beta \gamma }\mathcal{P} ^\gamma ,\quad
[\delta _\alpha ^{*},\mathcal{P}^\beta ]_{\mathcal{T}E^{*}}=-
\frac{\partial \mathcal{N}_{\alpha \gamma }}{\partial \mu _\beta
}\mathcal{P} ^\gamma ,\quad [\mathcal{P}^\alpha ,\mathcal{P}^\beta
]_{\mathcal{T} E^{*}}=0,
\]
\begin{equation}
\mathcal{R}_{\alpha \beta \gamma }=\sigma _\alpha ^i\frac{\partial
\mathcal{N }_{\beta \gamma }}{\partial q^i}-\sigma _\beta
^i\frac{\partial \mathcal{N} _{\alpha \gamma }}{\partial
q^i}+\mathcal{N}_{\alpha \delta }\frac{\partial \mathcal{N}_{\beta
\gamma }}{\partial \mu _\delta }\mathcal{-N}_{\beta \delta
}\frac{\partial \mathcal{N}_{\alpha \gamma }}{\partial \mu _\delta
} -L_{\alpha \beta }^\varepsilon \mathcal{N}_{\varepsilon \gamma
}. \tag{1.5.12}
\end{equation}
\end{proposition}

{}\textbf{Proof.} Using (1.5.8) we obtain
\[
\lbrack \delta _\alpha ^{*},\delta _\beta
^{*}]_{\mathcal{T}E^{*}}=\left( \sigma _\alpha ^i\frac{\partial
\mathcal{N}_{\beta \gamma }}{\partial q^i} -\sigma _\beta
^i\frac{\partial \mathcal{N}_{\alpha \gamma }}{\partial q^i}+
\mathcal{N}_{\alpha \delta }\frac{\partial \mathcal{N}_{\beta
\gamma }}{
\partial \mu _\delta }\mathcal{-N}_{\beta \delta }\frac{\partial \mathcal{N}
_{\alpha \gamma }}{\partial \mu _\delta }\right) \mathcal{P}^\gamma
+L_{\alpha \beta }^\varepsilon \mathcal{Q}_\varepsilon ,
\]
and putting $\mathcal{Q}_\varepsilon =\delta _\varepsilon
^{*}-\mathcal{N} _{\varepsilon \gamma }\mathcal{P}^\gamma $ we get
$[\delta _\alpha ^{*},\delta _\beta
^{*}]_{\mathcal{T}E^{*}}=L_{\alpha \beta }^\gamma \delta _\gamma
^{*}+\mathcal{R}_{\alpha \beta \gamma }\mathcal{P}^\gamma .$\hfill
\hbox{\rlap{$\sqcap$}$\sqcup$}

The curvature of a connection $\mathcal{N}$ on $\mathcal{T}E^{*}$
is given by $\Omega =-\mathbf{N}_h,$ where $h$ is the horizontal
projector and $ \mathbf{N}_h$ is the Nijenhuis tensor of $h$,
given by
\[
\mathrm{N}_h(\theta ,\rho )=[h\theta ,h\rho ]_{\mathcal{T}E^{*}}-h[h\theta
,\rho ]_{\mathcal{T}E^{*}}-h[\theta ,h\rho ]_{\mathcal{T}E^{*}}+h^2[\theta
,\rho ]_{\mathcal{T}E^{*}}.
\]

\begin{proposition}
In local coordinates we get
\begin{equation}
\Omega =-\frac 12\mathcal{R}_{\alpha \beta \gamma }\mathcal{Q}^\alpha \wedge
\mathcal{Q}^\beta \otimes \mathcal{P}^\gamma ,  \tag{1.5.13}
\end{equation}
where $\mathcal{R}_{\alpha \beta \gamma }$ is given by (1.5.12) and are
called the coefficients of the \textit{curvature tensor} of $\mathcal{N}$.
\end{proposition}

\textbf{Proof}. Since $h^2=h$ we obtain
\[
\Omega (h\rho _1,h\rho _2)=-v[h\rho _1,h\rho _2]_{\mathcal{T}E^{*}},\quad
\Omega (h\rho _1,v\rho _2)=\Omega (v\rho _1,v\rho _2)=0,
\]
and in local coordinates we get
\[
\Omega (\delta _\alpha ^{*},\delta _\beta ^{*})=-v[\delta _\alpha
^{*},\delta _\beta ^{*}]_{\mathcal{T}E^{*}}=-\mathcal{R}_{\alpha \beta
\gamma }\mathcal{P}^\gamma ,
\]
which concludes the proof.\hfill
\hbox{\rlap{$\sqcap$}$\sqcup$}

\begin{proposition}
The curvature satisfies the Bianchi identity
\[
\mathcal{R}_{\alpha \beta \gamma }+\mathcal{R}_{\beta \gamma
\alpha }+ \mathcal{R}_{\gamma \alpha \beta }=0.
\]
\end{proposition}

\textbf{Proof}. By direct computation, using the relation (1.5.12) and
structure equations given by (1.2.9), (1.2.10).

The curvature is an obstruction to the integrability of $H\mathcal{T}E^{*}$,
understanding that a vanishing curvature entails that horizontal sections
are closed under the Lie algebroid bracket of $\mathcal{T}E^{*}.$ We have:

\begin{remark}
{}$H\mathcal{T}E^{*}$ is integrable if and only if the curvature vanishes.
\end{remark}

Also, the integrability conditions for the almost product
structure $ \mathcal{N}$ is given by the vanishing of the
associated Nijenhuis tensor $ \mathrm{N}_{\mathcal{N}}$.By a
straightforward computation we obtain
\[
\mathrm{N}_{\mathcal{N}}(\mathcal{P}^\alpha ,\mathcal{P}^\beta
)=0,\quad \mathrm{N}_{\mathcal{N}}(\delta _\alpha
^{*},\mathcal{P}^\beta )=0,\quad \mathrm{N}_{\mathcal{N}}(\delta
_\alpha ^{*},\delta _\beta ^{*})=4\mathcal{R} _{\alpha \beta
\gamma }\mathcal{P}^\gamma .
\]
Thus
\[
\mathrm{N}_{\mathcal{N}}=-2\mathcal{R}_{\alpha \beta \gamma
}\mathcal{X} ^\alpha \wedge \mathcal{X}^\beta \otimes
\mathcal{P}^\gamma ,
\]
and it results that the distribution $H\mathcal{T}E^{*}$ is integrable if
and only if the almost product structure $\mathcal{N}$ is integrable.

We consider the connections $\mathcal{N}$ on $\mathcal{T}E^{*}$ and \ $N$ on
$TE^{*}$ which are $\sigma ^1$-related and the adapted basis $(\delta
_i,\frac \partial {\partial \mu _\alpha })$ of $N$ given by $\delta _i=\frac
\partial {\partial q^i}+N_{i\alpha }\frac \partial {\partial \mu _\alpha }$,
where $N_{i\beta }$ are the coefficients of $N$.

\begin{theorem}
The following relations hold
\begin{equation}
\sigma ^1(\delta _\alpha ^{*})=\sigma _\alpha ^i\delta _i,\quad
\mathcal{N} _{\alpha \beta }=\sigma _\alpha ^iN_{i\beta },\quad
\mathcal{R}_{\alpha \beta \gamma }=\sigma _\alpha ^i\sigma _\beta
^jR_{ij\gamma },  \tag{1.5.14}
\end{equation}
where
\[
R_{ij\gamma }=\delta _i(N_{j\gamma })-\delta _j(N_{i\gamma }),
\]
is the curvature of the nonlinear connection $N$ on $TE^{*}.$
\end{theorem}

\textbf{Proof}\textit{.} Since $\mathcal{N}(\sigma ^1(\delta
_\alpha ^{*}))=\sigma ^1(\delta _\alpha ^{*})$ the relation
$N\circ \sigma ^1=\sigma ^1\circ \mathcal{N}$ leads to $N(\sigma
^1(\delta _\alpha ^{*}))=\sigma ^1(\delta _\alpha ^{*})$. But
$N(\sigma ^1(\delta _\alpha ^{*}))=2\sigma _\alpha ^i\delta
_i-\sigma ^1(\delta _\alpha ^{*})$ which concludes that $ \sigma
^1(\delta _\alpha ^{*})=\sigma _\alpha ^i\delta _i,$ from which we
easily obtain $\mathcal{N}_{\alpha \beta }=\sigma _\alpha
^iN_{i\beta }$. By straightforward computation we get
\[
\mathcal{R}_{\alpha \beta \gamma }=\sigma _\alpha ^i\sigma _\beta ^j\left(
\delta _i(N_{j\gamma })-\delta _j(N_{i\gamma })\right) +N_{j\gamma }\left(
\sigma _\beta ^i\frac{\partial \sigma _\alpha ^j}{\partial q^i}-\sigma
_\alpha ^i\frac{\partial \sigma _\beta ^j}{\partial q^i}\right) +L_{\alpha
\beta }^\varepsilon \mathcal{N}_{\varepsilon \gamma },
\]
and using (1.2.10), the second term is given by $N_{j\gamma }\sigma
_\varepsilon ^jL_{\beta \alpha }^\varepsilon =-N_{\varepsilon \gamma
}L_{\alpha \beta }^\varepsilon ,$ which ends the proof.\hfill
\hbox{\rlap{$\sqcap$}$\sqcup$}

\begin{remark}
A $\sigma ^1$-related connection $N$ on $TE^{*}$ determines a
connection $ \mathcal{N}$ on $\mathcal{T}E^{*}$ with the
coefficients given by
\[
\mathcal{N}_{\alpha \beta }=\sigma _\alpha ^iN_{i\beta },
\]
and curvature
\[
\mathcal{R}_{\alpha \beta \gamma }=\sigma _\alpha ^i\sigma _\beta
^jR_{ij\gamma }.
\]
Conversely, it is not true, because $\sigma $ is only injective.
\end{remark}

\newpage\

\subsubsection{\textbf{Regular sections and connections}}

\begin{definition}
An almost tangent structure $\mathcal{J}$ on $\mathcal{T}E^{*}$ is
a bundle morphism $\mathcal{J}:\mathcal{T}E^{*}\rightarrow
\mathcal{T}E^{*}$ of $\tau _1:\mathcal{T}E^{*}\rightarrow E^{*}$
of rank $m$, such that $\mathcal{J} ^2=0.$ An almost tangent
structure $\mathcal{J}$ on $\mathcal{T}E^{*}$ is called
\textit{adapted} if
\[
\ Im\mathcal{J}=\ker \mathcal{J}=V\mathcal{T}E^{*}.
\]
Locally, an adapted almost tangent structure is given by
\begin{equation}
\mathcal{J}=t_{\alpha \beta }\mathcal{Q}^\alpha \otimes \mathcal{P}^\beta ,
\tag{1.5.15}
\end{equation}
where the tensor $t_{\alpha \beta }(x,\mu )$ is nondegenerate.
\end{definition}

Now, we find the integrability conditions for the adapted almost tangent
structure. From \cite{Hr4} it results

\begin{proposition}
$\mathcal{J}$ is an integrable structure if and only if
\begin{equation}
\frac{\partial t^{\alpha \gamma }}{\partial \mu _\beta }=\frac{\partial
t^{\beta \gamma }}{\partial \mu _\alpha },  \tag{1.5.16}
\end{equation}
where $t^{\alpha \gamma }t_{\gamma \beta }=\delta _\beta ^\alpha $.
\end{proposition}

\textbf{Proof}. $\mathcal{J}$ is integrable if and only if the associated
Nijenhuis tensor
\[
\mathrm{N}_{\mathcal{J}}(\rho ,\nu )=[\mathcal{J}\rho
,\mathcal{J}\nu ]_{ \mathcal{T}E^{*}}-\mathcal{J}[\mathcal{J}\rho
,\nu ]_{\mathcal{T}E^{*}}- \mathcal{J}[\rho ,\mathcal{J}\nu
]_{\mathcal{T}E^{*}},
\]
vanishes. This is locally equivalent with the relations
\[
\mathrm{N}_{\mathcal{J}}(\mathcal{X}_\alpha ,\mathcal{X}_\beta
)=\left( t_{\alpha \gamma }\frac{\partial t_{\beta \varepsilon
}}{\partial \mu _\gamma }-t_{\beta \gamma }\frac{\partial
t_{\alpha \varepsilon }}{\partial \mu _\gamma }\right)
\mathcal{P}^\varepsilon ,\ \mathrm{N}_{\mathcal{J}}(
\mathcal{X}_\alpha ,\mathcal{P}^\beta
)=\mathrm{N}_{\mathcal{J}}(\mathcal{P} ^\alpha ,\mathcal{P}^\beta
)=0.
\]
Therefore $\mathcal{J}$ is integrable if and only if
\[
t_{\alpha \gamma }\frac{\partial t_{\beta \varepsilon }}{\partial \mu
_\gamma }=t_{\beta \gamma }\frac{\partial t_{\alpha \varepsilon }}{\partial
\mu _\gamma },
\]
that is equivalent to (1.5.16).

\begin{definition}
An adapted almost tangent structure $\mathcal{J}$ on $\mathcal{T}E^{*}$ is
called symmetric \ if
\begin{equation}
\omega _E(\mathcal{J}\rho _1,\rho _2)=\omega _E(\mathcal{J}\rho _2,\rho _1).
\tag{1.5.17}
\end{equation}
\end{definition}

Locally, this requires the symmetry of the tensor $t_{\alpha \beta }.$

\begin{remark}
If $g$ is a pseudo-Riemannian metric on the vertical bundle
$V\mathcal{T} E^{*}$ then there exists a unique symmetric adapted
almost tangent structure on $\mathcal{T}E^{*}$ such that
\begin{equation}
g(\mathcal{J}\rho ,\mathcal{J}\upsilon )=-\omega _E(\mathcal{J}\rho
,\upsilon ),  \tag{1.5.18}
\end{equation}
and we say that $\mathcal{J}$ is induced by the metric $g$.
\end{remark}

Locally, if
\[
g(q,\mu )=g^{\alpha \beta }\mathcal{P}_\alpha \otimes \mathcal{P}_\beta ,
\]
then the relation (1.5.18) implies $t^{\alpha \beta }=g^{\alpha \beta }$.

\begin{remark}
Any symmetric adapted almost tangent structure $\mathcal{J}$ on
$\mathcal{T} E^{*}$ induces a pseudo-Riemannian metric on the
vertical bundle $V\mathcal{T }E^{*}$ as defined by (1.5.18).
\end{remark}

\begin{definition}
The torsion of a connection $\mathcal{N}$ is the vector valued two
form $T=[ \mathcal{J},h],$ where $\ h$ is the horizontal projector
and $[\mathcal{J} ,h] $ is the Fr\"olicher-Nijenhuis bracket
\[
\begin{array}{c}
[\mathcal{J},h](X,Y)=[\mathcal{J}X,hY]_{\mathcal{T}E^{*}}+[hX,\mathcal{J}Y]_{
\mathcal{T}E^{*}}+\mathcal{J}[X,Y]_{\mathcal{T}E^{*}}-\mathcal{J}[X,hY]_{
\mathcal{T}E^{*}}- \\
-\mathcal{J}[hX,Y]_{\mathcal{T}E^{*}}-h[X,\mathcal{J}Y]_{\mathcal{T}E^{*}}-h[
\mathcal{J}X,Y]_{\mathcal{T}E^{*}},
\end{array}
\]
\end{definition}

The torsion $T$ is a semibasic vector-valued form. Its local expression is\
\[
T=\frac 12T_{\alpha \beta \gamma }\mathcal{Q}^\alpha \wedge
\mathcal{Q} ^\beta \otimes \mathcal{P}^\gamma ,
\]
where
\begin{equation}
T_{\alpha \beta \gamma }=t_{\alpha \epsilon }\frac{\partial
\mathcal{N} _{\beta \gamma }}{\partial \mu _\varepsilon }-t_{\beta
\epsilon }\frac{
\partial \mathcal{N}_{\alpha \gamma }}{\partial \mu _\varepsilon }+\delta
_\alpha ^{*}(t_{\beta \gamma })-\delta _\beta ^{*}(t_{\alpha \gamma
})-L_{\alpha \beta }^\epsilon t_{\epsilon \gamma }.  \tag{1.5.19}
\end{equation}
Next, let us consider the linear mapping $\mathcal{F}:\mathcal{T}
E^{*}\rightarrow \mathcal{T}E^{*}$ given by
\begin{equation}
\mathcal{F}(h\rho )=\mathcal{J}\rho ,\quad \mathcal{F}(\mathcal{J}\rho
)=-h\rho ,  \tag{1.5.20}
\end{equation}
where $\rho \in \Gamma (\mathcal{T}E^{*})$, and $h$ is the horizontal
projector induced by the nonlinear connection.

\begin{proposition}
The mapping $\mathcal{F}$ has the properties:

a) $\mathcal{F}$ is an almost complex structure $\mathcal{F}\circ
\mathcal{F} =-Id.$

b) Locally it is given by
\begin{equation}
\mathcal{F}=t_{\alpha \beta }\mathcal{P}^\beta \otimes \mathcal{Q}^\alpha
-t^{\alpha \beta }\delta _\alpha ^{*}\otimes \delta \mathcal{P}_\beta .
\tag{1.5.21}
\end{equation}
\end{proposition}

\textbf{Proof}. It results by definition that
\[
(\mathcal{F}\circ \mathcal{F})(h\rho
)=\mathcal{F}(\mathcal{J}(\rho ))=-h\rho ,\quad (\mathcal{F}\circ
\mathcal{F})(\mathcal{J}\rho )=\mathcal{F} (-h\rho
)=-\mathcal{J}\rho ,
\]
\[
\mathcal{F}(\delta _\alpha ^{*})=t_{\alpha \beta }\mathcal{P}^\beta ,\quad
\mathcal{F}(\mathcal{P}^\alpha )=-t^{\alpha \beta }\delta _\beta ^{*},
\]
which concludes the proof.\hfill
\hbox{\rlap{$\sqcap$}$\sqcup$}

\begin{proposition}
The almost complex structure is integrable if and only if the torsion and
curvature of the connection satisfy the equations
\begin{equation}
T_{\alpha \beta \gamma }=0,\quad \mathcal{R}_{\alpha \beta \gamma
}=t_{\alpha \varepsilon }\frac{\partial t_{\beta \gamma
}}{\partial \mu _\varepsilon }-t_{\beta \varepsilon
}\frac{\partial t_{\alpha \gamma }}{
\partial \mu _\varepsilon }.  \tag{1.5.22}
\end{equation}
\end{proposition}

\textbf{Proof}. Let $\mathrm{N}_{\mathcal{F}}$ be the Nijenjuis tensor of
the almost complex structure. We get
\[
\mathrm{N}_{\mathcal{F}}(\delta _\alpha ^{*},\delta _\beta
^{*})=\mathrm{N} _{\alpha \beta }^\gamma \delta _\gamma
^{*}+\mathrm{N}_{\alpha \beta (\gamma )}\mathcal{P}^\gamma ,
\]
\[
\mathrm{N}_{\mathcal{F}}(\delta _\alpha ^{*},\mathcal{P}^\beta
)=\mathrm{N} _\alpha ^{(\beta )\gamma }\delta _\gamma
^{*}+\mathrm{N}_{\alpha (\gamma )}^{(\beta )}\mathcal{P}^\gamma ,
\]
\[
\mathrm{N}_{\mathcal{F}}(\mathcal{P}^\alpha ,\mathcal{P}^\beta
)=-t^{\varepsilon \alpha }t^{\beta \tau }\mathrm{N}_{\mathcal{F}}(\delta
_\varepsilon ^{*},\delta _\tau ^{*}),
\]
where
\[
\mathrm{N}_{\alpha \beta }^\gamma =T_{\alpha \beta \varepsilon
}t^{\varepsilon \gamma },\quad \mathrm{N}_{\alpha \beta (\gamma
)}=t_{\alpha \varepsilon }\frac{\partial t_{\beta \gamma
}}{\partial \mu _\varepsilon } -t_{\beta \varepsilon
}\frac{\partial t_{\alpha \gamma }}{\partial \mu _\varepsilon
}-\mathcal{R}_{\alpha \beta \gamma },
\]
\[
\mathrm{N}_{\alpha \beta (\gamma )}=\mathrm{N}_\alpha ^{(\varepsilon )\tau
}t_{\varepsilon \beta }t_{\gamma \tau },\quad \mathrm{N}_{\alpha \beta
}^\gamma =\mathrm{N}_{\alpha (\tau )}^{(\varepsilon )}t^{\tau \gamma
}t_{\beta \varepsilon },
\]
which ends the proof.\hfill\hbox{\rlap{$\sqcap$}$\sqcup$}

\begin{remark}
{}Let $\mathcal{N}$ be a bundle morphism of $\tau _1:\mathcal{T}
E^{*}\rightarrow E^{*}$, smooth on $\mathcal{T}E^{*}\backslash
\{0\}$. Then $ \mathcal{N}$ is a connection on $\mathcal{T}E^{*}$
if and only if there exists an adapted almost tangent structure
$\mathcal{J}$ on $\mathcal{T} E^{*} $ such that
\begin{equation}
\mathcal{JN}=\mathcal{J},\quad \mathcal{NJ}=-\mathcal{J}.  \tag{1.5.23}
\end{equation}
\end{remark}

\begin{definition}
Let $\mathcal{J}$ be an adapted tangent structure on
$\mathcal{T}E^{*}$. A section $\rho $ of $\mathcal{T}E^{*}$ is
called $\mathcal{J}-$\textit{regular } if
\begin{equation}
\mathcal{J}[\rho ,\mathcal{J}\nu ]_{\mathcal{T}E^{*}}=-\mathcal{J}\nu ,
\tag{1.5.24}
\end{equation}
for every section $\nu $ of $\mathcal{T}E^{*}$.
\end{definition}

Locally, the section $\rho =\xi ^\alpha \mathcal{Q}_\alpha +\rho _\beta
\mathcal{P}^\beta $\ is $\mathcal{J}-$regular if and only if
\[
t^{\alpha \beta }=\frac{\partial \xi ^\beta }{\partial \mu _\alpha }.
\]
where $t^{\alpha \beta }t_{\alpha \gamma }=\delta _\gamma ^\beta .$

We have to remark that if the equation (1.5.24) is satisfied for any section
$\nu $ $\in \Gamma (\mathcal{T}E^{*})$ with $rank[t^{\alpha \beta }]=m,$
then $\mathcal{J}$ is an \textit{integrable} structure. Indeed, we have
\[
\frac{\partial t^{\alpha \beta }}{\partial \mu _\gamma }=\frac{\partial
^2\xi ^\beta }{\partial \mu _\gamma \partial \mu _\alpha }=\frac{\partial
^2\xi ^\beta }{\partial \mu _\alpha \partial \mu _\gamma }=\frac{\partial
t^{\gamma \beta }}{\partial \mu _\alpha },
\]
and using (1.5.16) it follows that $\mathcal{J}$ is integrable.

From \cite{Hr4} we have:

\begin{theorem}
{}Let $\mathcal{J}$ be an adapted tangent structure on $\mathcal{T}E^{*}.$
If $\rho $ is a $\mathcal{J}-$regular section of $\mathcal{T}E^{*}$ then
\begin{equation}
\mathcal{N}=-\mathcal{L}_\rho \mathcal{J},  \tag{1.5.25}
\end{equation}
is a connection on $\mathcal{T}E^{*}.$
\end{theorem}

\textbf{Proof}. Since
\[
\mathcal{N}(\upsilon )=-\mathcal{L}_\rho \mathcal{J}(\upsilon
)=-[\rho , \mathcal{J}\upsilon
]_{\mathcal{T}E^{*}}+\mathcal{J}[\rho ,\upsilon ]_{
\mathcal{T}E^{*}},
\]
then
\[
\mathcal{JN(}\upsilon )=-\mathcal{J}[\rho ,\mathcal{J}\upsilon
]_{\mathcal{T} E^{*}}+\mathcal{J}^2[\rho ,\upsilon
]_{\mathcal{T}E^{*}}=\mathcal{J}\upsilon ,
\]
\[
\mathcal{NJ(}\upsilon )=-[\rho ,\mathcal{J}^2\upsilon
]_{\mathcal{T}E^{*}}+ \mathcal{J}[\rho ,\mathcal{J}\upsilon
]_{\mathcal{T}E^{*}}=-\mathcal{J} \upsilon .
\]
and using (1.5.23) it results the
conclusion.\hfill\hbox{\rlap{$\sqcap$}$ \sqcup$}

This connection is induced by $\mathcal{J}$ and $\rho $. Its local
coefficients are given by
\begin{equation}
\mathcal{N}_{\alpha \beta }=\frac 12\left( t_{\alpha \gamma
}\frac{\partial \rho _\beta }{\partial \mu _\gamma }-\sigma
_\alpha ^it_{\gamma \beta }\frac{
\partial \xi ^\gamma }{\partial q^i}-\rho (t_{\alpha \beta })+\xi ^\gamma
t_{\varepsilon \beta }L_{\gamma \alpha }^\varepsilon \right) .  \tag{1.5.26}
\end{equation}

\begin{definition}
An adapted tangent structure $\mathcal{J}$ on $\mathcal{T}E^{*}$ is called
homogeneous if
\[
\mathcal{L}_{\mathcal{C}}\mathcal{J}=-\mathcal{J}.
\]
\end{definition}

Notice that $\mathcal{J}$ is homogeneous if the local components $t_{\alpha
\beta }(x,\mu )$ are $0$-homogeneous with respect to $\mu $.

\begin{proposition}
Let $\mathcal{J}$ be a homogeneous adapted tangent structure. A
section $ \rho $ of $\mathcal{T}E^{*}$ is $\mathcal{J}-$regular if
and only if
\[
\mathcal{J}\rho =\mathcal{C}.
\]
\end{proposition}

\textbf{Proof}. If $\rho $ is $\mathcal{J}-$regular then
\[
t^{\alpha \beta }=\frac{\partial \xi ^\beta }{\partial \mu _\alpha },
\]
is $0$-homogeneous, hence $\xi ^\beta $ must be $1$-homogeneous
with respect to $\mu $, therefore $\xi ^\beta =\mu _\alpha
t^{\alpha \beta }$, that is equivalent to $\mathcal{J}\rho
=\mathcal{C}$. Vice versa, if $\mathcal{J} \rho =\mathcal{C}$ then
$\xi ^\beta =\mu _\alpha t^{\alpha \beta }$ and thus
\[
\frac{\partial \xi ^\beta }{\partial \mu _\gamma }=t^{\gamma \beta }+\mu
_\varepsilon \frac{t^{\gamma \beta }}{\partial \mu _\varepsilon }=t^{\gamma
\beta },
\]
which ends the proof. \hfill\hbox{\rlap{$\sqcap$}$\sqcup$}

\begin{remark}
(i) Based on the above result, the local expression for a
$\mathcal{J}-$ regular section with $\mathcal{J}$ a homogeneous
adapted tangent structure is
\begin{equation}
\rho =\mu _\alpha t^{\alpha \beta }\mathcal{X}_\beta +\rho _\alpha
\mathcal{P }^\alpha .  \tag{1.5.27}
\end{equation}

(ii) The coefficients (1.5.26) generated by $\rho $ from (1.5.27) can be
written in the following form
\[
\mathcal{N}_{\alpha \beta }=\frac 12\left( \left( \sigma _\alpha
^i\frac{
\partial t_{\gamma \beta }}{\partial x^i}-\sigma _\gamma ^i\frac{\partial
t_{\alpha \beta }}{\partial x^i}\right) t^{\varepsilon \gamma }\mu
_\varepsilon +t_{\alpha \varepsilon }\frac{\partial \rho _\beta
}{\partial \mu _\varepsilon }-\rho _\varepsilon \frac{\partial
t_{\alpha \beta }}{
\partial \mu _\varepsilon }+\mu _\varepsilon t^{\varepsilon \gamma
}L_{\gamma \alpha }^\lambda t_{\lambda \beta }\right) .
\]
\end{remark}

\newpage\

\subsubsection{\textbf{Hamilton sections}}

\begin{definition}
A section $\psi $ on $\mathcal{T}E^{*}$ is called a \textit{Hamilton section
}if it is $\mathcal{J}-$regular and
\[
\mathcal{L}_\psi \omega _E=0,
\]
where $\omega _E$ is the canonical symplectic section.
\end{definition}

If $\psi =\xi ^\alpha \mathcal{Q}_\alpha +\rho _\alpha \mathcal{P}^\alpha $
then the condition $\mathcal{L}_\psi \omega _E=0$ is expressed locally by
\cite{Po25}
\[
i)\quad \frac{\partial \xi ^\beta }{\partial \mu _\alpha }=\frac{\partial
\xi ^\alpha }{\partial \mu _\beta }.
\]
\begin{equation}
ii)\quad \sigma _\alpha ^i\frac{\partial \xi ^\beta }{\partial
q^i}+\frac{
\partial \rho _\alpha }{\partial \mu _\beta }=\xi ^\gamma L_{\gamma \alpha
}^\beta +\mu _\varepsilon L_{\gamma \alpha }^\varepsilon \frac{\partial \xi
^\gamma }{\partial \mu _\beta }.  \tag{1.5.28}
\end{equation}
\begin{eqnarray*}
iii)\quad \sigma _\beta ^i\frac{\partial \rho _\alpha }{\partial q^i}-\sigma
_\alpha ^i\frac{\partial \rho _\beta }{\partial q^i} &=&\mu _\varepsilon \xi
^\gamma \left( \sigma _\beta ^i\frac{\partial L_{\gamma \alpha }^\varepsilon
}{\partial q^i}-\sigma _\alpha ^i\frac{\partial L_{\gamma \beta
}^\varepsilon }{\partial q^i}+L_{\nu \gamma }^\varepsilon L_{\alpha \beta
}^\nu \right) + \\
&& \\
&&\ \ +\mu _\varepsilon \frac{\partial \xi ^\gamma }{\partial q^i}\left(
\sigma _\beta ^iL_{\gamma \alpha }^\varepsilon -\sigma _\alpha ^iL_{\gamma
\beta }^\varepsilon \right) -\rho _\gamma L_{\alpha \beta }^\gamma .
\end{eqnarray*}

Now, we deal with some generalizations of the Hamilton sections.

\begin{definition}
The section $\psi =\xi ^\alpha \mathcal{Q}_\alpha +\rho _\alpha
\mathcal{P} ^\alpha $ on $\mathcal{T}E^{*}$ defines a mechanical
structure if $\psi $ is $\mathcal{J}-$regular and
\begin{equation}
\omega _E(\mathcal{J}\rho _1,\rho _2)=\omega _E(\mathcal{J}\rho _2,\rho _1),
\tag{1.5.29}
\end{equation}
for any $\rho _1,\rho _2\in \Gamma (\mathcal{T}E^{*}).$
\end{definition}

The definition is equivalent with the symmetry of $t^{\alpha \beta
}=\frac{
\partial \xi ^\beta }{\partial \mu _\alpha },$ which means that the property
(1.5.28) $i)$ is fulfilled.

\begin{proposition}
The section $\psi $ on $\mathcal{T}E^{*}$ defines a mechanical structure if
and only if
\[
\mathcal{L}_\psi \omega _E(\nu _1,\nu _2)=0,
\]
whenever $\nu _1,\nu _2\in \Gamma (V\mathcal{T}E^{*}).$
\end{proposition}

\textbf{Proof}\textit{.} Let us consider $\rho _1,\rho _2\in
\Gamma ( \mathcal{T}E^{*})$ and $\nu _1=\mathcal{J}\rho _1,$ $\nu
_2=\mathcal{J}\rho _2$. Then

\begin{eqnarray*}
\mathcal{L}_\psi \omega _E(\mathcal{J}\rho _1,\mathcal{J}\rho _2)
&=&\psi (\omega _E(\mathcal{J}\rho _1,\mathcal{J}\rho _2)-\omega
_E((\mathcal{L}
_\psi \mathcal{J}\rho _1,\mathcal{J}\rho _2)- \\
-\omega _E((\mathcal{J}\rho _1,\mathcal{L}_\psi \mathcal{J}\rho
_2) &=&\omega _E((\mathcal{N}\rho _1,\mathcal{J}\rho _2)+\omega
_E((\mathcal{J} \rho _1,\mathcal{N}\rho _2).
\end{eqnarray*}
Considering $\mathcal{N}\rho _1$, $\mathcal{N}\rho _2$ for $\rho
_1,\rho _2$ in the previous relation and using the properties
(1.5.23) and $\mathcal{N} ^2=Id$ we obtain
\[
\mathcal{L}_\psi \omega _E(\mathcal{J}\rho _1,\mathcal{J}\rho
_2)=\mathcal{L} _\psi \omega _E(\mathcal{JN}\rho
_1,\mathcal{JN}\rho _2)=\omega _E(\rho _1, \mathcal{J}\rho
_2)+\omega _E(\mathcal{J}\rho _1,\rho _2),
\]
which ends the proof. \hfill
\hbox{\rlap{$\sqcap$}$\sqcup$}

\begin{definition}
The section $\psi $ on $\mathcal{T}E^{*}$ is a semi-Hamilton
section if $ \psi $ is $\mathcal{J}-$regular and
\begin{equation}
i_\upsilon (\mathcal{L}_\psi \omega _E)=0,  \tag{1.5.30}
\end{equation}
whenever $\upsilon \in \Gamma (V\mathcal{T}E^{*}).$
\end{definition}

In the case of a semi-Hamilton section on $\mathcal{T}E^{*}$ only the
conditions (1.5.28), $i)$ and $ii)$ are satisfied.

Let us consider a section $\psi $ defining a mechanical structure
on $ \mathcal{T}E^{*}$, then we can find the other connection.
Indeed, we consider
\[
g(\rho _1,\rho _2)=-\omega _E(\mathcal{J}\rho _1,\rho _2),
\]
for $\rho _1,\rho _2\in \Gamma (\mathcal{T}E^{*}).$ It results
that $g(\rho _1,\rho _2)=g(\rho _2,\rho _1)$ and $g(\upsilon ,\rho
)=0$ whenever $ \upsilon \in \Gamma (V\mathcal{T}E^{*}).$ The
local coordinate expression of $g$ is given by
\[
g=g_{\alpha \beta }\mathcal{Q}^\alpha \otimes \mathcal{Q}^\beta .
\]

\begin{proposition}
If $\psi $ defines a mechanical structure on $\mathcal{T}E^{*}$ then the
section $\mathcal{N}^{\prime }$ defined by
\begin{equation}
\omega _E(\mathcal{N}^{\prime }\rho _1,\rho _2)=(\mathcal{L}_\psi g)(\rho
_1,\rho _2),\quad \rho _1,\rho _2\in \Gamma (\mathcal{T}E^{*}),  \tag{1.5.31}
\end{equation}
determines a connection on $\mathcal{T}E^{*}.$
\end{proposition}

\textbf{Proof}. First, we show that the following relation holds
\[
(\mathcal{L}_\psi g)(\rho _1,\rho _2)=-(\mathcal{L}_\psi \omega
_E)(\mathcal{ J}\rho _1,\rho _2)+\omega _E(\mathcal{N}\rho _1,\rho
_2).
\]
Indeed, we get
\[
(\mathcal{L}_\psi g)(\rho _1,\rho _2)=\psi (g(\rho _1,\rho
_2))-g(\mathcal{L} _\psi \rho _1,\rho _2)-g(\rho
_1,\mathcal{L}_\psi \rho _2)=
\]
\[
=-\psi (\omega _E(\mathcal{J}\rho _1,\rho _2))+\omega
_E(\mathcal{J}( \mathcal{L}_\psi \rho _1),\rho _2)+\omega
_E(\mathcal{J}\rho _1,\mathcal{L} _\psi \rho _2).
\]
But the relation
\[
\psi (\omega _E(\mathcal{J}\rho _1,\rho _2))=(\mathcal{L}_\psi
\omega _E)( \mathcal{J}\rho _1,\rho _2)+\omega _E(\mathcal{L}_\psi
(\mathcal{J}\rho _1),\rho _2)+\omega _E(\mathcal{J}\rho
_1,\mathcal{L}_\psi \rho _2),
\]
yields
\[
(\mathcal{L}_\psi g)(\rho _1,\rho _2)=-(\mathcal{L}_\psi \omega
_E)(\mathcal{ J}\rho _1,\rho _2)+\omega _E(-(\mathcal{L}_\psi
\mathcal{J})\rho _1,\rho _2).
\]
Also, it results
\begin{equation}
(\mathcal{L}_\psi \omega _E)(\mathcal{J}\rho _1,\rho _2)=\omega
_E((\mathcal{ N}-\mathcal{N}^{\prime })\rho _1,\rho _2).
\tag{1.5.32}
\end{equation}
But,
\[
\omega _E(\mathcal{N}^{\prime }\mathcal{J}\rho _1,\rho _2)=(\mathcal{L}_\psi
g)(\mathcal{J}\rho _1,\rho _2)=
\]
\[
=\psi (g(\mathcal{J}\rho _1,\rho _2))-g((\mathcal{L}_\psi \mathcal{J})\rho
_1,\rho _2)-g(\mathcal{J}\rho _1,\mathcal{L}_\psi \rho _2)=
\]
\[
=g(\mathcal{N}\rho _1,\rho _2)=-\omega _E(\mathcal{J}\mathcal{N}\rho _1,\rho
_2)=-\omega _E(\mathcal{J}\rho _1,\rho _2),
\]
whence $\mathcal{N}^{\prime }\mathcal{J}=-\mathcal{J}$. Next
\[
\omega _E(\mathcal{J}\mathcal{N}^{\prime }\rho _1,\rho _2)=-\omega
_E( \mathcal{N}^{\prime }\rho _1,\mathcal{J}\rho
_2)=-(\mathcal{L}_\psi g)(\rho _1,\mathcal{J}\rho _2)=
\]
\[
(\mathcal{L}_\psi \omega _E)(\mathcal{J}\rho _1,\mathcal{J}\rho
_2)-\omega _E(\mathcal{N}\rho _1,\mathcal{J}\rho _2)=\omega
_E(\mathcal{J}\mathcal{N} \rho _1,\rho _2)=\omega
_E(\mathcal{J}\rho _1,\rho _2),
\]
whence $\mathcal{JN}^{\prime }=\mathcal{J},$ which ends the proof. \hfill
\hbox{\rlap{$\sqcap$}$\sqcup$}

In local coordinates the connection $\mathcal{N}^{\prime }$ has the
coefficients given by
\[
\mathcal{N}_{\alpha \beta }^{\prime }=\frac 12\left( -\psi
t_{\alpha \beta }-\sigma _\beta ^i\frac{\partial \xi ^\varepsilon
}{\partial q^i} t_{\varepsilon \alpha }-\sigma _\alpha
^i\frac{\partial \xi ^\varepsilon }{
\partial q^i}t_{\varepsilon \beta }+\xi ^\gamma L_{\gamma \beta
}^\varepsilon t_{\varepsilon \alpha }+\xi ^\gamma L_{\gamma \alpha
}^\varepsilon t_{\varepsilon \beta }+\mu _\varepsilon L_{\beta \alpha
}^\varepsilon \right) .
\]
Using (1.5.32) we obtain

\begin{remark}
In the case that $\psi $ is both semi-Hamiltonian and mechanical
sections and moreover Hamiltonian, then the connections
$\mathcal{N}$ and $\mathcal{N} ^{\prime }$ coincide.
\end{remark}

\newpage\

\subsubsection{\textbf{Hamiltonian formalism on Lie algebroids}}

Let us consider a differentiable and regular Hamiltonian
$\mathcal{H}$ $:$ $ E^{*}$ $\rightarrow \Bbb{R}$ i.e. the matrix
\[
g^{\alpha \beta }(q,\mu )=\frac{\partial ^2\mathcal{H}}{\partial \mu _\alpha
\partial \mu _\beta },
\]
is nondegenerate.

Any regular Hamiltonian $\mathcal{H}$\ on $E^{*}$ induces a
pseudo-Riemannian metric on $V\mathcal{T}E^{*}$ (the metric tensor
is $ g^{\alpha \beta }(q,\mu $)) therefore, it induces a unique
symmetric adapted \ almost tangent structure (denoted
$\mathcal{J}_{\mathcal{H}}$) such that (1.5.18) is verified.
Moreover, this is a tangent structure i.e., $\mathcal{J
}_{\mathcal{H}}$ is integrable.

A $\mathcal{J}-$regular section induced by the regular Hamiltonian
$\mathcal{ H}$ is
\[
\rho =\frac{\partial \mathcal{H}}{\partial \mu _\alpha }\mathcal{Q}_\alpha
+\rho _\alpha \mathcal{P}^\alpha .
\]
Since $\omega _E$ is a symplectic section on the Lie algebroid
$(\mathcal{T} E^{*},[\cdot ,\cdot ]_{\mathcal{T}E^{*}},\sigma ^1)$
and $d\mathcal{H}\in \Gamma (\mathcal{T}E^{*})$, we get

\begin{remark}
There exists a unique section $\rho _{\mathcal{H}}\in \Gamma
(\mathcal{T} E^{*})$ such that
\[
i_{\rho _{\mathcal{H}}}\omega _E=d^E\mathcal{H},
\]
\end{remark}

and $\rho _{\mathcal{H}}$ is a \textit{Hamilton section, i.e} the relations
(1.5.28) are fulfilled.

With respect to the local basis $\{\mathcal{Q}_\alpha
,\mathcal{P}^\alpha \}$ , the local expression of $\rho
_{\mathcal{H}}$ is
\begin{equation}
\rho _{\mathcal{H}}=\frac{\partial \mathcal{H}}{\partial \mu
_\alpha } \mathcal{Q}_\alpha -\left( \sigma _\alpha
^i\frac{\partial \mathcal{H}}{
\partial q^i}+\mu _\gamma L_{\alpha \beta }^\gamma \frac{\partial \mathcal{H}
}{\partial \mu _\beta }\right) \mathcal{P}^\alpha ,  \tag{1.5.33}
\end{equation}
Thus, the vector field $\sigma ^1(\rho _{\mathcal{H}})$ on $E^{*}$ is given
by \cite{Le}
\[
\sigma ^1(\rho _{\mathcal{H}})=\sigma _\alpha ^i\frac{\partial
\mathcal{H}}{
\partial \mu _\alpha }\frac \partial {\partial q^i}-\left( \sigma _\alpha ^i
\frac{\partial \mathcal{H}}{\partial q^i}+\mu _\gamma L_{\alpha \beta
}^\gamma \frac{\partial \mathcal{H}}{\partial \mu _\beta }\right) \frac
\partial {\partial \mu _\alpha },
\]
and consequently, the integral curves of $\rho _{\mathcal{H}}$ (i.e. the
integral curves of the vector field $\sigma ^1(\rho _{\mathcal{H}})$)
satisfy the Hamilton equations
\begin{equation}
\frac{dq^i}{dt}=\sigma _\alpha ^i\frac{\partial
\mathcal{H}}{\partial \mu _\alpha },\quad \frac{d\mu _\alpha
}{dt}=-\sigma _\alpha ^i\frac{\partial \mathcal{H}}{\partial
q^i}-\mu _\gamma L_{\alpha \beta }^\gamma \frac{
\partial \mathcal{H}}{\partial \mu _\beta }.  \tag{1.5.34}
\end{equation}
The \textit{Theorem 1.5.2} yields \cite{Po19}
\begin{corollary}
The symmetric nonlinear connection $\mathcal{N}=-\mathcal{L}_{\rho
_{ \mathcal{H}}}\mathcal{J}_{\mathcal{H}}$ has the coefficients
given by
\begin{equation*}
\mathcal{N}_{\alpha \beta }=\frac 12(\sigma _\gamma ^i\{g_{\alpha
\beta }, \mathcal{H}\}-\frac{\partial ^2\mathcal{H}}{\partial
q^i\partial \mu _\varepsilon }(\sigma _\beta ^ig_{\alpha
\varepsilon }+\sigma _\alpha
^ig_{\beta \varepsilon })+
\end{equation*}
\begin{equation}
 +\ \mu _\gamma L_{\varepsilon \kappa }^\gamma \frac{\partial
\mathcal{H }}{\partial \mu _\varepsilon }\frac{\partial g_{\alpha
\beta }}{\partial \mu _\kappa }+\mu _\gamma L_{\alpha \beta
}^\gamma +\frac{\partial \mathcal{H}}{
\partial \mu _\delta }(g_{\alpha \varepsilon }L_{\delta \beta }^\varepsilon
+g_{\beta \varepsilon }L_{\delta \alpha }^\varepsilon )),
\tag{1.5.35}
\end{equation}
where
\[
\{g_{\alpha \beta },\mathcal{H}\}=\frac{\partial g_{\alpha \beta }}{\partial
\mu _\gamma }\frac{\partial \mathcal{H}}{\partial q^i}-\frac{\partial
g_{\alpha \beta }}{\partial q^i}\frac{\partial \mathcal{H}}{\partial \mu
_\gamma },
\]
is the Poisson bracket.
\end{corollary}

For the particular case of the cotangent bundle see \cite{Mi2,
Op3, Hr1}
\newpage\

\subsubsection{\textbf{Duality between the Lagrangian and Hamiltonian
formalism}}

We have a local diffeomorphism $\Phi $ from $E^{*}$ to $E,$ locally given by
\begin{equation}
x^i=q^i,\quad y^\alpha =\xi ^\alpha (q,\mu )=\frac{\partial
\mathcal{H}}{
\partial \mu _\alpha }.  \tag{1.5.36}
\end{equation}
and its inverse $\Phi ^{-1}$ has the following local coordinates expression
\begin{equation}
q^i=x^i,\quad \mu _\alpha =\zeta _\alpha (x,y)=\frac{\partial
\mathcal{L}}{
\partial y^\alpha },  \tag{1.5.37}
\end{equation}
where
\[
\mathcal{L}(x,y)=\zeta _\alpha y^\alpha -\mathcal{H}(x,\mu ),
\]
and the components $\zeta _\alpha (x,y)$ define an 1-section on $E$.

From the condition for $\Phi ^{-1}$ to be the inverse of $\Phi $ we get the
following formulas
\[
\mathcal{V}_\beta (\zeta _\alpha )\circ \Phi =g_{\alpha \beta },
\]
\[
\mathcal{X}_\beta (\zeta _\alpha )\circ \Phi =-g_{\alpha \gamma
}\mathcal{Q} _\beta (\xi ^\gamma ),
\]
\begin{equation}
\Phi _{*}\mathcal{P}^\alpha =(g^{\alpha \beta }\circ \Phi
^{-1})\mathcal{V} _\beta ,  \tag{1.5.38}
\end{equation}
\[
\Phi _{*}(\mathcal{Q}_\alpha )=\mathcal{X}_\alpha +(\mathcal{Q}_\alpha (\xi
^\beta )\circ \Phi ^{-1})\mathcal{V}_\beta ,
\]
\[
\Phi _{*}^{-1}(\mathcal{V}_\alpha )=g_{\alpha \beta }\mathcal{P}^\beta ,
\]
\[
\Phi _{*}^{-1}(\mathcal{X}_\alpha )=\mathcal{Q}_\alpha -g_{\gamma
\varepsilon }\mathcal{Q}_\alpha (\xi ^\varepsilon )\mathcal{P}^\gamma ,
\]
where $\Phi _{*}$ is tangent map of $\Phi $ (see \cite{Le}) and $g^{\alpha
\beta }=\partial \xi ^\alpha /\partial \mu _\beta ,$ $g^{\alpha \beta
}g_{\beta \gamma }=\delta _\gamma ^\alpha $.

\begin{theorem}
Let $\rho $ be a $\mathcal{J}-$regular section on $\mathcal{T}E^{*}$ induced
by the hyperregular Hamiltonian $\mathcal{H}$, and $\Phi :E^{*}\rightarrow E$
the global diffeomorpfism given by (1.5.36), then the section $\Phi _{*}\rho
$ is a semispray on $\mathcal{T}E$ whose induced connection $N$ is the image
by $\Phi $ of the connection $\mathcal{N}$ defined by $\rho $.
\end{theorem}

\textbf{Proof.} Let us consider $\rho =\xi ^\alpha \mathcal{Q}_\alpha +\rho
_\alpha \mathcal{P}^\alpha $ and we obtain the semispray $\vartheta $
\[
\vartheta =\Phi _{*}\rho =\xi ^\alpha (\mathcal{X}_\alpha
+(\mathcal{Q} _\alpha (\xi ^\beta )\circ \Phi
^{-1})\mathcal{V}_\beta )+\rho _\alpha (g^{\alpha \beta }\circ
\Phi ^{-1})\mathcal{V}_\beta =y^\alpha \mathcal{X} _\alpha
+\vartheta ^\alpha \mathcal{V}_\alpha ,
\]
where $\vartheta ^\alpha $ is given by
\[
\vartheta ^\alpha \circ \Phi =\xi ^\beta \mathcal{Q}_\beta (\xi ^\alpha
)+\rho _\beta \mathcal{P}^\beta (\xi ^\alpha )=\rho (\xi ^\alpha ).
\]
Considering $\widetilde{\Phi }$ the map induced by $\Phi $ at the level of
tensors and using (1.5.38) we obtain
\[
\widetilde{\Phi }\mathcal{J}=J,
\]
which yields
\[
N=-\mathcal{L}_\vartheta J=-\mathcal{L}_{\Phi _{*}\rho
}\widetilde{\Phi } \mathcal{J}=-\widetilde{\Phi }(\mathcal{L}_\rho
\mathcal{J})=\widetilde{\Phi }\mathcal{N},
\]
thus $N$ is the image of $\mathcal{N}$ by $\widetilde{\Phi }$.\hfill
\hbox{\rlap{$\sqcap$}$\sqcup$}

This theorem shows that the splitting $V\mathcal{T}E^{*}\oplus
H\mathcal{T} E^{*}$ defined by $\mathcal{N}$ is mapped by $\Phi
_{*}$ into the splitting $V \mathcal{T}E\oplus H\mathcal{T}E$
defined by $N$, and it results:

\begin{corollary}
The following equalities hold
\[
\Phi _{*}\delta _\alpha ^{*}=\delta _\alpha ,\quad \Phi _{*}^{-1}\delta
_\alpha =\delta _\alpha ^{*},
\]
\[
\mathcal{N}_{\alpha \beta }(q,\mu )=-\left( N_\alpha ^\gamma
(x,y)+\mathcal{Q }_\alpha (\xi ^\gamma )\right) g_{\beta \gamma }.
\]
\[
\mathcal{R}_{\alpha \beta \gamma }g^{\alpha \varepsilon
}=\mathcal{R}_{\beta \gamma }^\varepsilon \circ \Phi ,\qquad
\mathcal{R}_{\beta \gamma }^\varepsilon \frac{\partial \zeta
_\alpha }{\partial y^\varepsilon }= \mathcal{R}_{\alpha \beta
\gamma }\circ \Phi ^{-1},
\]
\[
N_\alpha ^\beta \circ \Phi =-\delta _\alpha ^{*}(\xi ^\beta
),\quad \mathcal{ N}_{\alpha \beta }\circ \Phi ^{-1}=-\delta
_\alpha (\zeta _\beta ),
\]
\end{corollary}

\textbf{Proof.} We have on one hand
\[
\Phi _{*}^{-1}(\delta _\alpha )=\delta _\alpha
^{*}=\mathcal{Q}_\alpha + \mathcal{N}_{\alpha \gamma
}\mathcal{P}^\gamma ,
\]
and on the other hand
\[
\Phi _{*}^{-1}(\delta _\alpha )=\Phi _{*}^{-1}(\mathcal{X}_\alpha
-N_\alpha ^\beta \mathcal{V}_\beta )=\mathcal{Q}_\alpha -g_{\gamma
\varepsilon } \mathcal{Q}_\alpha (\xi ^\alpha )\mathcal{P}^\gamma
-N_\alpha ^\beta g_{\beta \gamma }\mathcal{P}^\gamma .
\]
Therefore, we get $\mathcal{N}_{\alpha \gamma }=-\left( N_\alpha
^\varepsilon +\mathcal{Q}_\alpha (\xi ^\varepsilon )\right) g_{\gamma
\varepsilon }.$ Next, since $\Phi ^{-1}$ is a diffeomorphism, we obtain
\[
\lbrack \Phi _{*}^{-1}(\delta _\alpha ),\Phi _{*}^{-1}(\delta
_\beta )]_{ \mathcal{T}E^{*}}=L_{\alpha \beta }^\gamma \Phi
_{*}^{-1}(\delta _\gamma )+ \mathcal{R}_{\alpha \beta }^\gamma
\Phi _{*}^{-1}(\mathcal{V}_\gamma ),
\]
\[
\lbrack \Phi _{*}^{-1}(\delta _\alpha ),\Phi _{*}^{-1}(\delta
_\beta )]_{ \mathcal{T}E^{*}}=[\delta _\alpha ^{*},\delta _\beta
^{*}]_{\mathcal{T} E^{*}}=L_{\alpha \beta }^\gamma \delta _\gamma
^{*}+\mathcal{R}_{\alpha \beta \gamma }\mathcal{P}^\gamma ,
\]
and using (1.5.38) it results the relation between the curvature tensors of
connections $N$ on $\mathcal{T}E$ and $\mathcal{N}$ on $\mathcal{T}E^{*}$.
By direct computations the other relations are obtained. \hfill
\hbox{\rlap{$\sqcap$}$\sqcup$}

Let us consider $\mathcal{S}$ a semispray on $\mathcal{T}E$ and $\Phi
^{-1}:E\rightarrow E^{*}$ the diffeomorphism given by (1.5.37). Then we set
\cite{Po25}:

\begin{theorem}
The section $\rho =\Phi _{*}^{-1}\mathcal{S}$ is a
semi-Hamiltonian section on $\mathcal{T}E^{*}$ if and only if the
nonlinear connection on $\mathcal{T} E$ induced by semispray is
the canonical nonlinear connection induced by the regular
Lagrangian.
\end{theorem}

\textbf{Proof}. Considering $\mathcal{S}=y^\alpha
\mathcal{X}_\alpha + \mathcal{S}^\alpha \mathcal{V}_\alpha $ we
have from (1.5.38)
\[
\Phi _{*}^{-1}\mathcal{S}=\xi ^\alpha \mathcal{Q}_\alpha +\left(
-\xi ^\alpha g_{\gamma \varepsilon }\mathcal{Q}_\alpha (\xi
^\varepsilon )+ \mathcal{S}^\alpha g_{\alpha \gamma }\right)
\mathcal{P}^\gamma .
\]
The conditions (1.5.28) $i),$ $ii)$ and (1.5.38) lead to
\[
g^{\theta \beta }\left( \mathcal{X}_\theta (\zeta _\alpha
)-\mathcal{X} _\alpha (\zeta _\theta )+\mathcal{V}_\theta
(\mathcal{S}^\varepsilon )g_{\varepsilon \alpha }+\xi ^\varepsilon
\mathcal{X}_\varepsilon (g_{\alpha \theta
})+\mathcal{S}^\varepsilon \mathcal{V}_\varepsilon (g_{\alpha
\theta })+\zeta _\varepsilon L_{\alpha \theta }^\varepsilon
\right) +y^\varepsilon L_{\alpha \varepsilon }^\beta =0,
\]
and using (1.3.25) it follows
\[
\mathcal{X}_\alpha (\zeta _\theta )-\mathcal{X}_\theta (\zeta
_\alpha )= \mathcal{S}(g_{\alpha \theta })-2N_\theta ^\varepsilon
g_{\varepsilon \alpha }+y^\beta L_{\theta \beta }^\varepsilon
g_{\varepsilon \alpha }+y^\beta L_{\alpha \beta }^\varepsilon
g_{\varepsilon \theta }+\zeta _\varepsilon L_{\alpha \theta
}^\varepsilon ,
\]
which concludes that (change $\alpha $ with $\theta $ and totalizing)
\[
\mathcal{S}(g_{\alpha \theta })-N_\theta ^\varepsilon g_{\varepsilon \alpha
}-N_\alpha ^\varepsilon g_{\varepsilon \theta }+y^\beta \left( L_{\theta
\beta }^\varepsilon g_{\varepsilon \alpha }+L_{\alpha \beta }^\varepsilon
g_{\varepsilon \theta }\right) =0\stackrel{(1.4.19)}{\Leftrightarrow }\nabla
g\left( \mathcal{V}_\alpha ,\mathcal{V}_\theta \right) =0,
\]
and
\[
\mathcal{X}_\alpha (\zeta _\theta )-\mathcal{X}_\theta (\zeta _\alpha
)=N_\theta ^\varepsilon g_{\varepsilon \alpha }-N_\alpha ^\varepsilon
g_{\varepsilon \theta }+\zeta _\varepsilon L_{\alpha \theta }^\varepsilon ,
\]
which is equivalent with
\[
N_\theta ^\varepsilon g_{\varepsilon \alpha }-N_\alpha
^\varepsilon g_{\varepsilon \theta }+\frac{\partial ^2L}{\partial
x^i\partial y^\alpha } \sigma _\beta ^i-\frac{\partial
^2L}{\partial x^i\partial y^\beta }\sigma _\alpha
^i-\frac{\partial L}{\partial y^\varepsilon }L_{\alpha \beta
}^\varepsilon =0\stackrel{(1.4.20)}{\Leftrightarrow }\omega
_L(\mathrm{h} \rho ,\mathrm{h}\nu )=0.
\]

From the Theorem 1.4.2 we get that the section $\rho =\Phi
_{*}^{-1}\mathcal{ S}$ is semi-Hamiltonian if and only if the
connection induced by semispray is the canonical nonlinear
connection (1.4.18). When $\rho =\Phi _{*}^{-1} \mathcal{S}$ is a
Hamiltonian section we have the same result, but moreover, the
condition (1.5.28) $iii)$ seems to lead to the third Helmholtz
condition on a Lie algebroid (1.4.22) and the work is in progress.
\hfill \hbox{\rlap{$\sqcap$}$\sqcup$}\newpage

\subsection{\textbf{Dynamical covariant derivative and metric non-linear
connection on} $\mathcal{T}E^{*}$}

\begin{definition}
A map $\nabla :\frak{T}(\mathcal{T}E^{*}\backslash
\{0\})\rightarrow \frak{T} (\mathcal{T}E^{*}\backslash \{0\})$ is
said to be a tensor derivation on $ \mathcal{T}E^{*}\backslash
\{0\}$ if the following conditions are satisfied:
\\i) $\nabla $ is $\Bbb{R}$-linear\\ii) $\nabla $ is type preserving, i.e. $
\nabla \frak{T}_s^r(\mathcal{T}E^{*}\backslash \{0\})\subset
\frak{T}_s^r( \mathcal{T}E^{*}\backslash \{0\})$, for each
$(r,s)\in \Bbb{N}\times \Bbb{N}$
\\iii) $\nabla $ obeys the Leibnitz rule $\nabla (P\otimes S)=\nabla
P\otimes S+P\otimes \nabla S$ for any tensors $P,S$ on
$\mathcal{T} E^{*}\backslash \{0\}$ \\iv) $\nabla \,$commutes with
any contractions.
\end{definition}

We consider the $\Bbb{R}$-linear map $\nabla _\rho :\Gamma
(\mathcal{T} E^{*}\backslash \{0\})\rightarrow \Gamma
(\mathcal{T}E^{*}\backslash \{0\})$ by
\begin{equation}
\nabla _\rho X=\mathrm{h}[\rho
,\mathrm{h}X]_{\mathcal{T}E^{*}}+\mathrm{v} [\rho
,\mathrm{v}X]_{\mathcal{T}E^{*}},\quad \forall X\in \Gamma
(\mathcal{T} E^{*}\backslash \{0\}).  \tag{1.6.1}
\end{equation}
where $\rho $ is a $\mathcal{J}-$regular section and it follows that
\begin{equation}
\nabla _\rho (fX)=\rho (f)X+f\nabla _\rho X,\quad \forall f\in
C^\infty ( \mathcal{T}E^{*}\backslash \{0\}),  \tag{1.6.2}
\end{equation}
Any tensor derivation on $\mathcal{T}E^{*}\backslash \{0\}$ is
completely determined by its actions on smooth functions and
sections on $\mathcal{T} E^{*}\backslash \{0\}$ (see \cite{Si}
generalized Willmore's theorem, p. 1217). Therefore, there exists
a unique tensor derivation $\nabla _\rho $ on
$\mathcal{T}E^{*}\backslash \{0\}$ such that
\[
\nabla _\rho \mid _{C^\infty (\mathcal{T}E^{*}\backslash \{0\})}=\rho ,\quad
\nabla _\rho \mid _{\Gamma (\mathcal{T}E^{*}\backslash \{0\})}=\nabla _0.
\]
We will call the tensor derivation $\nabla _\rho $, the \textit{dynamical
covariant derivative} induced by the $\mathcal{J}$-regular section $\rho $
and a nonlinear connection $\mathcal{N}$.

\begin{proposition}
The following equations hold
\begin{equation}
\lbrack \rho ,\mathcal{P}^\beta ]_{\mathcal{T}E^{*}}=-t^{\alpha
\beta }\delta _\alpha ^{*}+\left( t^{\alpha \beta
}\mathcal{N}_{\alpha \gamma }- \frac{\partial \rho _\gamma
}{\partial \mu _\beta }\right) \mathcal{P} ^\gamma ,  \tag{1.6.3}
\end{equation}
\begin{equation}
\lbrack \rho ,\delta _\beta ^{*}]_{\mathcal{T}E^{*}}=-\left( \delta _\beta
^{*}(\xi ^\alpha )+\xi ^\varepsilon L_{\varepsilon \beta }^\alpha \right)
\delta _\alpha ^{*}+\mathcal{R}_{\beta \gamma }\mathcal{P}^\gamma ,
\tag{1.6.4}
\end{equation}
\end{proposition}

where
\begin{equation}
\mathcal{R}_{\beta \gamma }=\delta _\beta ^{*}(\xi ^\alpha
)\mathcal{N} _{\alpha \gamma }+\rho (\mathcal{N}_{\beta \gamma
})-\delta _\beta ^{*}(\rho _\gamma )-\xi ^\varepsilon
L_{\varepsilon \beta }^\alpha \mathcal{N}_{\alpha \gamma },
\tag{1.6.5}
\end{equation}

The action of $\nabla _\rho $ on the Berwald basis has the form
\[
\nabla _\rho \mathcal{P}^\beta =\mathrm{v}[\rho ,\mathcal{P}^\beta
]_{ \mathcal{T}E^{*}}=\left( t^{\alpha \beta }\mathcal{N}_{\alpha
\gamma }-\frac{
\partial \rho _\gamma }{\partial \mu _\beta }\right) \mathcal{P}^\gamma ,
\]
\[
\nabla _\rho \delta _\beta ^{*}=\mathrm{h}[\rho ,\delta _\beta
^{*}]_{ \mathcal{T}E^{*}}=-\left( \delta _\beta ^{*}(\xi ^\alpha
)+\xi ^\varepsilon L_{\varepsilon \beta }^\alpha \right) \delta
_\alpha ^{*}.
\]
For a pseudo-Riemannian metric $g$ on $\mathcal{T}E^{*}$ the
action of $ \nabla _\rho $ is given by
\begin{equation}
\nabla _\rho g(X,Y)=\rho (g(X,Y))-g(\nabla _\rho X,Y)-g(X,\nabla _\rho Y),
\tag{1.6.6}
\end{equation}
which in local coordinates leads to
\[
g_{/}^{\alpha \beta }:=\nabla _\rho g\left( \mathcal{P}^\alpha
,\mathcal{P} ^\beta \right) =\rho (g^{\alpha \beta
})-g^{\varepsilon \beta }\left( t^{\alpha \gamma
}\mathcal{N}_{\gamma \varepsilon }-\frac{\partial \rho
_\varepsilon }{\partial \mu _\alpha }\right) -g^{\varepsilon
\alpha }\left( t^{\beta \gamma }\mathcal{N}_{\gamma \varepsilon
}-\frac{\partial \rho _\varepsilon }{\partial \mu _\beta }\right)
.
\]

\begin{definition}
A nonlinear connection is called metric or compatible with the metric tensor
$g$ if $\nabla _\rho g=0,$ for all $\mathcal{J}$-regular sections $\rho $,
that is
\begin{equation}
\rho (g(X,Y))=g(\nabla _\rho X,Y)+g(X,\nabla _\rho Y),\quad \forall X,Y\in
\Gamma (V\mathcal{T}E^{*}).  \tag{1.6.7}
\end{equation}
\end{definition}

\begin{theorem}
The connection $\widetilde{\mathcal{N}}$ with the coefficients
\begin{equation}
\widetilde{\mathcal{N}_{\alpha \beta }}=\mathcal{N}_{\alpha \beta }+\frac
12t_{\alpha \varepsilon }g_{\beta \gamma }g_{/}^{\varepsilon \gamma },
\tag{1.6.8}
\end{equation}
is a metric nonlinear connection.
\end{theorem}

\textbf{Proof}. Let us consider the dynamical covariant derivative induced
by $\rho $ and $\widetilde{\mathcal{N}}$ given by
\[
\nabla _\rho g\left( \mathcal{P}^\alpha ,\mathcal{P}^\beta \right)
=\rho (g^{\alpha \beta })-g^{\varepsilon \beta }\left( t^{\alpha
\gamma } \widetilde{\mathcal{N}}_{\gamma \varepsilon
}-\frac{\partial \rho _\varepsilon }{\partial \mu _\alpha }\right)
-g^{\varepsilon \alpha }\left( t^{\beta \gamma
}\widetilde{\mathcal{N}}_{\gamma \varepsilon }-\frac{
\partial \rho _\varepsilon }{\partial \mu _\beta }\right) ,
\]
and using (1.6.8) it follows
\[
\nabla _\rho g\left( \mathcal{P}^\alpha ,\mathcal{P}^\beta \right)
=g_{/}^{\alpha \beta }-\frac 12g^{\varepsilon \beta }t^{\alpha \gamma
}t_{\gamma r}g_{\varepsilon p}g_{/}^{rp}-\frac 12g^{\alpha \varepsilon
}t^{\beta \gamma }t_{\gamma r}g_{\varepsilon p}g_{/}^{rp}=0,
\]
that is $\widetilde{\mathcal{N}}$ is a metric nonlinear connection.\hfill
\hbox{\rlap{$\sqcap$}$\sqcup$}
\newpage

\subsubsection{\textbf{Nonlinear connection induced by a }$\mathcal{J}$
\textbf{-regular section}}

If $\mathcal{J}$ is an adapted tangent structure and $\rho $ is a
$\mathcal{J }$- regular section then \cite{Hr4}
\[
\mathcal{N}=-\mathcal{L}_\rho \mathcal{J},
\]
is a nonlinear connection on $\mathcal{T}E^{*}$ with local coefficients
given by (1.5.26)
\[
\mathcal{N}_{\alpha \beta }=\frac 12\left( t_{\alpha \gamma
}\frac{\partial \rho _\beta }{\partial \mu _\gamma }-\sigma
_\alpha ^it_{\gamma \beta }\frac{
\partial \xi ^\gamma }{\partial q^i}-\rho (t_{\alpha \beta })+\xi ^\gamma
t_{\varepsilon \beta }L_{\gamma \alpha }^\varepsilon \right) .
\]

\begin{definition}
The Jacobi endomorphism $\psi $ is given by
\[
\psi =\mathrm{v}[\rho ,\mathrm{h}X]_{\mathcal{T}E^{*}}.
\]
\end{definition}

Locally, from (1.6.4) we obtain that $\psi =\mathcal{R}_{\alpha
\beta } \mathcal{Q}^\alpha \otimes \mathcal{P}^\beta ,$ where
$\mathcal{R}_{\alpha \beta }$ is given by (1.6.8).
$\mathcal{R}_{\alpha \beta }$ are the local coefficients of the
Jacobi endomorphism.

\begin{proposition}
The following result holds
\[
\psi =i_\rho \Omega +\mathrm{v}[\mathrm{v}\rho
,\mathrm{h}X]_{\mathcal{T} E^{*}}.
\]
\end{proposition}

\textbf{Proof}. Indeed, $\psi (X)=\mathrm{v}[\rho
,\mathrm{h}X]_{\mathcal{T} E^{*}}=\mathrm{v}[\mathrm{h}\rho
,\mathrm{h}X]_{\mathcal{T}E^{*}}+\mathrm{v}[ \mathrm{v}\rho
,\mathrm{h}X]_{\mathcal{T}E^{*}}$ and $\Omega (\rho ,X)=
\mathrm{v}[\mathrm{h}\rho ,\mathrm{h}X]_{\mathcal{T}E^{*}},$ that
is $\psi (X)=\Omega (\rho ,X)+\mathrm{v}[\mathrm{v}\rho
,\mathrm{h}X]_{\mathcal{T} E^{*}}.$\hfill
\hbox{\rlap{$\sqcap$}$\sqcup$}

\begin{remark}
If $\rho $ is a horizontal section $\rho =\mathrm{h}\rho $, then we obtain
\[
\psi =i_\rho \Omega ,
\]
and locally it follows
\[
\rho _\gamma =\xi ^\alpha N_{\alpha \gamma },
\]
which yields
\begin{equation}
\mathcal{R}_{\alpha \beta }=\mathcal{R}_{\varepsilon \alpha \beta }\xi
^\varepsilon .  \tag{1.6.9}
\end{equation}
and the local relation between the Jacobi endomorphism and the curvature
tensor of the nonlinear connection is obtained.\newpage\
\end{remark}

\subsubsection{\textbf{Hamiltonian case}}

In what follows, we consider a regular Hamiltonian $\mathcal{H}$ $:$ $E^{*}$
$\rightarrow \Bbb{R}$, that is the matrix
\[
g^{\alpha \beta }(q,\mu )=\frac{\partial ^2\mathcal{H}}{\partial \mu _\alpha
\partial \mu _\beta },
\]
is nondegenerate. This regular Hamiltonian determines a symmetric canonical
nonlinear connection given by (1.5.35). We have the following result:

\begin{theorem}
The canonical nonlinear connection induced by a regular
Hamiltonian $ \mathcal{H}$ is a metric nonlinear connection.
\end{theorem}

\textbf{Proof}. Introducing the coefficients (1.5.35) into the expression of
the dynamical covariant derivative and using (1.5.33) we obtain
\[
\nabla _\rho g\left( \mathcal{P}^\alpha ,\mathcal{P}^\beta \right)
=\sigma _\varepsilon ^i\frac{\partial g^{\alpha \beta }}{\partial
q^i}\frac{\partial \mathcal{H}}{\partial \mu _\varepsilon }-\sigma
_\varepsilon ^i\frac{
\partial \mathcal{H}}{\partial q^i}\frac{\partial g^{\alpha \beta }}{
\partial \mu _\varepsilon }-\mu _\gamma L_{\varepsilon \delta }^\gamma \frac{
\partial \mathcal{H}}{\partial \mu _\delta }\frac{\partial g^{\alpha \beta }
}{\partial \mu _\varepsilon }
\]
\[
-\frac 12\left( g^{\varepsilon \beta }g^{\alpha \gamma
}+g^{\varepsilon \alpha }g^{\beta \gamma }\right) \left[ \sigma
_l^i\left( \frac{\partial g_{\gamma \varepsilon }}{\partial \mu
_l}\frac{\partial \mathcal{H}}{
\partial q^i}-\frac{\partial g_{\gamma \varepsilon }}{\partial q^i}\frac{
\partial \mathcal{H}}{\partial \mu _l}\right) -\sigma _\varepsilon ^i\frac{
\partial ^2\mathcal{H}}{\partial q^i\partial \mu _l}g_{\gamma l}\right.
\]
\[
\ \left. -\sigma _\gamma ^i\frac{\partial ^2\mathcal{H}}{\partial
q^i\partial \mu _l}g_{\varepsilon l}+\mu _lL_{rs}^l\frac{\partial
\mathcal{H} }{\partial \mu _r}\frac{\partial g_{\gamma \varepsilon
}}{\partial \mu _s} +\mu _lL_{\gamma \varepsilon
}^l+\frac{\partial \mathcal{H}}{\partial \mu _\delta }\left(
g_{\gamma l}L_{\delta \varepsilon }^l+g_{\varepsilon l}L_{\delta
\gamma }^l\right) \right]
\]
\[
-g^{\varepsilon \beta }\sigma _\varepsilon ^i\frac{\partial
^2\mathcal{H}}{
\partial \mu _\alpha \partial q^i}-g^{\varepsilon \beta }L_{\varepsilon
\delta }^\alpha \frac{\partial \mathcal{H}}{\partial \mu _\delta }
-g^{\varepsilon \beta }g^{\alpha \delta }\mu _\gamma
L_{\varepsilon \delta }^\gamma
\]
\[
\ -g^{\varepsilon \alpha }\sigma _\varepsilon ^i\frac{\partial
^2\mathcal{H} }{\partial \mu _\beta \partial q^i}-g^{\varepsilon
\alpha }L_{\varepsilon \delta }^\beta \frac{\partial
\mathcal{H}}{\partial \mu _\delta } -g^{\varepsilon \alpha
}g^{\beta \delta }\mu _\gamma L_{\varepsilon \delta }^\gamma .
\]
From the equalities
\begin{equation}
g^{\varepsilon \beta }g^{\alpha \gamma }\frac{\partial g_{\gamma
\varepsilon }}{\partial \mu _l}=-g^{\varepsilon \beta }g_{\gamma
\varepsilon }\frac{
\partial g^{\alpha \gamma }}{\partial \mu _l}=-\frac{\partial g^{\alpha
\beta }}{\partial \mu _l},  \tag{1.6.10}
\end{equation}
\begin{equation}
g^{\varepsilon \beta }g^{\alpha \gamma }\frac{\partial g_{\gamma
\varepsilon }}{\partial q^i}=-g^{\varepsilon \beta }g_{\gamma
\varepsilon }\frac{
\partial g^{\alpha \gamma }}{\partial q^i}=-\frac{\partial g^{\alpha \beta }
}{\partial q^i},\quad L_{\varepsilon \delta }^\gamma =-L_{\delta \varepsilon
}^\gamma  \tag{1.6.11}
\end{equation}
by direct computation, it follows that $\nabla _\rho g\left(
\mathcal{P} ^\alpha ,\mathcal{P}^\beta \right) =0$, which ends the
proof.\hfill \hbox{\rlap{$\sqcap$}$\sqcup$}

\begin{theorem}
The canonical nonlinear connection induced by a regular Hamiltonian is the
unique metric and symmetric nonlinear connection.
\end{theorem}

\textbf{Proof}. Let us consider a metric and symmetric nonlinear
connection $ \mathcal{N}$ with the coefficients
$\mathcal{N}_{\gamma \varepsilon }$. Then we have
\[
\rho _{\mathcal{H}}(g^{\alpha \beta })=g^{\varepsilon \beta }\left(
g^{\alpha \gamma }\mathcal{N}_{\gamma \varepsilon }-\frac{\partial \rho
_\varepsilon }{\partial \mu _\alpha }\right) +g^{\varepsilon \alpha }\left(
g^{\beta \gamma }\mathcal{N}_{\gamma \varepsilon }-\frac{\partial \rho
_\varepsilon }{\partial \mu _\beta }\right)
\]
and using (1.5.33) we obtain
\[
\frac{\partial \mathcal{H}}{\partial \mu _l}\sigma
_l^i\frac{\partial g^{\alpha \beta }}{\partial q^i}-\left( \sigma
_l^i\frac{\partial \mathcal{H} }{\partial q^i}+\mu _\gamma
L_{l\varepsilon }^\gamma \frac{\partial \mathcal{ H}}{\partial \mu
_\varepsilon }\right) \frac{\partial g^{\alpha \beta }}{
\partial \mu _l}=
\]
\[
g^{\varepsilon \beta }g^{\alpha \gamma }\mathcal{N}_{\gamma
\varepsilon }+g^{\varepsilon \alpha }g^{\beta \gamma
}\mathcal{N}_{\gamma \varepsilon }+g^{\gamma \beta }\sigma _\gamma
^i\frac{\partial ^2\mathcal{H}}{\partial \mu _\alpha \partial
q^i}+g^{\gamma \beta }L_{\gamma \tau }^\alpha \frac{
\partial \mathcal{H}}{\partial \mu _\tau }+
\]
\[
g^{\gamma \beta }\mu _lL_{\gamma \tau }^l\frac{\partial
^2\mathcal{H}}{
\partial \mu _\alpha \partial \mu _\tau }+g^{\varepsilon \alpha }\sigma
_\varepsilon ^i\frac{\partial ^2\mathcal{H}}{\partial \mu _\beta
\partial q^i
}+g^{\varepsilon \alpha }L_{\varepsilon \tau }^\beta
\frac{\partial \mathcal{ H}}{\partial \mu _\tau }+g^{\varepsilon
\alpha }\mu _lL_{\varepsilon \tau }^l \frac{\partial
^2\mathcal{H}}{\partial \mu _\beta \partial \mu _\tau }.
\]
But the connection is symmetric (1.5.10) and using (1.6.10),
(1.6.11) we get $\pi $
\[
\sigma _l^i\frac{\partial \mathcal{H}}{\partial q^i}\frac{\partial
g_{\varepsilon \gamma }}{\partial \mu _l}-\sigma
_l^i\frac{\partial \mathcal{ H}}{\partial \mu _l}\frac{\partial
g_{\varepsilon \gamma }}{\partial q^i} +\mu _\tau L_{l\alpha
}^\tau \frac{\partial \mathcal{H}}{\partial \mu _\alpha
}\frac{\partial g_{\varepsilon \gamma }}{\partial \mu _l}=
\]
\[
2\mathcal{N}_{\varepsilon \gamma }+\mu _\tau L_{\gamma \varepsilon
}^\tau -g_{\alpha \varepsilon }\sigma _\gamma ^i\frac{\partial
^2\mathcal{H}}{
\partial \mu _\alpha \partial q^i}-g_{\alpha \gamma }\sigma _\varepsilon ^i
\frac{\partial ^2\mathcal{H}}{\partial \mu _\alpha \partial
q^i}-g_{\alpha \varepsilon }L_{\gamma l}^\alpha \frac{\partial
\mathcal{H}}{\partial \mu _l} -g_{\alpha \gamma }L_{\varepsilon
l}^\alpha \frac{\partial \mathcal{H}}{
\partial \mu _l},
\]
and we get the coefficients (1.5.35), which ends the proof.\hfill
\hbox{\rlap{$\sqcap$}$\sqcup$}

For the particular case of the cotangent bundle see \cite{Po16,
Po21}\newpage \

\subsection{\textbf{Poisson-Lie algebroids}}

Poisson manifolds were introduced by A. Lichnerowicz in his famous paper
\cite{Lic} and their properties were later investigated by A. Weinstein \cite
{We1}. The Poisson manifolds are the smooth manifolds equipped with a
Poisson bracket on their ring of functions. We remark that the cotangent
bundle of a Poisson manifold has a natural structure of Lie algebroid. In
this chapter we study the Poisson structures on Lie algebroids \cite{Po13,
Po15, Po23}.

\subsubsection{\textbf{Poisson structures on Lie algebroids}}

The Schouten-Nijenhuis bracket is given by \cite{Va1}

\begin{eqnarray*}
\lbrack X_1\wedge ...\wedge X_p,Y_1\wedge ...\wedge Y_q]
&=&(-1)^{p+1}\sum_{i=1}^p\sum_{j=1}^q(-1)^{i+j}[X_i,Y_j]\wedge X_1\wedge
...\wedge \\
&&\quad \stackrel{\symbol{94}}{X}_i\wedge ...\wedge X_p\wedge Y_1\wedge
...\wedge \stackrel{\symbol{94}}{Y}_j\wedge ...\wedge Y_q
\end{eqnarray*}

Let us consider the bivector (i.e. contravariant, skew-symmetric, 2-section)
$\Pi \in \Gamma (\wedge ^2E)$ given by the expression
\begin{equation}
\Pi =\frac 12\pi ^{\alpha \beta }(x)s_\alpha \wedge s_\beta .  \tag{1.7.1}
\end{equation}

\begin{definition}
The bivector $\Pi $ is a Poisson bivector on $E$ if and only if $[\Pi ,\Pi
]=0$, where [$\cdot $,$\cdot $] is Schouten-Nijenhuis bracket.
\end{definition}

\begin{proposition}
Locally the condition $[\Pi ,\Pi ]=0$ is expressed as
\begin{equation}
\sum_{(\alpha ,\varepsilon ,\delta )}(\pi ^{\alpha \beta }\sigma
_\beta ^i \frac{\partial \pi ^{\varepsilon \delta }}{\partial
x^i}+\pi ^{\alpha \beta }\pi ^{\gamma \delta }L_{\beta \gamma
}^\varepsilon )=0  \tag{1.7.2}
\end{equation}
\end{proposition}

If $\Pi $ is a Poisson bivector then the pair $(E,\Pi )$ is called the
\textit{Poisson-Lie} algebroid. The Poisson bracket on $M$ is given by
\[
\{f_1,f_2\}=\Pi (d^Ef_1,d^Ef_2),\quad f_1,f_2\in C^\infty (M)
\]
We also have the bundle map $\pi ^{\#}:E^{*}\rightarrow E$ defined by
\[
\pi ^{\#}\rho =i_\rho \Pi ,\quad \rho \in \Gamma (E^{*}).
\]
Let us consider the bracket
\[
\lbrack \rho ,\theta ]_\pi =\mathcal{L}_{\pi ^{\#}\rho }\theta
-\mathcal{L} _{\pi ^{\#}\theta }\rho -d^E(\Pi (\rho ,\theta )),
\]
where $\mathcal{L}$ is Lie derivative and $\rho ,\theta \in \Gamma (E^{*}).$
With respect to this bracket and the usual Lie bracket on vector fields, the
map $\widetilde{\sigma }:E^{*}\rightarrow TM$ given by
\[
\widetilde{\sigma }=\sigma \circ \pi ^{\#},
\]
is a Lie algebra homomorphism
\[
\widetilde{\sigma }[\rho ,\theta ]_\pi =[\widetilde{\sigma }\rho
,\widetilde{ \sigma }\theta ].
\]
The bracket $[.,.]_\pi $ satisfies also the Leibniz rule
\[
\lbrack \rho ,f\theta ]_\pi =f[\rho ,\theta ]_\pi +\widetilde{\sigma }(\rho
)(f)\theta ,
\]
and it results that $(E^{*},[.,.]_\pi ,\widetilde{\sigma })$ is a Lie
algebroid \cite{Kos1}.

Next, we can define the contravariant exterior differential $d^\pi
:$ $ \bigwedge^k(E^{*})\rightarrow $ $\bigwedge^{k+1}(E^{*})$ by
\begin{eqnarray*}
d^\pi \omega (s_1,...,s_{k+1})
&=&\stackrel{k+1}{\sum_{i=1}}(-1)^{i+1} \widetilde{\sigma
}(s_i)\omega (s_1,...,\stackrel{\symbol{94}}{s}
_i,...,s_{k+1})+ \\
&&\ \ \ \ \ \ \ +\sum_{1\leq i<j\leq k+1}(-1)^{i+j}\omega
([s_{i,}s_j]_\pi
,s_1,...,\stackrel{\symbol{94}}{s_i},...,\stackrel{\symbol{94}}{s_j}
,...s_{k+1}).
\end{eqnarray*}

Accordingly, we get the cohomology of Lie algebroid $E^{*}$ with
the anchor $ \widetilde{\sigma }$ and the bracket $[.,.]_\pi $
which generalize the Poisson cohomology of Lichnerowicz for
Poisson manifolds \cite{Va1}.

In the following we deal with the notion of contravariant connection on Lie
algebroids, which generalize the similar notion on Poisson manifolds \cite
{Fe1}, \cite{Lic}.

\begin{definition}
If $\rho ,\theta $ $\in \Gamma (E^{*})$ and $\Phi ,\Psi \in \Gamma (E)$ then
the linear contravariant connection is an application $D:\Gamma
(E^{*})\times \Gamma (E)\rightarrow \Gamma (E)$ which satisfies the relations

i) $\ \ D_{\rho +\theta }\Phi =D_\rho \Phi +D_\theta \Phi ,$

ii) \ $D_\rho (\Phi +\Psi )=D_\rho \Phi +D_\rho \Psi ,$

iii) $D_{f\rho }\Phi =fD_\rho \Phi ,\quad $

iv) $D_\rho (f\Phi )=fD_\rho \Phi +\widetilde{\sigma }(\rho )(f)\Phi ,\quad
f\in C^\infty (M).$
\end{definition}

The contravariant connection induces a contravariant derivative $D_\alpha
:\Gamma (E)\rightarrow \Gamma (E)$ such that the following equalities are
fulfilled
\[
D_{f_1\alpha _1+f_2\alpha _2}=f_1D_{\alpha _1}+f_2D_{\alpha _2},\quad f_i\in
C^\infty (M),\quad \alpha _i\in \Gamma (E^{*})
\]
\[
D_\rho (f\theta )=fD_\rho \theta +\widetilde{\sigma }(\rho )(f)\theta ,\quad
f\in C^\infty (M),\quad \rho ,\theta \in \Gamma (E^{*})
\]

In the case where the contravariant connection $D$ is induced by a covariant
connection $\nabla $ on a Lie algebroid $E$ (see \cite{Co}) we have $D_\rho
=\nabla _{\pi ^{\#}\rho }$.

\begin{definition}
The torsion and curvature of linear contravariant connection are given by
\[
T(\rho ,\theta )=D_\rho \theta -D_\theta \rho -[\rho ,\theta ]_\pi ,
\]
\[
R(\rho ,\theta )\mu =D_\rho D_\theta \mu -D_\theta D_\rho \mu -D_{[\rho
,\theta ]_\pi }\mu ,
\]
\end{definition}

where $\rho ,\theta ,\mu \in \Gamma (E^{*}).$

The curvature tensor satisfies the equalities
\[
R(\rho ,\theta )=-R(\theta ,\rho ),\ R(f\rho ,\theta )=fR(\rho ,\theta ).
\]
The Bianchi identities have the following form
\[
\sum_{\rho,\theta,\mu}\left( D_\rho R(\theta ,\mu )+R(T(\rho
,\theta ),\mu )\right) =0,
\]
\[
\sum_{\rho,\theta,\mu}\left( R(\rho ,\theta )\mu -T(T(\rho ,\theta
),\mu )-D_\rho T(\theta ,\mu )\right) =0.
\]
In local coordinates we define the Christoffel symbols $\Gamma _\gamma
^{\alpha \beta }$ by the formula
\[
D_{s^\alpha }s^\beta =\Gamma _\gamma ^{\alpha \beta }s^\gamma ,
\]
and under a change of coordinates
\begin{equation}
\left\{
\begin{array}{l}
x^{i^{\prime }}=x^{i^{\prime }}(x^i),\quad i,i^{\prime }=\overline{1,n}\quad
on\quad M \\
y^{\alpha ^{\prime }}=A_\alpha ^{\alpha ^{\prime }}y^\alpha ,\quad \alpha
,\alpha ^{\prime }=\overline{1,m}\quad on\quad E,
\end{array}
\right.  \tag{1.7.3}
\end{equation}
corresponding to a new base $s^{\alpha ^{\prime }}=A_\alpha ^{\alpha
^{\prime }}s^\alpha $, these symbols transform according to
\begin{equation}
\Gamma _{\gamma ^{\prime }}^{\alpha ^{\prime }\beta ^{\prime }}=A_\alpha
^{\alpha ^{\prime }}A_\beta ^{\beta ^{\prime }}A_{\gamma ^{\prime }}^\gamma
\Gamma _\gamma ^{\alpha \beta }+A_\alpha ^{\alpha ^{\prime }}A_{\gamma
^{\prime }}^\gamma \sigma _\varepsilon ^i\frac{\partial A_\gamma ^{\beta
^{\prime }}}{\partial x^i}\pi ^{\alpha \varepsilon }.  \tag{1.7.4}
\end{equation}
If we denote $T(s^\alpha ,s^\beta )=T_\gamma ^{\alpha \beta
}s^\gamma $ and $ R(s^\alpha ,s^\beta )s^\gamma =R_\delta ^{\alpha
\beta \gamma }s^\delta $ then, under a change of coordinates, we
obtain that

\[
T_{\gamma ^{\prime }}^{\alpha ^{\prime }\beta \prime }=A_\alpha ^{\alpha
^{\prime }}A_\beta ^{\beta ^{\prime }}A_{\gamma ^{\prime }}^\gamma T_\gamma
^{\alpha \beta },\qquad R_{\delta ^{\prime }}^{\alpha ^{\prime }\beta
^{\prime }\gamma ^{\prime }}=A_\alpha ^{\alpha ^{\prime }}A_\beta ^{\beta
^{\prime }}A_\gamma ^{\gamma ^{\prime }}A_{\delta ^{\prime }}^\delta
R_\delta ^{\alpha \beta \gamma }.
\]

\begin{proposition}
The local components of torsion and curvature of linear contravariant
connection are
\[
T_\varepsilon ^{\alpha \beta }=\Gamma _\varepsilon ^{\alpha \beta }-\Gamma
_\varepsilon ^{\beta \alpha }-\pi ^{\alpha \gamma }L_{\gamma \varepsilon
}^\beta +\pi ^{\beta \gamma }L_{\gamma \varepsilon }^\alpha -\sigma
_\varepsilon ^i\frac{\partial \pi ^{\alpha \beta }}{\partial x^i},
\]
\[
R_\delta ^{\alpha \beta \gamma }=\Gamma _\delta ^{\alpha
\varepsilon }\Gamma _\varepsilon ^{\beta \gamma }-\Gamma _\delta
^{\beta \varepsilon }\Gamma _\varepsilon ^{\alpha \gamma }+\pi
^{\alpha \varepsilon }\sigma _\varepsilon ^i\frac{\partial \Gamma
_\delta ^{\beta \gamma }}{\partial x^i}-\pi ^{\beta \varepsilon
}\sigma _\varepsilon ^i\frac{\partial \Gamma _\delta ^{\alpha
\gamma }}{\partial x^i}+(\pi ^{\beta \nu }L_{\nu \varepsilon
}^\alpha -\pi ^{\alpha \nu }L_{\nu \varepsilon }^\beta -\sigma
_\varepsilon ^i\frac{
\partial \pi ^{\alpha \beta }}{\partial x^i})\Gamma _\delta ^{\varepsilon
\gamma }.
\]
\end{proposition}

\begin{definition}
A tensor field $T$ on $E$ is called \textit{parallel }with respect
to\textit{ \ }$D$ if and only if $DT=0.$
\end{definition}

\begin{definition}
A contravariant connection $D$ is called a Poisson connection if the Poisson
bivector $\Pi $ is parallel with respect to $D$.
\end{definition}

\begin{remark}
If the Poisson connection $D$ is induced by a covariant connection $\nabla $
(i.e. $D\Pi =0$, $D_\rho =\nabla _{\pi ^{\#}\rho }$, $\pi ^{\#}D_\rho \phi
=\nabla _{\pi ^{\#}\rho }\pi ^{\#}\phi $) then the torsion and curvature
tensors of the both connections are related by the following equalities
\begin{eqnarray*}
T^\nabla (\pi ^{\#}\rho ,\pi ^{\#}\theta )=\pi ^{\#}T^D(\rho ,\theta ),\qquad
\\
R^\nabla (\pi ^{\#}\rho ,\pi ^{\#}\theta )\pi ^{\#}\mu =\pi ^{\#}R^D(\rho
,\theta )\mu ,\text{ }\forall \rho ,\theta ,\mu \in \Gamma (E^{*}).
\end{eqnarray*}
\end{remark}

Let $T$ be a tensor of type $(r,s)$ with the components $
T_{j_1...j_s}^{i_1...i_r}$ and $\theta =\theta _\alpha s^\alpha $
a section of $E^{*}$. The local coordinates expression of
contravariant derivative is given by
\[
D_\theta T=\theta _\alpha T_{j_1...j_s}^{i_1...i_r}/^\alpha s_{i_1}\otimes
\cdot \cdot \cdot \otimes s_{i_r}\otimes s^{j_1}\otimes \cdot \cdot \cdot
\otimes s^{j_s},
\]
where
\[
T_{j_1...j_s}^{i_1...i_r}/^\alpha =\pi ^{\alpha \varepsilon
}\sigma _\varepsilon ^i\frac{\partial
T_{j_1...j_s}^{i_1...i_r}}{\partial x^i} +\sum_{a=1}^r\left(
\Gamma _\varepsilon ^{i_a\alpha }T_{j_1...j_s}^{i_1...\varepsilon
...i_r}\right) -\sum_{b=1}^s\left( \Gamma _{j_b}^{\varepsilon
\alpha }T_{j_1...\varepsilon ...j_s}^{i_1...i_r}\right) ,
\]
and $/$ denote the \textit{contravariat derivative operator}.

Let us consider a contravariant connection $D$ with the coefficients $\Gamma
_\gamma ^{\alpha \beta }$. We have

\begin{proposition}
The contravariant connection $\overline{D}$ with the coefficients given by
\begin{equation}
\overline{\Gamma }_\gamma ^{\alpha \beta }=\Gamma _\gamma ^{\alpha \beta
}-\frac 12\pi _{\gamma \varepsilon }\pi ^{\alpha \varepsilon }/^\beta ,
\tag{1.7.5}
\end{equation}
is a Poisson connection.
\end{proposition}

\textbf{Proof}. Considering $\overline{/}$ the contravariant derivative
operator with respect to contravariant connection $\overline{D}$, we get
\[
\pi ^{\beta \gamma }\overline{/}^\alpha =\pi ^{\alpha \varepsilon
}\sigma _\varepsilon ^i\frac{\partial \pi ^{\beta \gamma
}}{\partial x^i}+\overline{ \Gamma }_\varepsilon ^{\beta \alpha
}\pi ^{\varepsilon \gamma }+\overline{ \Gamma }_\varepsilon
^{\gamma \alpha }\pi ^{\beta \varepsilon }=
\]
\[
=\pi ^{\alpha \varepsilon }\sigma _\varepsilon ^i\frac{\partial \pi ^{\beta
\gamma }}{\partial x^i}+\left( \Gamma _\varepsilon ^{\beta \alpha }-\frac
12\pi _{\varepsilon \tau }\pi ^{\beta \tau }/^\alpha \right) \pi
^{\varepsilon \gamma }+\left( \Gamma _\varepsilon ^{\gamma \alpha }-\frac
12\pi _{\varepsilon \tau }\pi ^{\gamma \tau }/^\alpha \right) \pi ^{\beta
\varepsilon }
\]
\[
=\pi ^{\beta \gamma }/^\alpha -\frac 12\pi ^{\beta \gamma }/^\alpha -\frac
12\pi ^{\beta \gamma }/^\alpha =0.
\]
\hfill\hbox{\rlap{$\sqcap$}$\sqcup$}

\begin{proposition}
a) The functions
\begin{equation}
\Gamma _\gamma ^{\alpha \beta }=\sigma _\gamma ^i\frac{\partial \pi ^{\alpha
\beta }}{\partial x^i},  \tag{1.7.6}
\end{equation}
are the coefficients of a contravariant connection.

b) The contravariant connection with the coefficients given by (1.7.6) is a
Poisson connection if and only if
\begin{equation}
\sum_{\alpha,\varepsilon,\delta}\pi ^{\alpha \beta }\pi ^{\gamma
\delta }L_{\beta \gamma }^\varepsilon =0. \tag{1.7.7}
\end{equation}
\end{proposition}

\textbf{Proof}. a) Using the change of coordinates (1.7.3) and the fact that
the structure function $\sigma _\alpha ^i$ change by the rule \cite{Ma1}
\[
\sigma _{\alpha ^{\prime }}^{i^{\prime }}A_\alpha ^{\alpha
^{\prime }}=\frac{
\partial x^{i^{\prime }}}{\partial x^i}\sigma _\alpha ^i,
\]
we obtain that the coefficients (1.7.6) satisfy the transformation law
(1.7.4).

b) We obtain that
\[
\pi ^{\beta \gamma }/^\alpha =\pi ^{\alpha \varepsilon }\sigma _\varepsilon
^i\frac{\partial \pi ^{^{\beta \gamma }}}{\partial x^i}-\Gamma _\varepsilon
^{\beta \alpha }\pi ^{\varepsilon \gamma }-\Gamma _\varepsilon ^{\gamma
\alpha }\pi ^{\beta \varepsilon },
\]
and using the equality (1.7.6), it results
\[
\pi ^{\beta \gamma }/^\alpha =\sum_{(\alpha ,\beta ,\gamma )}\pi
^{\alpha \varepsilon }\sigma _\varepsilon ^i\frac{\partial \pi
^{\beta \gamma }}{
\partial x^i}.
\]
From the condition $[\Pi ,\Pi ]=0$, locally given by the equation (1.7.2),
it results that $\pi ^{\beta \gamma }/^\alpha =0$ if and only if the
required relation is fulfilled.\hfill\hbox{\rlap{$\sqcap$}$\sqcup$}

\begin{remark}
Considering
\[
\Gamma _\gamma ^{\alpha \beta }=\sigma _\gamma ^i\frac{\partial \pi ^{\alpha
\beta }}{\partial x^i}
\]
in relation (1.7.5) we obtain a Poisson connection $\overline{D}$ with the
coefficients
\[
\overline{\Gamma }_\gamma ^{\alpha \beta }=\sigma _\gamma ^i\frac{\partial
\pi ^{\alpha \beta }}{\partial x^i}-\frac 12\pi _{\gamma \varepsilon }\pi
^{\alpha \varepsilon }/^\beta
\]
which depends only on the Poisson bivector and structural functions of Lie
algebroid.
\end{remark}

\begin{proposition}
The set of Poisson connections on Lie algebroid are given by
\[
\overline{\Gamma }_\gamma ^{\alpha \beta }=\Gamma _\gamma ^{\alpha \beta
}+\Omega _{\gamma \nu }^{\alpha \varepsilon }X_\varepsilon ^{\nu \beta },
\]
where
\[
\Omega _{\gamma \nu }^{\alpha \varepsilon }=\frac 12\left( \delta _\nu
^\alpha \delta _\gamma ^\varepsilon -\pi _{\gamma \nu }\pi ^{\alpha
\varepsilon }\right) ,
\]
and $D(\Gamma _\gamma ^{\alpha \beta })$ is a Poisson connection
with $ X_\varepsilon ^{\delta \beta }$ an arbitrary tensor.
\end{proposition}

\textbf{Proof}. By straightforward computation it results
\[
\pi ^{\beta \gamma }\overline{/}^\alpha =\pi ^{\alpha \varepsilon
}\sigma _\varepsilon ^i\frac{\partial \pi ^{\beta \gamma
}}{\partial x^i}+\overline{ \Gamma }_\varepsilon ^{\beta \alpha
}\pi ^{\varepsilon \gamma }+\overline{ \Gamma }_\varepsilon
^{\gamma \alpha }\pi ^{\beta \varepsilon }=
\]
\[
\pi ^{\beta \gamma }/^\alpha +\frac 12\pi ^{\varepsilon \gamma }(\delta _\nu
^\beta \delta _\varepsilon ^\theta -\pi _{\varepsilon \nu }\pi ^{\beta
\theta })X_\theta ^{\nu \alpha }+\frac 12\pi ^{\beta \varepsilon }(\delta
_\nu ^\gamma \delta _\varepsilon ^\theta -\pi _{\varepsilon \nu }\pi
^{\gamma \theta })X_\theta ^{\nu \alpha }=
\]
\[
\pi ^{\beta \gamma }/^\alpha +\frac 12\pi ^{\theta \gamma }X_\theta ^{\beta
\alpha }-\frac 12\pi ^{\beta \theta }X_\theta ^{\gamma \alpha }+\frac 12\pi
^{\beta \theta }X_\theta ^{\gamma \alpha }-\frac 12\pi ^{\theta \gamma
}X_\theta ^{\beta \alpha }=0,
\]
because $\pi ^{\beta \gamma }/^\alpha =0,$ which ends the
proof.\hfill \hbox{\rlap{$\sqcap$}$\sqcup$}

The image of the anchor map $\sigma (E)\subseteq TM$ defines an
integrable smooth distribution on $M$. Therefore, the manifold $M$
is foliated by the integral leaves of $\sigma (E),$ which are
called the leaves of the Lie algebroid. A curve
$u:[t_0,t_1]\rightarrow E$ is called admissible if $ \sigma
(u(t))=\dot c(t)$, where $c(t)=\pi (u(t))$ is the base curve on
$M$. It follows that $u(t)$ is admissible if and only if the base
curve $c(t)$ lies on a leaf of the Lie algebroid whereas two
points can be joint by an admissible curve if and only if they are
situated on the same leaf. We can choose a smooth family
$t\rightarrow \theta (t)\in E^{*}$ of $1$-form such that
$\widetilde{\sigma }\theta (t)=\dot c(t)$. We shall call the pair
($ u(t),$ $\theta (t)$) the \textit{dual curve}.

\begin{definition}
Let ($u(t),$ $\theta (t)$) a \textit{dual curve} on $E$. We say that ($u(t),$
$\theta (t)$) is a geodesic if
\[
(D_\theta \theta )_{u(t)}=0.
\]
\end{definition}

In local coordinates we obtain that a curve
\[
(u(t),\theta (t))=(x^1(t),...,x^n(t),\theta _1(t),...,\theta _m(t))
\]
is geodesic if and only if it satisfies the following system of differential
equations
\begin{equation}
\left\{
\begin{array}{c}
\frac{dx^i(t)}{dt}=\sigma _\gamma ^i\pi ^{\alpha \gamma
}(x^1(t),...,x^n(t))\theta _\alpha (t) \\
\\
\frac{d\theta _\alpha (t)}{dt}=-\Gamma _\alpha ^{\beta \gamma
}(x^1(t),...,x^n(t))\theta _\beta \theta _\gamma .
\end{array}
\right.  \tag{1.7.8}
\end{equation}

\newpage\

\subsubsection{\textbf{Compatible Poisson structures}}

Let us consider the Poisson bivector on Lie algebroids given by the relation
(1.7.1). We obtain

\begin{proposition}
The complete lift of $\Pi $ on $\mathcal{T}E$ is given by
\begin{equation}
\Pi ^{\mathrm{c}}=\pi ^{\alpha \beta }\mathcal{X}_\alpha \wedge
\mathcal{V} _\beta +\left( \frac 12\sigma _\gamma ^i\frac{\partial
\pi ^{\alpha \beta }}{
\partial x^i}-\pi ^{\delta \beta }L_{\delta \gamma }^\alpha \right) y^\gamma
\mathcal{V}_\alpha \wedge \mathcal{V}_\beta .  \tag{1.7.9}
\end{equation}
\end{proposition}

\textbf{Proof}. Using the properties of the vertical and complete lifts we
obtain

\[
\begin{array}{l}
\Pi ^{\mathrm{c}}=(\frac 12\pi ^{\alpha \beta }s_\alpha \wedge
s_\beta )^{ \mathrm{c}}=(\frac 12\pi ^{\alpha \beta
})^{\mathrm{c}}(s_\alpha \wedge s_\beta )^{\mathrm{v}}+(\frac
12\pi ^{\alpha \beta })^{\mathrm{v}}(s_\alpha
\wedge s_\beta )^{\mathrm{c}}= \\
\quad \quad =\frac 12\stackrel{\cdot }{\pi }^{\alpha \beta
}s_\alpha ^{ \mathrm{v}}\wedge s_\beta ^{\mathrm{v}}+\frac 12\pi
^{\alpha \beta }(s_\alpha ^{\mathrm{c}}\wedge s_\beta
^{\mathrm{v}}+s_\alpha ^{\mathrm{v} }\wedge s_\beta
^{\mathrm{c}})=\frac 12\stackrel{\cdot }{\pi }^{\alpha \beta
}\mathcal{V}_\alpha \wedge \mathcal{V}_\beta + \\
\quad \quad +\frac 12\pi ^{\alpha \beta }\left(
(\mathcal{X}_\alpha -L_{\alpha \gamma }^\delta y^\gamma
\mathcal{V}_\delta )\wedge \mathcal{V} _\beta +\mathcal{V}_\alpha
\wedge (\mathcal{X}_\beta -L_{\beta \gamma
}^\delta y^\gamma )\mathcal{V}_\delta \right) = \\
\quad \quad =\pi ^{\alpha \beta }\mathcal{X}_\alpha \wedge
\mathcal{V}_\beta +\left( \frac 12\sigma _\gamma ^i\frac{\partial
\pi ^{\alpha \beta }}{
\partial x^i}-\pi ^{\delta \beta }L_{\delta \gamma }^\alpha \right) y^\gamma
\mathcal{V}_\alpha \wedge \mathcal{V}_\beta .
\end{array}
\]
\hfill\hbox{\rlap{$\sqcap$}$\sqcup$}

\begin{proposition}
The complete lift $\Pi ^{\mathrm{c}}$ is a Poisson bivector on
$\mathcal{T}E$ .
\end{proposition}

\textbf{Proof}. Using the relation (1.7.9) by straightforward
computation we obtain $[\Pi ^{\mathrm{c}},\Pi ^{\mathrm{c}}]=0$
which ends the proof.\hfill \hbox{\rlap{$\sqcap$}$\sqcup$}

\begin{proposition}
The Poisson structure $\Pi ^{\mathrm{c}}$ has the following property
\[
\Pi ^{\mathrm{c}}=-\mathcal{L}_C\Pi ^{\mathrm{c}},
\]
which means that $(\mathcal{T}E,\Pi ^{\mathrm{c}})$ is a homogeneous Poisson
manifold.
\end{proposition}

\textbf{Proof}. A direct computation in local coordinates yields

\[
\begin{array}{ll}
\mathcal{L}_{\mathcal{C}}\Pi ^c & =\mathcal{L}_{y^\varepsilon
\mathcal{V} _\varepsilon }\left( \pi ^{\alpha \beta
}\mathcal{X}_\alpha \wedge \mathcal{V }_\beta +(\frac 12\sigma
_\gamma ^i\frac{\partial \pi ^{\alpha \beta }}{
\partial x^i}-\pi ^{\delta \beta }L_{\delta \gamma }^\alpha )y^\gamma
\mathcal{V}_\alpha \wedge \mathcal{V}_\beta \right) \\
& =\mathcal{L}_{y^\varepsilon \mathcal{V}_\varepsilon }(\pi
^{\alpha \beta } \mathcal{X}_\alpha )\wedge \mathcal{V}_\beta -\pi
^{\alpha \beta }\mathcal{X} _\alpha \wedge \mathcal{V}_\beta
+\mathcal{L}_{y^\varepsilon \mathcal{V} _\varepsilon }(\frac
12\sigma _\gamma ^i\frac{\partial \pi ^{\alpha \beta }}{
\partial x^i}y^\gamma \mathcal{V}_\alpha )\wedge \mathcal{V}_\beta \\
& -\frac 12\sigma _\gamma ^i\frac{\partial \pi ^{\alpha \beta
}}{\partial x^i }y^\gamma \mathcal{V}_\alpha \wedge
\mathcal{V}_\beta -(\mathcal{L} _{y^\varepsilon
\mathcal{V}_\varepsilon }\pi ^{\delta \beta }L_{\delta \gamma
}^\alpha y^\gamma \mathcal{V}_\alpha )\wedge \mathcal{V}_\beta
+\pi ^{\delta \beta }L_{\delta \gamma }^\alpha y^\gamma
\mathcal{V}_\alpha \wedge
\mathcal{V}_\beta \\
& =-\pi ^{\alpha \beta }\mathcal{X}_\alpha \wedge
\mathcal{V}_\beta -\frac 12\sigma _\gamma ^i\frac{\partial \pi
^{\alpha \beta }}{\partial x^i} y^\gamma \mathcal{V}_\alpha \wedge
\mathcal{V}_\beta +\pi ^{\delta \beta }L_{\delta \gamma }^\alpha
y^\gamma \mathcal{V}_\alpha \wedge \mathcal{V}
_\beta \\
& =-\Pi ^c
\end{array}
\]
\hfill \hbox{\rlap{$\sqcap$}$\sqcup$}\\

\begin{definition}
Let us consider a Poisson bivector on $E$ given by (1.7.1) then the
horizontal lift of $\Pi $ to $\mathcal{T}E$ is the bivector defined by
\[
\Pi ^H=\frac 12\pi ^{\alpha \beta }(x)\delta _\alpha \wedge \delta _\beta .
\]
\end{definition}

\begin{proposition}
The horizontal lift $\Pi ^H$ is a Poisson bivector if and only if $\Pi $ is
a Poisson bivector on $E$ and
\begin{equation}
\mathcal{R}\left( (\pi ^{\#}\rho )^h,(\pi ^{\#}\theta )^h\right) =0,\quad
\forall \rho ,\theta \in \Gamma (E).  \tag{1.7.10}
\end{equation}
\end{proposition}

\textbf{Proof}. The Poisson condition $[\Pi ,\Pi ]=0$ leads to the relation
(1.7.2) and $[\Pi ^H,\Pi ^H]=0$ yields
\[
\sum_{(\varepsilon ,\delta ,\alpha )}\left( \pi ^{\alpha \beta
}\pi ^{\gamma \delta }L_{\beta \gamma }^\varepsilon +\pi ^{\alpha
\beta }\sigma _\beta ^i \frac{\partial \pi ^{\varepsilon \delta
}}{\partial x^i}\right) \delta _\varepsilon \wedge \delta _\alpha
\wedge \delta _\delta +\pi ^{\alpha \beta }\pi ^{\gamma \delta
}\mathcal{R}_{\beta \gamma }^\varepsilon \mathcal{V} _\varepsilon
\wedge \delta _\alpha \wedge \delta _\gamma =0,
\]
and $\pi ^{\alpha \beta }\pi ^{\gamma \delta }\mathcal{R}_{\beta
\gamma }^\varepsilon =0,$ which is the local expression of the
relation (1.7.10). \hfill\hbox{\rlap{$\sqcap$}$\sqcup$}

We recall that two Poisson structures are compatible if the
bivectors $ \omega _1$ and $\omega _2$ satisfy the condition
\[
\lbrack \omega _1,\omega _2]=0
\]

By straightforward computation in local coordinates we get:

\begin{proposition}
The Poisson bivector $\Pi ^H$ is compatible with the complete lift
$\Pi ^{ \mathrm{c}}$ if and only if the following relations hold
\[
\pi ^{r\beta }\pi ^{\alpha s}(\frac{\partial \mathcal{N}_r^\gamma }{\partial
y^s}-\frac{\partial \mathcal{N}_s^\gamma }{\partial y^r})-\pi ^{r\gamma }\pi
^{s\alpha }L_{sr}^\beta =0,
\]
\[
\pi ^{rs}\left( \delta _r(a^{\alpha \beta })-a^{l\alpha
}\frac{\partial N_r^\beta }{\partial y^l}+a^{l\beta
}\frac{\partial N_r^\alpha }{\partial y^l }\right. \mathcal{-}
\]
\[
-\pi ^{\theta \beta }R_{r\theta }^\alpha +(\pi ^{\varepsilon \beta
}L_{\varepsilon \gamma }^\theta -\pi ^{\varepsilon \theta }L_{\varepsilon
\gamma }^\beta )y^\gamma \frac{\partial N_r^\alpha }{\partial y^\theta }+
\]
\[
+\left. \sigma _r^i\frac{\partial \pi ^{\varepsilon \beta
}}{\partial x^i} y^\gamma L_{\varepsilon \gamma }^\alpha +\pi
^{\varepsilon \beta }L_{\varepsilon \gamma }^\alpha N_r^\gamma
-\pi ^{\varepsilon \beta }\sigma _r^i\frac{\partial L_{\varepsilon
\gamma }^\alpha }{\partial x^i}y^\gamma \right) =0.
\]
\end{proposition}

where we have denoted
\[
a^{\alpha \beta }=\sigma _\varepsilon ^i\frac{\partial \pi
^{\alpha \beta }}{
\partial x^i}y^\varepsilon +N_\varepsilon ^\alpha \pi ^{\varepsilon \beta
}-N_\varepsilon ^\beta \pi ^{\varepsilon \alpha }.
\]

See \cite{Mit} for the particular case of tangent bundle.

\newpage\

\subsubsection{\textbf{Canonical Poisson structure}}

On Lie algebroid $(\mathcal{T}E^{*},[\cdot ,\cdot \pi
]_{\mathcal{T} E^{*}},\sigma ^1)$ we have the canonical symplectic
section $\omega _E$ given by (1.5.6) which induces a vector bundle
isomorphism
\[
\natural _{\omega _E}:E^{*}\rightarrow E,\quad i_\zeta \omega _E\in
E^{*}\rightarrow \zeta \in E.
\]

\begin{definition}
The canonical Poisson bivector is given by
\[
\Lambda =\natural _{\omega _E}\omega _E.
\]
\end{definition}

It follows that
\[
\Lambda (dF,dG)=-\omega _E(\natural (dF),\natural (dG)),\quad F,G\in
C^\infty (E^{*})
\]
and in local coordinates we get
\[
\Lambda =\mathcal{P}^\alpha \wedge \mathcal{X}_\alpha +\frac 12\mu _\alpha
L_{\beta \gamma }^\alpha \mathcal{P}^\beta \wedge \mathcal{P}^\gamma .
\]

\begin{remark}
The Schouten-Nijenhuis bracket $[\Lambda ,\Lambda ]$ leads, locally, to the
expression
\[
\frac 13\sum_{(\alpha ,\beta ,\gamma )}(\sigma _\alpha
^i\frac{\partial L_{\beta \gamma }^\varepsilon }{\partial
x^i}+L_{\alpha \delta }^\varepsilon L_{\beta \gamma }^\delta )\mu
_\varepsilon \mathcal{P}^\beta \wedge \mathcal{ P}^\alpha \wedge
\mathcal{P}^\gamma
\]
and $[\Lambda ,\Lambda ]=0$ follows from the structure equations on Lie
algebroids (1.2.9).
\end{remark}

\begin{definition}
Let us consider a Poisson bivector on $E$ given by
\[
\Pi =\frac 12\pi ^{\alpha \beta }(x)s_\alpha \wedge s_\beta
\]
then the horizontal lift of $\Pi $ to $\mathcal{T}E^{*}$ is the bivector
defined by
\[
\Pi ^H=\frac 12\pi ^{\alpha \beta }(x)\delta _\alpha ^{*}\wedge \delta
_\beta ^{*}.
\]
\end{definition}

\begin{proposition}
The horizontal lift $\Pi ^H$ is a Poisson bivector if and only if $\Pi $ is
a Poisson bivector on $E$ and
\begin{equation}
\mathcal{R}\left( (\pi ^{\#}\rho )^h,(\pi ^{\#}\theta )^h\right) =0,\quad
\forall \rho ,\theta \in \Gamma (E^{*})  \tag{1.7.11}
\end{equation}
\end{proposition}

\textbf{Proof}. The Poisson condition $[W,W]=0$ leads to the relation
\[
\sum_{(\alpha ,\varepsilon ,\delta )}(\pi ^{\alpha \beta }\pi
^{\gamma \delta }L_{\beta \gamma }^\varepsilon +\pi ^{\alpha \beta
}\sigma _\beta ^i \frac{\partial \pi ^{\varepsilon \delta
}}{\partial x^i})=0
\]
and $[\Pi ^H,\Pi ^H]=0$ yields
\[
\sum_{(\varepsilon ,\delta ,\alpha )}\left( \pi ^{\alpha \beta
}\pi ^{\gamma \delta }L_{\beta \gamma }^\varepsilon +\pi ^{\alpha
\beta }\sigma _\beta ^i \frac{\partial \pi ^{\varepsilon \delta
}}{\partial x^i}\right) \delta _\varepsilon ^{*}\wedge \delta
_\alpha ^{*}\wedge \delta _\delta ^{*}+\pi ^{\alpha \beta }\pi
^{\gamma \delta }\mathcal{R}_{\beta \gamma \varepsilon }
\mathcal{P}^\varepsilon \wedge \delta _\alpha ^{*}\wedge \delta
_\gamma ^{*}=0
\]
and it results $\pi ^{\alpha \beta }\pi ^{\gamma \delta }\mathcal{R}_{\beta
\gamma \varepsilon }=0,$ which is the local expression of the condition
(1.7.11).\hfill \hbox{\rlap{$\sqcap$}$\sqcup$}\\

\begin{proposition}
If the connection $\mathcal{N}$ on $\mathcal{T}E^{*}$ is defined by a linear
connection $\nabla $ with the coefficients $\Gamma _{\alpha \beta }^\gamma $
on the Lie algebroid $E$, the the bivector $\Pi ^H$ has the following form
\begin{equation}
\Pi ^H=\frac 12\pi ^{\alpha \beta }\mathcal{Q}_\alpha \wedge
\mathcal{Q} _\beta +\frac 12\pi ^{\alpha \beta }\mu _\gamma \mu
_\theta \Gamma _{\alpha \varepsilon }^\gamma \Gamma _{\beta \delta
}^\theta \mathcal{P}^\varepsilon \wedge \mathcal{P}^\delta +\pi
^{\alpha \beta }\mu _\gamma \Gamma _{\beta \varepsilon }^\gamma
\mathcal{Q}_\alpha \wedge \mathcal{P}^\varepsilon . \tag{1.7.12}
\end{equation}
\end{proposition}

\textbf{Proof}. The coefficients of the nonlinear connection have
the form $ N_{\alpha \beta }=\mu _\gamma \Gamma _{\alpha \beta
}^\gamma $ and introducing the relation $\delta _\alpha
^{*}=\mathcal{Q}_\alpha +\mathcal{N} _{\alpha \beta
}\mathcal{P}^\beta $ into the expression of $\Pi ^H$ it results
(1.7.12).

\begin{proposition}
If $\mathcal{N}$ is a symmetric nonlinear connection then the canonical
Poisson bivector has the form
\[
\Lambda =\mathcal{P}^\alpha \wedge \delta _\alpha ^{*}.
\]
\end{proposition}

\textbf{Proof}. We have
\[
\begin{array}{l}
\Lambda =\mathcal{P}^\alpha \wedge \mathcal{Q}_\alpha +\frac 12\mu
_\alpha L_{\beta \gamma }^\alpha \mathcal{P}^\beta \wedge
\mathcal{P}^\gamma = \mathcal{P}^\alpha \wedge (\delta _\alpha
^{*}-\mathcal{N}_{\alpha \beta } \mathcal{P}^\beta )+\frac 12\mu
_\gamma L_{\beta \alpha }^\gamma \mathcal{P}
^\beta \wedge \mathcal{P}^\alpha \\
\quad =\mathcal{P}^\alpha \wedge \delta _\alpha ^{*}-\frac 12\mathcal{N}
_{\alpha \beta }\mathcal{P}^\alpha \wedge \mathcal{P}^\beta -\frac 12\left(
\mathcal{N}_{\alpha \beta }+\mu _\gamma L_{\beta \alpha }^\gamma \right)
\mathcal{P}^\alpha \wedge \mathcal{P}^\beta \\
\quad =\mathcal{P}^\alpha \wedge \delta _\alpha ^{*}-\frac
12\mathcal{N} _{\alpha \beta }\mathcal{P}^\alpha \wedge
\mathcal{P}^\beta -\frac 12
\mathcal{N}_{\beta \alpha }\mathcal{P}^\alpha \wedge \mathcal{P}^\beta \\
\quad =\mathcal{P}^\alpha \wedge \delta _\alpha ^{*}-\frac
12\mathcal{N} _{\alpha \beta }\mathcal{P}^\alpha \wedge
\mathcal{P}^\beta -\frac 12
\mathcal{N}_{\alpha \beta }\mathcal{P}^\beta \wedge \mathcal{P}^\alpha \\
\quad =\mathcal{P}^\alpha \wedge \delta _\alpha ^{*}-\frac
12\mathcal{N} _{\alpha \beta }\mathcal{P}^\alpha \wedge
\mathcal{P}^\beta +\frac 12 \mathcal{N}_{\alpha \beta
}\mathcal{P}^\alpha \wedge \mathcal{P}^\beta = \mathcal{P}^\alpha
\wedge \delta _\alpha ^{*}.\quad \hfill
\hbox{\rlap{$\sqcap$}$\sqcup$}
\end{array}
\]

\begin{proposition}
If $\Pi ^H$ is a Poisson bivector and $\mathcal{N}$ is a symmetric nonlinear
connection, then $\Pi ^H$ is compatible with the canonical Poisson structure
$\Lambda $ if and only if the following relations fulfilled
\begin{equation}
\sigma _\gamma ^i\frac{\partial \pi ^{\alpha \beta }}{\partial
x^i}+\pi ^{\varepsilon \alpha }\left( \frac{\partial
\mathcal{N}_{\varepsilon \gamma } }{\partial \mu _\beta
}-L_{\varepsilon \gamma }^\beta \right) -\pi ^{\varepsilon \beta
}\left( \frac{\partial \mathcal{N}_{\varepsilon \gamma }
}{\partial \mu _\alpha }-L_{\varepsilon \gamma }^\alpha \right)
=0, \tag{1.7.13}
\end{equation}
\begin{equation}
\pi ^{\alpha \beta }\mathcal{R}_{\alpha \gamma \varepsilon }=0.  \tag{1.7.14}
\end{equation}
\end{proposition}

\textbf{Proof}. If $\mathcal{N}$ is symmetric then $\mathcal{N}_{\alpha
\beta }-\mathcal{N}_{\beta \alpha }=\mu _\gamma L_{\alpha \beta }^\gamma $
and with respect with the basis \{$\delta _\alpha ^{*},\mathcal{P}^\alpha \}$
it results $\Lambda =\mathcal{P}^\alpha \wedge \delta _\alpha ^{*}.$ By a
straightforward computation we obtain
\begin{eqnarray*}
\lbrack \Pi ^H,\Lambda ] &=&-\frac 12\left( \sigma _\gamma
^i\frac{\partial \pi ^{\alpha \beta }}{\partial x^i}+\pi
^{\varepsilon \alpha }\left( \frac{
\partial \mathcal{N}_{\varepsilon \gamma }}{\partial \mu _\beta }
-L_{\varepsilon \gamma }^\beta \right) \right) \delta _\alpha ^{*}\wedge
\delta _\beta ^{*}\wedge \mathcal{P}^\gamma \\
&&\ \ +\frac{\pi ^{\varepsilon \beta }}2\left( \frac{\partial
\mathcal{N} _{\varepsilon \gamma }}{\partial \mu _\alpha
}-L_{\varepsilon \gamma }^\alpha \right) \delta _\alpha ^{*}\wedge
\delta _\beta ^{*}\wedge \mathcal{
P}^\gamma \\
&&\ \ +\pi ^{\alpha \beta }\mathcal{R}_{\alpha \gamma \varepsilon
}\mathcal{P }^\varepsilon \wedge \delta _\beta ^{*}\wedge
\mathcal{P}^\gamma .
\end{eqnarray*}
and $[\Pi ^H,\Lambda ]=0$ is equivalent with the relations (1.7.13),
(1.7.14).\hfill\hbox{\rlap{$\sqcap$}$\sqcup$}

\begin{remark}
If the nonlinear connection $\mathcal{N}$ is defined by a linear
connection $ \nabla $ with the coefficients $\Gamma _{\alpha \beta
}^\gamma $ on the Lie algebroid $E$ then we obtain the conditions
\[
\sigma _\gamma ^i\frac{\partial \pi ^{\alpha \beta }}{\partial x^i}+\pi
^{\varepsilon \alpha }\left( \Gamma _{\varepsilon \gamma }^\beta
-L_{\varepsilon \gamma }^\beta \right) -\pi ^{\varepsilon \beta }\left(
\Gamma _{\varepsilon \gamma }^\alpha -L_{\varepsilon \gamma }^\alpha \right)
=0,
\]
\[
\pi ^{\alpha \beta }\mu _\varepsilon \left( \sigma _\alpha
^i\frac{\partial \Gamma _{\beta \gamma }^\varepsilon }{\partial
q^i}-\sigma _\beta ^i\frac{
\partial \Gamma _{\alpha \gamma }^\varepsilon }{\partial q^i}+\Gamma
_{\alpha \theta }^\varepsilon \Gamma _{\beta \gamma }^\theta -\Gamma _{\beta
\theta }^\varepsilon \Gamma _{\alpha \gamma }^\theta -L_{\alpha \beta
}^\theta \Gamma _{\theta \gamma }^\varepsilon \right) =0.
\]
\end{remark}

\newpage\

\section{\textbf{Optimal Control }}

\quad

The Lie geometric methods in control theory have been applied by many
authors (see for instance. \cite{Br, Ju, Bl, Ma3}). One of the most
important issues in the geometric approach is the analysis of the solution
to the optimal control problem as provided by Pontryagin's Maximum
Principle; that is, the curve $c(t)=(x(t),u(t))$ is an optimal trajectory if
there exists a lifting of $x(t)$ to the dual space $(x(t),p(t))$ satisfying
the Hamilton equations, together with a maximization condition for the
Hamiltonian with respect to the control variables $u(t)$. The purpose of
this part is to study the drift less control affine systems (distributional
systems) with positive homogeneous cost, using the Pontryagin Maximum
Principle at the level of a Lie algebroid in the case of constant rank of
distribution. Author's papers \cite{Hr3, Po12, Po14, Po17, Po18, Po22, Po24}
are used in writting this part.

We prove that the framework of Lie algebroids is better than cotangent
bundle in order to solve some problems of drift less control affine systems.
In the first chapter the known results on the optimal control systems are
recalled by geometric viewpoint. In the next chapter the distributional
systems are presented and the relation between the Hamiltonians on $E^{*}$
and $T^{*}M$ is given. We investigate the cases of holonomic and
nonholonomic distributions with constant rank.

In the holonomic case, we will consider the Lie algebroid being
just the distribution whereas in the nonholonomic case (i.e.,
strong bracket generating distribution) the Lie algebroid is the
tangent bundle with the basis given by vectors of distribution
completed by the first Lie brackets. In the both cases
illustrative examples are presented. Also, the case of
distribution $D$ with non-constant rank is studied and some
interesting examples are given. In the last chapter we present the
intrinsic relation between the distributional systems and
sub-Riemannian geometry. Thus, the optimal trajectory of our
distributional systems are the geodesics in the framework of
sub-Riemannian geometry. We investigate two classical cases:
Grusin plan and Heisenberg group \cite{Hr3}, equipped with
positive homogeneous costs (Randers metric \cite{Bao}). \newpage\

\subsection{\textbf{Geometric viewpoint of the optimal control}}

Let $M$ be a smooth $n$-dimensional manifold. We consider the control system
\begin{eqnarray*}
\frac{dx^i}{dt}=f^i(x,u),
\end{eqnarray*}
where $x\in M$ and the control $u$ takes values in an open subset $\Omega $
of ${R}^m$. Let $x_0$ and $x_1$ be two points of $M$. An optimal control
problem consists of finding the trajectories of our control system which
connects $x_0$ and $x_1$ and minimizing the cost
\begin{eqnarray*}
\min \int_0^TL(x(t),u(t))dt,\quad x(0)=x_0,\ x(T)=x_1,
\end{eqnarray*}
where $L$ is the \textit{Lagrangian} or \textit{running cost}.\\Necessary
conditions for a trajectory to be an extreme are given by Pontryagin Maximum
Principle. The Hamiltonian reads as
\begin{eqnarray*}
H(x,p,u)=\left\langle p,f(x,u))\right\rangle -L(x,u),\quad p\in T^{*}M,
\end{eqnarray*}
while the maximization condition with respect to the control variables $u$,
namely
\begin{eqnarray*}
H(x(t),p(t),u(t))=\max_vH(x(t),p(t),v),
\end{eqnarray*}
leads to
\begin{eqnarray*}
\frac{\partial H}{\partial u}=0.
\end{eqnarray*}
The extreme trajectories satisfy the Hamilton equations
\begin{equation}
\dot x=\frac{\partial H}{\partial p},\quad \dot p=-\frac{\partial
H}{
\partial x}.  \tag{2.1.1}
\end{equation}
From a geometric viewpoint the pair $(x^i,u^a)$, where
$i=\overline{1,n}$ and $a=\overline{1,m}$, can be understood as a
local coordinate pair of a manifold $E,$ that is fibered over $M$
by the projection $\pi :E\rightarrow M $. The functions $f^i(x,u)$
are the components of a vector field $ X=f^i(x,u)\frac \partial
{\partial x^i}$ along $\pi $, that is, of a fibered mapping
$X:E\rightarrow TM$ from the bundle $(E,\pi ,M)$ to the tangent
bundle $(TM,\tau ,M)$ such that $\tau \circ X=\pi $. The
admissible curves of the control system are the curves $\gamma
:I\subset {R}\rightarrow E$ such that
\begin{eqnarray*}
X(\gamma (t))=\frac d{dt}(\pi (\gamma (t)).
\end{eqnarray*}
The optimal control problem consists of obtaining the admissible curves that
minimize the cost
\begin{eqnarray*}
\min \int_IL(\gamma (t))dt,\quad L\in C^\infty (E),
\end{eqnarray*}
and satisfy certain boundary conditions not to be considered
here.\\The Hamiltonian $H$ is a real-valued function defined on
the fibered product $ T^{*}M\times _ME$ that is given by
\begin{eqnarray*}
H(\mu ,v)=\left\langle \mu ,X(v)\right\rangle -L(v),
\end{eqnarray*}
for any $(\mu ,v)\in T^{*}M\times _ME.$\\The critical equations
follow from asking a vector field $X_H$ defined along a map
$pr_1:T^{*}M\times _ME\rightarrow T^{*}M$ to satisfy the
symplectic equations
\begin{eqnarray*}
i_{X_H}\omega =dH,
\end{eqnarray*}
where $\omega =dx^i\wedge dp_i$. Since
\begin{eqnarray*}
dH=\frac{\partial H}{\partial x^i}dx^i+\frac{\partial H}{\partial
p_i}dp_i+ \frac{\partial H}{\partial u^a}du^a,
\end{eqnarray*}
we obtain that the solution of the above equations is the vector field
\begin{eqnarray*}
X_H=\frac{\partial H}{\partial p_i}\frac \partial {\partial
x^i}-\frac{
\partial H}{\partial x^i}\frac \partial {\partial p_i},
\end{eqnarray*}
defined on the subset
\[
\frac{\partial H}{\partial u^a}=0
\]
of $T^{*}M\times _ME$, and therefore, the critical trajectories are the
integral curves of the above vector field, namely
\begin{equation}
\frac{\partial H}{\partial u^a}=0,\quad \dot {x^i}=\frac{\partial
H}{
\partial p_i},\quad \dot {p_i}=-\frac{\partial H}{\partial x^i}.  \tag{2.1.2}
\end{equation}

By a control system on the Lie algebroid (see \cite{Ma3}) $\pi
:E\rightarrow M$ with the control space $\tau :A\rightarrow M$ we
mean a section $\rho $ of $E$ along $\tau $. A trajectory of the
system $\rho $ is an integral curve of the vector field $\sigma
(\rho )$. Given the cost function $ \mathcal{L}\in C^\infty (A)$,
we have to minimize the integral of $\mathcal{L }$ over the set of
those system trajectories which satisfy certain boundary
conditions. The Hamiltonian function $\mathcal{H}\in C^\infty
(E^{*}\times _MA)$ is defined by
\begin{eqnarray*}
\mathcal{H}(\mu ,u)=\left\langle \mu ,\rho (u)\right\rangle -\mathcal{L}(u),
\end{eqnarray*}
whereas the associated Hamiltonian control system $\rho _H$ is given by the
symplectic equation
\begin{eqnarray*}
i_{\rho _H}\omega _E=d^E\mathcal{H}.
\end{eqnarray*}
In local coordinates, the solution of the previous equation reads as
\begin{eqnarray*}
\rho _H=\frac{\partial \mathcal{H}}{\partial \mu _\alpha }\mathcal{X}_\alpha
-\left( \sigma _\alpha ^i\frac{\partial \mathcal{H}}{\partial q^i}+\mu
_\gamma L_{\alpha \beta }^\gamma \frac{\partial \mathcal{H}}{\partial \mu
_\beta }\right) \mathcal{P}^\alpha ,
\end{eqnarray*}
on the subset where
\[
\frac{\partial \mathcal{H}}{\partial u^A}=0.
\]
Therefore, the critical trajectories are given by
\begin{equation}
\frac{\partial \mathcal{H}}{\partial u^A}=0,\quad
\frac{dq^i}{dt}=\sigma _\alpha ^i\frac{\partial
\mathcal{H}}{\partial \mu _\alpha },\quad \frac{ d\mu _\alpha
}{dt}=-\sigma _\alpha ^i\frac{\partial \mathcal{H}}{\partial q^i
}-\mu _\gamma L_{\alpha \beta }^\gamma \frac{\partial
\mathcal{H}}{\partial \mu _\beta }.  \tag{2.1.3}
\end{equation}

\newpage\

\subsection{\textbf{Distributional systems}}

Let $M$ be a smooth $n$-dimensional manifold. We consider the distributional
system (drift less control-affine system)
\begin{equation}
\dot x=\sum_{i=1}^mu_iX_i(x),  \tag{2.2.1}
\end{equation}
where $x\in M\,$, $X_1,X_2,...,X_m$ are smooth vector fields on
$M$ and the control $u=(u_1,u_2,...,u_m)$ takes values in an open
subset $\Omega $ of ${R }^m$. The vector fields
$X_i,i=\overline{1,m}$, generate a distribution $ D\subset TM$
such that the rank of $D$ is constant. Let $x_0$ and $x_1$ be two
points of $M$. An optimal control problem consists of finding
those trajectories of the distributional system which connect
$x_0$ and $x_1,$ while minimizing the cost
\begin{equation}
\min_{u(\cdot )}\int_I\mathcal{F}(u(t))dt,  \tag{2.2.2}
\end{equation}
where $\mathcal{F}$ is a Minkowski norm (positive homogeneous) on $D$.

\begin{remark}
We can associate to any positive homogeneous cost $\mathcal{F}$ on the Lie
algebroid $E$ a cost $F$ on $Im\sigma \subset TM$ defined by
\begin{equation}
F(v)=\{\mathcal{F}(u)\left| u\in E_x,\quad \sigma (u)=v\right. \},
\tag{2.2.3}
\end{equation}
where $v\in (Im\sigma )_x\subset T_xM$, $x\in M$.
\end{remark}

A piecewise smooth curve $c:I\subset {R}\rightarrow M$ is called horizontal
if the tangent vectors are in $D$, i.e. $\dot c(t)\in D_{c(t)}\subset TM$
for almost every $t\in I$. Let $u:I\rightarrow E$ be an admissible curve
projected by $\pi $ onto the horizontal curve $c:I\rightarrow M.$ The length
of the horizontal curve $c$ is defined by
\begin{eqnarray*}
length(c)=\int_I\mathcal{F}(u(t))dt=\int_IF(\dot c(t))dt,
\end{eqnarray*}
and the distance is given by $d(a,b)=\inf length(c)$ where the infimum is
taken over all the horizontal curves connecting $a$ and $b$. The distance is
infinite if there is no admissible curve that connects these two points.

\begin{remark}
The energy of a horizontal curve is
\begin{eqnarray*}
E(c)=\frac 12\int_IF^2(\dot c(t))dt,
\end{eqnarray*}
and it can easily be proved that if a curve is parameterized to a
constant speed, then it minimizes the length integral if and only
if it minimizes the energy integral.
\end{remark}

For the 2-homogeneous Lagrangian $L=\frac 12F^2$ and
$\mathcal{L}=\frac 12 \mathcal{F}^2$ we have
\[
\mathcal{L}=L\circ \sigma .
\]
Further, $L$ on $TM$ is a Lagrangian with constraints. According to the
Pontryagin Maximum Principle, the Hamiltonian is recast as
\begin{eqnarray*}
H(x,p,u)=\left\langle p,\dot x\right\rangle -L(x,u).
\end{eqnarray*}
If the equations
\begin{eqnarray*}
\frac{\partial H(x,p,u)}{\partial u}=0,
\end{eqnarray*}
permit us to find in a unique way $u$ as a smooth function of $(x,p)$ then
we can write the Hamiltonian system without any dependence on the control.
This nice situation happens always for distributional systems with quadratic
cost
\begin{eqnarray*}
\min \int_I\sum_{i=1}^mu_i^2dt.
\end{eqnarray*}
If the cost is not quadratic, then we cannot guarantee that the Hamiltonian
can be calculated without dependence on the control. However, there exist
several situations when the Hamiltonian can still be found.

\begin{proposition}
The relation between the Hamiltonian $H$ on cotangent bundle $T^{*}M$ and
the Hamiltonian $\mathcal{H}$ on dual bundle $E^{*}$ is given by
\begin{equation}
H(p)=\mathcal{H}(\mu ),\quad \mu =\sigma ^{\star }(p),\quad p\in
T_x^{*}M,\quad \mu \in E_x^{*}.  \tag{2.2.4}
\end{equation}
\end{proposition}

\textbf{Proof.} The Fenchel-Legendre dual of Lagrangian $L$ is the
Hamiltonian $H$ given by
\begin{eqnarray*}
H(p) &=&\sup_v\left\{ \left\langle p,v\right\rangle -L(v)\right\}
={\sup_v}
\left\{ \left\langle p,v\right\rangle -\mathcal{L}(u);\sigma (u)=v\right\} \\
\ &=&\sup_u\left\{ \left\langle p,\sigma (u)\right\rangle
-\mathcal{L} (u)\right\} ={\sup_u}\left\{ \left\langle \sigma
^{\star }(p),u\right\rangle -\mathcal{L}(u)\right\}
=\mathcal{H}(\sigma ^{\star }(p)),
\end{eqnarray*}
and we get
\[
H(p)=\mathcal{H}(\mu ),\mu =\sigma ^{\star }(p),p\in T_x^{*}M,\mu \in
E_x^{*},
\]
or locally
\begin{equation}
\mu _\alpha =\sigma _\alpha ^ip_i,  \tag{2.2.5}
\end{equation}
where the Hamiltonian $H(p)$ is degenerate on $Ker\sigma ^{\star }\subset
T^{*}M.$\hfill \hbox{\rlap{$\sqcap$}$\sqcup$}\newpage\

\subsubsection{\textbf{Holonomic distribution}}

We assume for the beginning that the distribution $D=\left\langle
X_1,X_2,...,X_m\right\rangle $ is holonomic, which means that
$[X_i,X_j]\in D $ for every $i,j=\overline{1,m}$, $i\neq j.$ In
order to apply the theory of Lie algebroids we consider $E=D$ with
the inclusion as anchor $\sigma :E\rightarrow TM$. From the
Frobenius theorem, the distribution $D$ is integrable, it
determines a foliation on $M$ and two points can be joined if and
only if they are situated on the same leaf.\\We consider the
following distributional system with positive homogeneous cost
\cite{Po14}:
\begin{eqnarray*}
\dot x=u^1X_1+u^2X_2,\quad x=\left(
\begin{array}{c}
x^1 \\
x^2 \\
x^3
\end{array}
\right) \in {R}^3,\ X_1=\left(
\begin{array}{c}
1 \\
0 \\
0
\end{array}
\right) ,\ X_2=\left(
\begin{array}{c}
x^1 \\
x^2 \\
1
\end{array}
\right) ,
\end{eqnarray*}
and
\begin{eqnarray*}
\min_{u(\cdot )}\int_0^T\mathcal{F}(u(t))dt,\quad
\mathcal{F}(u)=\sqrt{ (u^1)^2+(u^2)^2}+\varepsilon u^1,\quad 0\leq
\varepsilon <1.
\end{eqnarray*}
We are looking for the trajectories starting from the point
$(1,1,0)^t$ and parameterized by arclength. The associated
distribution $D=\left\langle X_1,X_2\right\rangle $ is holonomic,
because
\[
X_1=\frac \partial {\partial x^1},\quad X_2=x^1\frac \partial {\partial
x^1}+x^2\frac \partial {\partial x^2}+\frac \partial {\partial x^3},
\]
and therefore $[X_1,X_2]=X_1$. In the case of the Lie algebroid,
we consider $E=\left\langle X_1,X_2\right\rangle $ and the anchor
$\sigma :E\rightarrow T \Bbb{R}^3$ has the components
\begin{equation}
\sigma _\alpha ^i=\left(
\begin{array}{cc}
1 & x^1 \\
0 & x^2 \\
0 & 1
\end{array}
\right) ,  \tag{2.2.6}
\end{equation}
and we get the Lagrangian
\[
\mathcal{L}=\frac 12\left( \sqrt{(u_1)^2+(u_2)^2}+\varepsilon u_1\right) ^2.
\]
Using \cite{Hr2},\cite{Mi5} we can find the Hamiltonian on $E^{*}$ given by
\begin{equation}
\mathcal{H}(\mu )=\frac 12\left( \sqrt{\frac{(\mu
_1)^2}{(1-\varepsilon ^2)^2 }+\frac{(\mu _2)^2}{1-\varepsilon
^2}}-\frac{\varepsilon \mu _1}{ 1-\varepsilon ^2}\right) ^2.
\tag{2.2.7}
\end{equation}

\begin{remark}
Using (2.2.4) we can calculate the Hamiltonian $H$ on $T^{*}M$
given by $ H(x,p)=\mathcal{H}(\mu ),$ $\mu =\sigma ^{\star }(p)$,
where
\begin{eqnarray*}
\left(
\begin{array}{c}
\mu _1 \\
\mu _2
\end{array}
\right) =\left(
\begin{array}{ccc}
1 & 0 & 0 \\
x^1 & x^2 & 1
\end{array}
\right) \left(
\begin{array}{c}
p_1 \\
p_2 \\
p_3
\end{array}
\right) .
\end{eqnarray*}
We get that
\begin{equation}
H(x,p)=\frac 12\left( \sqrt{\frac{(p_1)^2}{\left( 1-\varepsilon
^2\right) ^2} +\frac{(p_1x^1+p_2x^2+p_3)^2}{1-\varepsilon
^2}}-\frac{\varepsilon p_1}{ 1-\varepsilon ^2}\right) ^2.
\tag{2.2.8}
\end{equation}
Unfortunately, with $H(x,p)$ from (2.2.8) the Hamilton equations on $T^{*}M$
lead to a complicated system of implicit differential equations.
\end{remark}

We will use the geometric model of a Lie algebroid. From the
relation $ [X_\alpha ,X_\beta ]=L_{\alpha \beta }^\gamma X_\gamma
$ we obtain the non-zero components $L_{12}^1=1,$ $L_{21}^1=-1$
while from (2.1.3) we deduce that
\begin{equation}
\dot x^1 =\frac{\partial \mathcal{H}}{\partial \mu
_1}+x^1\frac{\partial \mathcal{H}}{\partial \mu _2},\ \dot
x^2=x^2\frac{\partial \mathcal{H}}{
\partial \mu _2},\ \dot x^3=\frac{\partial \mathcal{H}}{\partial \mu _2},
\tag{2.2.9}
\end{equation}
\begin{equation*}
\dot \mu _1 =-\mu _1\frac{\partial \mathcal{H}}{\partial \mu _2},\
\dot \mu _2=\mu _1\frac{\partial \mathcal{H}}{\partial \mu _1},
\end{equation*}
where
\begin{eqnarray*}
\frac{\partial \mathcal{H}}{\partial \mu _1}=\frac{\left(
1+\varepsilon ^2\right) \mu _1}{(1-\varepsilon
^2)^2}-\frac{\varepsilon \sqrt{\frac{(\mu _1)^2}{(1-\varepsilon
^2)^2}+\frac{(\mu _2)^2}{1-\varepsilon ^2}}}{ 1-\varepsilon
^2}-\frac{\varepsilon \mu _1^2}{(1-\varepsilon ^2)^3\sqrt{
\frac{(\mu _1)^2}{(1-\varepsilon ^2)^2}+\frac{(\mu
_2)^2}{1-\varepsilon ^2}}} ,
\end{eqnarray*}
\begin{equation}
\frac{\partial \mathcal{H}}{\partial \mu _2}=\frac{\mu
_2}{1-\varepsilon ^2}- \frac{\varepsilon \mu _1\mu
_2}{(1-\varepsilon ^2)^2\sqrt{\frac{(\mu _1)^2}{ (1-\varepsilon
^2)^2}+\frac{(\mu _2)^2}{1-\varepsilon ^2}}}.  \tag{2.2.10}
\end{equation}
The form of the last relations leads to the following change of variables
\begin{eqnarray*}
\mu _1(t) &=&(1-\varepsilon ^2)r(t)\,\text{sech}\theta (t),\quad \\
&& \\
\mu _2(t) &=&\sqrt{1-\varepsilon ^2}r(t)\tanh \theta (t).
\end{eqnarray*}
In these circumstances we have
\[
\sqrt{\frac{(\mu _1)^2}{(1-\varepsilon ^2)^2}+\frac{(\mu _2)^2}{
1-\varepsilon ^2}}=\left| r\right| ,
\]
whereas
\[
\dot \mu _1=-\mu _1\frac{\partial \mathcal{H}}{\partial \mu _2},
\]
yields
\begin{equation}
\sqrt{1-\varepsilon ^2}\left( \frac{\dot r}r-\dot \theta \tanh \theta
\right) =r(-\tanh \theta +\varepsilon \sec h\theta \tanh \theta ).
\tag{2.2.11}
\end{equation}
From
\[
\dot \mu _2=\mu _1\frac{\partial \mathcal{H}}{\partial \mu _1},
\]
we get
\begin{equation}
\sqrt{1-\varepsilon ^2}\left( \frac{\dot r}r\tanh \theta +\dot
\theta \text{ sech}^2\theta \right) =r((1+\varepsilon
)^2\text{sech}^2\theta -\varepsilon \text{sech}\theta -\varepsilon
\text{sech}^3\theta ).  \tag{2.2.12}
\end{equation}
Now, reducing $\dot \theta $ and $\frac{\dot r}r$ from the relations
(2.2.11) and (2.2.12), we obtain
\begin{eqnarray*}
\sqrt{1-\varepsilon ^2}\dot r=r^2\varepsilon \text{sech}\theta \tanh \theta
(\varepsilon \text{sech}\theta -1),
\end{eqnarray*}
and
\begin{equation}
\sqrt{1-\varepsilon ^2}\dot \theta =r(\varepsilon \text{sech}\theta -1)^2.
\tag{2.2.13}
\end{equation}
The last two relations lead to
\begin{eqnarray*}
\frac{\dot r}{\dot \theta }=\frac{r\varepsilon \text{sech}\theta \tanh
\theta }{\varepsilon \text{sech}\theta -1},
\end{eqnarray*}
and respectively to
\begin{eqnarray*}
\frac 1rdr=\frac{\varepsilon \text{sech}\theta \tanh \theta }{\varepsilon
\text{sech}\theta -1}d\theta ,
\end{eqnarray*}
with the solution
\[
\ln \left| r\right| =-\ln (\varepsilon \sec h\theta -1)-\ln c.
\]
Therefore
\begin{eqnarray*}
\left| r\right| =\frac 1{c\left( \varepsilon \text{sech}\theta -1\right) }.
\end{eqnarray*}
Since the geodesics are parameterized by arclength, the conclusion
corresponds exactly to the $1/2$ level of the Hamiltonian and so we have
\begin{eqnarray*}
\mathcal{H}=\frac{r^2}2(1-\varepsilon \text{sech}\theta )^2=\frac 1{2c^2}.
\end{eqnarray*}
Now, $c=\pm 1$ and
\begin{eqnarray*}
r=\pm \frac 1{\varepsilon \text{sech}\theta -1}.
\end{eqnarray*}
From the relation (2.2.13) we have
\begin{eqnarray*}
\frac{d\theta }{dt}=\frac{\sqrt{1-\varepsilon ^2}}{1-\varepsilon
\text{sech} \theta },
\end{eqnarray*}
and respectively
\begin{eqnarray*}
t=\sqrt{1-\varepsilon ^2}\int \frac 1{1-\varepsilon \text{sech}\theta
}d\theta ,
\end{eqnarray*}
The relation
\[
\dot \mu _1=-\mu _1\dot x^3,
\]
implies that
\begin{eqnarray*}
x^3(\theta )=\ln \frac{c_1(1-\varepsilon \text{sech}\theta )}{(1-\varepsilon
^2)\text{sech}\theta },\quad c_1\in {R}.
\end{eqnarray*}
Since we are looking for the trajectories starting from the point
$(1,1,0)^t$ , we have
\begin{eqnarray*}
\ln \frac{c_1}{1+\varepsilon }=0\Rightarrow c_1=1+\varepsilon ,
\end{eqnarray*}
and so
\begin{eqnarray*}
x^3(\theta )=\ln \frac{1-\varepsilon \text{sech}\theta
}{(1-\varepsilon ) \text{sech}\theta }=\ln \frac{\cosh \theta
-\varepsilon }{1-\varepsilon }.
\end{eqnarray*}
The relation
\begin{eqnarray*}
\frac{\dot x^2}{x^2}=-\frac{\dot \mu _1}{\mu _1},
\end{eqnarray*}
leads to
\begin{eqnarray*}
x^2(\theta )=\frac{c_2(1-\varepsilon \text{sech}\theta
)}{(1-\varepsilon ^2) \text{sech}\theta },
\end{eqnarray*}
whereas from $x^2(0)=1$ we get $c_2=1+\varepsilon $. These lead to
\begin{eqnarray*}
x^2(\theta )=\frac{\cosh \theta -\varepsilon }{1-\varepsilon }.
\end{eqnarray*}
We obtain also that
\begin{eqnarray*}
\dot \mu _2=\mu _1\left( \dot x^1-x^1\frac{\partial \mathcal{H}}{\partial
\mu _2}\right) =\mu _1\dot x^1+x^1\dot \mu _1,
\end{eqnarray*}
and, consequently, $\mu _2=\mu _1x^1+c_3$. Further,
\begin{eqnarray*}
x^1(\theta )=\frac{\sinh \theta }{\sqrt{1-\varepsilon ^2}}\pm
\frac{ c_3(1-\varepsilon \text{sech}\theta )}{\left( 1-\varepsilon
^2\right) \text{ sech}\theta }.
\end{eqnarray*}
From $x^1(0)=1$ we obtain that $c_3=1+\varepsilon $ and this yields
\begin{eqnarray*}
x^1(\theta )=\frac{\sinh \theta }{\sqrt{1-\varepsilon ^2}}+\frac{\cosh
\theta -\varepsilon }{1-\varepsilon }.
\end{eqnarray*}

\begin{remark}
If $\varepsilon =0$ we regain the case of distributional systems with
quadratic cost with the solution
\begin{eqnarray*}
x^1(t)=\sinh t+\cosh t,\quad x^2(t)=\cosh t,\quad x^3(t)=\ln \cosh t.
\end{eqnarray*}
\end{remark}

\newpage\

\subsubsection{\textbf{Nonholonomic distribution}}

We assume that the distribution $D=\left\langle
X_1,X_2,...,X_m\right\rangle $ is nonholonomic and is strong
bracket generating \cite{Str1}, i.e. sections of $D$ and first
iterated brackets span the entire tangent space $ TM $. By a
well-known theorem of Chow, the system is controllable, that is
any two points are connected through a horizontal curve (M is
assumed to be connected). We also suppose that the vectors
$\mathcal{B} =\{X_1,X_2,...,X_m,[X_i,X_j]\}$ determine a base in
$TM.$ The space $E=TM$ with the base $\mathcal{B}$ is a Lie
algebroid over $M$ with at least one structural function nonzero.
The anchor $\sigma :E\rightarrow TM$ is just the identity and the
matrix corresponding to $\sigma $ is determined by the base
vectors.

The control system can be written as
\begin{eqnarray*}
\dot x=\sum_{i=1}^mu_iX_i(x)+0X_{m+1}+\dots 0X_n,
\end{eqnarray*}
with
\begin{eqnarray*}
\min_{u(\cdot )}\int_0^T\mathcal{L}(u(t))dt.
\end{eqnarray*}
To solve this minimization problem we consider the Lagrangian
\begin{eqnarray*}
\widetilde{\mathcal{L}}=\mathcal{L}+\sum_{k=m+1}^n\lambda _ku^k,
\end{eqnarray*}
($\lambda _k$ are the Lagrange multipliers) and still work via the maximum
principle but at the level of Lie algebroids. We set $\mu _i=\frac{\partial
\widetilde{\mathcal{L}}}{\partial u^i}$ and
\begin{eqnarray*}
\mathcal{H}=\sum_{i=1}^n\mu _iu^i-\widetilde{\mathcal{L}},
\end{eqnarray*}
thus obtaining $\mu _j=\frac{\partial \mathcal{L}}{\partial u^j}$,
$\mu _k=\lambda _k$ with $j=\overline{1,m}$ and
$k=\overline{m+1,n}.$ Since $ \mathcal{L}$ is 2-homogeneous with
respect to $u^i$, $i=\overline{1,m}$, we get
\begin{eqnarray*}
\mathcal{H}=\sum_{i=1}^m\frac{\partial \mathcal{L}}{\partial
u^i}u^i- \mathcal{L}=\mathcal{L,}
\end{eqnarray*}
with the constrains $\mu _k=\lambda _k$, $k=\overline{m+1,n}$.

We consider the following distributional system with positive homogeneous
cost \cite{Po14}:

\begin{eqnarray*}
\dot x=u^1X_1+u^2X_2,\quad x=\left(
\begin{array}{c}
x^1 \\
x^2 \\
x^3
\end{array}
\right) \in {R}^3,\ X_1=\left(
\begin{array}{c}
1 \\
0 \\
0
\end{array}
\right) ,\ X_2=\left(
\begin{array}{c}
0 \\
1 \\
x^1
\end{array}
\right) ,
\end{eqnarray*}
and
\begin{eqnarray*}
\min_{u(\cdot )}\int_0^T\mathcal{F}(u(t))dt,\quad
\mathcal{F}(u)=\sqrt{ (u^1)^2+(u^2)^2}+\varepsilon u^1,\quad 0\leq
\varepsilon <1.
\end{eqnarray*}
We are looking for the trajectories starting from the origin and
parameterized by arclength. We have
\[
X_1=\frac \partial {\partial x^1},\quad X_2=\frac \partial {\partial
x^2}+x^1\frac \partial {\partial x^3},
\]
\[
X_3=[X_1,X_2]=\frac \partial {\partial x^3},
\]
and $\langle X_1,X_2,X_3\rangle \equiv T{R}^3,$ hence the distribution $D=<$
$X_1,X_2>$ of constant rank is strong bracket generating.

\begin{remark}
We can work, as in classical case, directly on the cotangent bundle by
computing the Hamiltonian $H(x,p)=\mathcal{H}(x,\mu )$ , $\mu =\sigma
^{*}(p).$ Since
\begin{eqnarray*}
\mathcal{H}=\frac 12\left( \sqrt{\frac{(\mu _1)^2}{(1-\varepsilon
^2)^2}+ \frac{(\mu _2)^2}{1-\varepsilon ^2}}-\frac{\varepsilon \mu
_1}{1-\varepsilon ^2}\right) ^2,\quad \mu _3=\lambda ,
\end{eqnarray*}
\begin{eqnarray*}
\left(
\begin{array}{c}
\mu _1 \\
\mu _2
\end{array}
\right) =\left(
\begin{array}{ccc}
1 & 0 & 0 \\
0 & 1 & x^1
\end{array}
\right) \left(
\begin{array}{c}
p_1 \\
p_2 \\
p_3
\end{array}
\right) ,
\end{eqnarray*}
we obtain
\begin{eqnarray*}
H=\frac 12\left( \sqrt{\frac{\left( p_1\right) ^2}{(1-\varepsilon
^2)^2}+ \frac{\left( p_2+p_3x^1\right) ^2}{1-\varepsilon
^2}}-\frac{\varepsilon p_1}{ 1-\varepsilon ^2}\right) ^2.
\end{eqnarray*}
Unfortunately, from the Hamilton equations on $T^{*}M$ a very complicated
system of implicit differential equations is obtained.
\end{remark}

We will use a different approach. Let us take $M={R}^3$ and $E=TM$ with the
basis $\{X_{1,}X_2,X_3\}$. $E$ is a Lie algebroid over $M$ with at least one
structural function nonzero. The anchor $\sigma :E\rightarrow TM$ is just
the identity and the matrix of $\sigma $ with respect to the basis of $E$
and $TM$ basis is
\begin{eqnarray*}
\sigma _\alpha ^i=\left(
\begin{array}{ccc}
1 & 0 & 0 \\
0 & 1 & 0 \\
0 & x^1 & 1
\end{array}
\right) .
\end{eqnarray*}
From the structure equations of Lie algebroids and the relation $[X_\alpha
,X_\beta ]=L_{\alpha \beta }^\gamma X_\gamma $ we obtain the non-zero
structural functions $L_{12}^3=1$, $L_{21}^3=-1$. Now, using the Hamilton
equations on Lie algebroids (2.1.3) we get the following systems of
differential equations
\begin{equation}
\dot x^1=\frac{\partial \mathcal{H}}{\partial \mu _1},\quad \dot
x^2=\frac{
\partial \mathcal{H}}{\partial \mu _2},\quad \dot x^3=x^1\frac{\partial
\mathcal{H}}{\partial \mu _2},  \tag{2.2.14}
\end{equation}
and
\begin{equation}
\left\{
\begin{array}{l}
\dot \mu _1=-L_{12}^3\mu _3\frac{\partial \mathcal{H}}{\partial
\mu _2}
=-\lambda \frac{\partial \mathcal{H}}{\partial \mu _2}, \\
\\
\dot \mu _2=-L_{21}^3\mu _3\frac{\partial \mathcal{H}}{\partial
\mu _1}
=\lambda \frac{\partial \mathcal{H}}{\partial \mu _1}, \\
\\
\dot \mu _3=0\Rightarrow \lambda =ct,
\end{array}
\right.  \tag{2.2.15}
\end{equation}
where
\begin{eqnarray*}
&&\ \ \ \ \ \ \ \ \ \ \ \left. \frac{\partial
\mathcal{H}}{\partial \mu _1}= \frac{\left( 1+\varepsilon
^2\right) \mu _1}{(1-\varepsilon ^2)^2}-\frac{ \varepsilon
\sqrt{\frac{(\mu _1)^2}{(1-\varepsilon ^2)^2}+\frac{(\mu _2)^2}{
1-\varepsilon ^2}}}{1-\varepsilon ^2}-\frac{\varepsilon \mu _1^2}{
(1-\varepsilon ^2)^3\sqrt{\frac{(\mu _1)^2}{(1-\varepsilon
^2)^2}+\frac{(\mu
_2)^2}{1-\varepsilon ^2}}},\right. \\
&&\ \ \ \ \ \ \ \ \ \ \ \left. \frac{\partial
\mathcal{H}}{\partial \mu _2}= \frac{\mu _2}{1-\varepsilon
^2}-\frac{\varepsilon \mu _1\mu _2}{ (1-\varepsilon
^2)^2\sqrt{\frac{(\mu _1)^2}{(1-\varepsilon ^2)^2}+\frac{(\mu
_2)^2}{1-\varepsilon ^2}}}.\right.
\end{eqnarray*}
We may use to the following transformations
\begin{equation}
\mu _1(t)=(1-\varepsilon ^2)r(t)(a\cos A\theta (t)-b\sin A\theta
(t)), \tag{2.2.16}
\end{equation}
\begin{equation*}
 \mu _2(t)=\sqrt{1-\varepsilon ^2}r(t)(a\sin
A\theta (t)+b\cos A\theta (t)), \tag{2.2.16'}
\end{equation*}
with $a^2+b^2=1$. We also have $\sqrt{\frac{(\mu
_1)^2}{(1-\varepsilon ^2)^2} +\frac{(\mu _2)^2}{(1-\varepsilon
^2)}}=\left| r\right| $. Further, (2.2.15) yields
\begin{equation}
\begin{array}{c}
c_1\left( \frac{\dot r}r(a\sin A\theta +b\cos A\theta )+A\dot \theta (a\cos
A\theta -b\sin A\theta )\right) = \\
\\
=(1+\varepsilon ^2)(a\cos A\theta -b\sin A\theta )-\varepsilon (1+(a\cos
A\theta -b\sin A\theta )^2)
\end{array}
\tag{2.2.17}
\end{equation}
and
\begin{equation}
\begin{array}{c}
c_1\left( \frac{\dot r}r(a\cos A\theta -b\sin A\theta )-A\dot \theta (a\sin
A\theta +b\cos A\theta )\right) = \\
\\
=-(a\sin A\theta +b\cos A\theta )+\varepsilon (a\cos A\theta -b\sin A\theta
)(a\sin A\theta +b\cos A\theta ),
\end{array}
\tag{2.2.18}
\end{equation}
where $c_1=\frac{(1-\varepsilon ^2)\sqrt{1-\varepsilon ^2}}\lambda $.
Reducing $\dot \theta $ and $\frac{\dot r}r$ from (2.2.17) and (2.2.18), we
get
\begin{equation}
c_1\frac{\dot r}r=\varepsilon \left( a\sin A\theta +b\cos A\theta
)(\varepsilon a\cos A\theta -\varepsilon b\sin A\theta -1\right)
\tag{2.2.19}
\end{equation}
and
\begin{equation}
c_1A\dot \theta =\left( 1-\varepsilon (a\cos A\theta -b\sin A\theta )\right)
^2.  \tag{2.2.20}
\end{equation}
These lead to
\begin{eqnarray*}
t=\frac{(1-\varepsilon ^2)\sqrt{1-\varepsilon ^2}A}\lambda \int
\frac{ d\theta }{\left( 1-\varepsilon \left( a\cos A\theta -b\sin
A\theta \right) \right) ^2}.
\end{eqnarray*}
From (2.2.19) and (2.2.20) we obtain
\begin{eqnarray*}
\frac{dr}r=\frac{A\varepsilon (a\sin A\theta +b\cos A\theta )}{\varepsilon
(a\cos A\theta -b\sin A\theta )-1}d\theta
\end{eqnarray*}
and
\begin{eqnarray*}
r=\frac 1{c\left( 1-\varepsilon \left( a\cos A\theta -b\sin A\theta \right)
\right) }.
\end{eqnarray*}
Since the geodesics are parameterized by arclength this
corresponds exactly to the $1/2$ level of the Hamiltonian and we
have
\begin{eqnarray*}
\mathcal{H}=\frac{r^2}2\left( 1-\varepsilon \left( a\cos A\theta -b\sin
A\theta \right) \right) ^2=\frac 1{2c^2}.
\end{eqnarray*}
So, $c=\pm 1$ and
\begin{eqnarray*}
r=\pm \frac 1{1-\varepsilon \left( a\cos A\theta -b\sin A\theta \right) }.
\end{eqnarray*}
From (2.2.16) we obtain
\begin{eqnarray*}
\mu _1(t) &=&\pm \frac{(1-\varepsilon ^2)(a\cos A\theta -b\sin
A\theta )}{
1-\varepsilon \left( a\cos A\theta -b\sin A\theta \right) } \\
&& \\
\mu _2(t) &=&\frac{\sqrt{1-\varepsilon ^2}(a\sin A\theta +b\cos
A\theta )}{ 1-\varepsilon \left( a\cos A\theta -b\sin A\theta
\right) }.
\end{eqnarray*}
Since $\dot \mu _2=\lambda \dot x^1$, we also have $x^1(\theta
)=\frac{\mu _2 }\lambda -a_1$. As we are looking for geodesics
with start from the origin, we have $a_1=\frac{\sqrt{1-\varepsilon
^2}b}{\lambda (1-\varepsilon a)}$ and therefore
\begin{eqnarray*}
x^1(\theta )=\frac{\sqrt{1-\varepsilon ^2}(a\sin A\theta +b\cos
A\theta )}{ \lambda \left( 1-\varepsilon \left( a\cos A\theta
-b\sin A\theta \right) \right) }-\frac{\sqrt{1-\varepsilon
^2}b}{\lambda (1-\varepsilon a)}.
\end{eqnarray*}
From $\dot \mu _1=-\lambda \dot x^2$ we get
\begin{eqnarray*}
x^2(\theta )=\frac{(1-\varepsilon ^2)(b\sin A\theta -a\cos A\theta
)}{ \lambda \left( 1-\varepsilon \left( a\cos A\theta -b\sin
A\theta \right) \right) }+\frac{(1-\varepsilon ^2)a}{\lambda
(1-\varepsilon a)}.
\end{eqnarray*}
Finally, because $\dot x^3=x^1\dot x^2$ a straightforward computation leads
to
\begin{eqnarray*}
x^3(\theta ) &=&\frac{(1-\varepsilon ^2)\sqrt{1-\varepsilon
^2}A}{\lambda ^2} \int \frac{(a\sin A\theta +b\cos A\theta
)^2}{\left( 1-\varepsilon \left(
a\cos A\theta -b\sin A\theta \right) \right) ^3}d\theta - \\
&&\ \ \ -\frac{(1-\varepsilon ^2)Ab}{\lambda ^2}\int \frac{a\sin A\theta
+b\cos A\theta }{\left( 1-\varepsilon \left( a\cos A\theta -b\sin A\theta
\right) \right) ^2}d\theta
\end{eqnarray*}

\begin{remark}
For $\varepsilon =0$ we obtain the distributional systems with
quadratic cost with the solution $\cdot $
\[
x^1(t)=\frac{a\sin \lambda t-b(1-\cos \lambda \theta )}\lambda ,\quad
a^2+b^2=1.
\]
\[
x^2(t)=\frac{b\sin \lambda t+a(1-\cos \lambda t)}\lambda ,
\]
\[
x^3(t)=\frac t{2\lambda }+\frac{b^2-a^2}{4\lambda ^2}\sin 2\lambda
t-\frac{ab }{\lambda ^2}\cos ^2\lambda t+\frac{ab}{\lambda ^2}\cos
\lambda t-\frac{b^2}{ a^2}\sin \lambda t.
\]
\end{remark}

\newpage\

\subsubsection{\textbf{Distributional systems with no constant rank of
distribution}}

Let us consider in the three dimensional space $R^3$ the drift less control
affine system \cite{Po18}
\begin{equation}
\dot x(t)=u^1X_1+u^2X_2+u^3X_3,  \tag{2.2.21}
\end{equation}
with
\[
X_1=\left(
\begin{array}{c}
1 \\
0 \\
0
\end{array}
\right) ,\quad X_2=\left(
\begin{array}{c}
0 \\
x \\
0
\end{array}
\right) ,\quad X_3=\left(
\begin{array}{c}
0 \\
0 \\
x
\end{array}
\right) ,
\]
and minimizing the cost
\begin{equation}
\min_{u(.)}\int F(u(t))dt,  \tag{2.2.22}
\end{equation}
where $F=\sqrt{\left( u^1\right) ^2+\left( u^2\right) ^2+\left(
u^3\right) ^2 }+\varepsilon u^1$, $0\leq \varepsilon <1$ is the
positive homogeneous cost (Randers metric). The distribution $D$
is generated by the vectors $X_1$, $ X_2$, $X_3$ and we can write
$D=\{X_1,X_2,X_3\}$. We observe that
\[
rankD=\left\{
\begin{array}{c}
3\quad if\quad x\not =0 \\
1\quad if\quad x=0
\end{array}
\right.
\]
In the canonical base $\left( \frac \partial {\partial x},\frac
\partial {\partial y},\frac \partial {\partial z}\right) $ of
$R^3$ we have $ X_1=\frac \partial {\partial x}$, $X_2=x\frac
\partial {\partial y}$, $ X_3=x\frac \partial {\partial z}$ and
the Lie brackets are given by
\[
\left[ X_1,X_2\right] =\frac \partial {\partial y}=X_4\notin D,\quad \left[
X_1,X_3\right] =\frac \partial {\partial z}=X_5\notin D,\quad \left[
X_2,X_3\right] =0.
\]
It results that the distribution is nonholonomic, but is bracket generating,
because the vector fields $\{X_1,X_2,X_3,X_4=\left[ X_1,X_2\right]
,X_5=\left[ X_1,X_3\right] \}$ generate the entire space $R^3$.

From (2.2.21) we obtain
\[
\frac{dx}{dt}=u^1=s^1,\quad \frac{dy}{dt}=u^2x=s^2,\quad
\frac{dz}{dt} =u^3x=s^3.
\]
The cost function can be written in the form ($x\not =0$)
\begin{eqnarray*}
F &=&\sqrt{\left( u^1\right) ^2+\left( u^2\right) ^2+\left(
u^3\right) ^2} +\varepsilon u^1=\sqrt{\left( s^1\right)
^2+\frac{\left( s^2\right) ^2}{x^2}+
\frac{\left( s^3\right) ^2}{x^2}}+\varepsilon s^1 \\
\ &=&\sqrt{g_{ij}s^is^j}+\sum_{i=1}^3b^is^i
\end{eqnarray*}
(Einstein's summation) where $b^1=\varepsilon $, $b^2=0$, $b^3=0$ and
\[
g_{ij}=\left(
\begin{array}{ccc}
1 & 0 & 0 \\
0 & 1/x^2 & 0 \\
0 & 0 & 1/x^2
\end{array}
\right) .
\]
The Lagrangian has the form $L=\frac 12F^2$ and using \cite{Mi5} (Th. 4.5
pp. 191) we obtain the Hamiltonian in the form
\begin{equation}
H=\frac 12\left( \sqrt{\widetilde{g}_{ij}p_ip_j}-\widetilde{b}^ip_i\right) ,
\tag{2.2.23}
\end{equation}
where
\[
\widetilde{g}_{ij}=\frac 1{1-b^2}g^{ij}+\frac 1{(1-b^2)^2}b^ib^j,\quad
\widetilde{b}^i=\frac 1{1-b^2}b^i,\quad b=\sqrt{g_{ij}b^ib^j},
\]
and $g^{ij}$ is the inverse of the matrix $g_{ij}$. In these conditions we
obtain that
\[
b^2=\varepsilon ^2,\quad \widetilde{b}^1=\frac \varepsilon {1-\varepsilon
^2},\quad \widetilde{b}^2=0,\quad \widetilde{b}^3=0,
\]
\[
g^{ij}=\left(
\begin{array}{ccc}
1 & 0 & 0 \\
0 & x^2 & 0 \\
0 & 0 & x^2
\end{array}
\right)
\]
and and it results
\[
\widetilde{g}^{ij}=\left(
\begin{array}{ccc}
\frac 1{\left( 1-\varepsilon ^2\right) ^2} & 0 & 0 \\
0 & \frac{x^2}{1-\varepsilon ^2} & 0 \\
0 & 0 & \frac{x^2}{1-\varepsilon ^2}
\end{array}
\right)
\]
From (2.2.23) we obtain
\begin{equation}
H=\frac 12\left( \sqrt{\frac{p_1^2}{(1-\varepsilon
^2)^2}+\frac{p_2^2x^2}{ 1-\varepsilon
^2}+\frac{p_3^2x^2}{1-\varepsilon ^2}}-\frac{\varepsilon p_1}{
1-\varepsilon ^2}\right) ^2  \tag{2.2.24}
\end{equation}
or, in the equivalent form
\[
H=\frac{\left( 1+\varepsilon ^2\right) p_1^2}{2\left(
1-\varepsilon ^2\right) ^2}+\frac{\left( p_2^2+p_3^2\right)
x^2}{2\left( 1-\varepsilon ^2\right) }-\frac{\varepsilon
p_1}{1-\varepsilon ^2}\sqrt{\frac{p_1^2}{ (1-\varepsilon
^2)^2}+\frac{\left( p_2^2+p_3^2\right) x^2}{1-\varepsilon ^2}} .
\]
In the case $x=0$ we obtain
\[
L=\frac 12F^2=\frac{(1+\varepsilon )^2u_1^2}2,
\]
with the constraints $\stackrel{\cdot }{y}=0$, $\stackrel{\cdot }{z}=0.$
Using Lagrange multipliers we obtain
\[
L_1=L+\lambda _1\stackrel{\cdot }{y}+\lambda _2\stackrel{\cdot }{z},
\]
and from Legendre transformation it results
\[
H_1=\frac 12\frac{p_1^2}{(1+\varepsilon )^2}.
\]
For $x=0$ from (2.2.24) we have
\[
H=\frac 12\left( \frac{p_1}{1-\varepsilon ^2}-\frac{\varepsilon
p_1}{ 1-\varepsilon ^2}\right) ^2=\frac 12\frac{p_1^2}{\left(
1+\varepsilon \right) ^2},
\]
which leads to the equality
\[
\left. H\right| _{x=0}=H_1.
\]
Next, if we denote
\[
\Theta =\frac{p_1^2}{(1-\varepsilon ^2)^2}+\frac{\left( p_2^2+p_3^2\right)
x^2}{1-\varepsilon ^2},
\]
then the Hamilton's equations (2.1.1) lead to the following differential
equations
\begin{equation}
\frac{dx}{dt}=\frac{\partial H}{\partial p_1}=\frac{(1+\varepsilon
^2)p_1}{ (1-\varepsilon ^2)^2}-\frac \varepsilon {1-\varepsilon
^2}\sqrt{\Theta }- \frac{\varepsilon p_1^2}{\left( 1-\varepsilon
^2\right) ^3}\frac 1{\sqrt{ \Theta }},  \tag{2.2.25}
\end{equation}
\begin{equation}
\frac{dy}{dt}=\frac{\partial H}{\partial
p_2}=\frac{p_2x^2}{1-\varepsilon ^2} -\frac{\varepsilon
p_1p_2x^2}{(1-\varepsilon ^2)^2}\frac 1{\sqrt{\Theta }},
\tag{2.2.26}
\end{equation}
\begin{equation}
\frac{dz}{dt}=\frac{\partial H}{\partial
p_3}=\frac{p_3x^2}{1-\varepsilon ^2} -\frac{\varepsilon
p_1p_3x^2}{(1-\varepsilon ^2)^2}\frac 1{\sqrt{\Theta }},
\tag{2.2.27}
\end{equation}
\begin{equation}
\frac{dp_1}{dt}=-\frac{\partial H}{\partial x}=-\frac{\left(
p_2^2+p_3^2\right) x}{1-\varepsilon ^2}+\frac{\varepsilon p_1\left(
p_2^2+p_3^2\right) x}{\left( 1-\varepsilon ^2\right) ^2}\frac 1{\sqrt{\Theta
}},  \tag{2.2.28}
\end{equation}
\[
\frac{dp_2}{dt}=-\frac{\partial H}{\partial y}=0\Rightarrow p_2=a=const.
\]
\[
\frac{dp_3}{dt}=-\frac{\partial H}{\partial z}=0\Rightarrow \Rightarrow
p_3=b=const.
\]
In these conditions the relation $\Theta
=\frac{p_1^2}{(1-\varepsilon ^2)^2}+ \frac{\left(
p_2^2+p_3^2\right) x^2}{1-\varepsilon ^2}$ leads to the following
change of variables:
\[
x(t)=\frac{\sqrt{1-\varepsilon ^2}r(t)\sin A\theta
(t)}{\sqrt{a^2+b^2}} ,\quad p_1(t)=\left( 1-\varepsilon ^2\right)
r(t)\cos A\theta (t).
\]
It results $\Theta =r^2(t)$ and from (2.2.25) we get
\[
\frac{dx}{dt}=\frac{(1+\varepsilon ^2)r\cos A\theta
}{1-\varepsilon ^2}- \frac{\varepsilon r}{1-\varepsilon
^2}-\frac{\varepsilon r\cos ^2A\theta }{ 1-\varepsilon ^2}.
\]
But
\[
\frac{dx}{dt}=\frac{\sqrt{1-\varepsilon ^2}}{\sqrt{a^2+b^2}}\left(
\stackrel{ \cdot }{r}\sin A\theta +rA\stackrel{.}{\theta }\cos
A\theta \right) ,
\]
and it results
\begin{equation}
c_1\left( \stackrel{\cdot }{r}\sin A\theta +rA\stackrel{\cdot }{\theta }\cos
A\theta \right) =(1+\varepsilon ^2)r\cos A\theta -\varepsilon r(1+\cos
^2A\theta ),  \tag{2.2.29}
\end{equation}
where we have denoted
\[
c_1=\frac{\left( 1-\varepsilon ^2\right) \sqrt{1-\varepsilon
^2}}{\sqrt{ a^2+b^2}}.
\]
The equation (2.2.28) yields
\[
\frac{dp_1}{dt}=\frac{\sqrt{a^2+b^2}}{\sqrt{1-\varepsilon ^2}}\left( -r\sin
A\theta +\varepsilon r\cos A\theta \sin A\theta \right) .
\]
But
\[
\frac{dp_1}{dt}=\left( 1-\varepsilon ^2\right) \left( \stackrel{.}{r}\cos
A\theta -rA\stackrel{.}{\theta }\sin A\theta \right) ,
\]
which leads to
\begin{equation}
c_1\left( \stackrel{\cdot }{r}\cos A\theta -rA\stackrel{\cdot }{\theta }\sin
A\theta \right) =-r\sin A\theta +\varepsilon r\cos A\theta \sin A\theta .
\tag{2.2.30}
\end{equation}
The equation (2.2.29) multiplied by $\cos A\theta $, minus equation (2.2.30)
multiplied by $\sin A\theta $ lead to the equation
\begin{equation}
c_1A\frac{d\theta }{dt}=\left( 1-\varepsilon \cos A\theta \right) ^2,
\tag{2.2.31}
\end{equation}
and it results
\[
t=\frac{A\left( 1-\varepsilon ^2\right) \sqrt{1-\varepsilon
^2}}{\sqrt{ a^2+b^2}}\int \frac 1{\left( 1-\varepsilon \cos
A\theta \right) ^2}d\theta
\]
Moreover, the equation (2.2.29) multiplied by $\sin A\theta $, plus equation
(2.2.30) multiplied by $\cos A\theta $ lead to the equation
\begin{equation}
c_1\frac{dr}{dt}=\varepsilon r\sin A\theta \left( \varepsilon \cos A\theta
-1\right) .  \tag{2.2.32}
\end{equation}
From the equations (2.2.31) and (2.2.32) it results
\[
\frac{dr}r=-A\frac{\varepsilon \sin A\theta }{1-\varepsilon \cos
A\theta } d\theta ,
\]
which leads to the following result
\[
\ln r=-\ln c(1-\varepsilon \cos A\theta ),\quad c\in R,
\]
and we get
\begin{equation}
r(t)=\frac 1{c(1-\varepsilon \cos A\theta (t))}.  \tag{2.2.33}
\end{equation}
Using (2.2.33) the Hamiltonian become
\begin{eqnarray*}
H &=&\frac 12(1+\varepsilon ^2)r^2\cos ^2A\theta +\frac 12r^2\sin ^2A\theta
-\varepsilon r^2\cos A\theta \\
\ &=&\frac 12r^2+\frac 12\varepsilon ^2r^2\cos ^2A\theta -\varepsilon
r^2\cos A\theta \\
\ &=&\frac{r^2}2\left( 1-\varepsilon \cos A\theta \right) \\
\ &=&\frac 1{2c^2}
\end{eqnarray*}
Considering the integral curves parameterized by arclength, that
corresponds to fix the level $\frac 12$ of the Hamiltonian, we
have $c=\pm 1$ and it results
\begin{equation}
r=\pm \frac 1{1-\varepsilon \cos A\theta }.  \tag{2.2.34}
\end{equation}
In these conditions, from (2.2.34) we obtain
\begin{equation}
x(t)=\pm \frac{\sqrt{1-\varepsilon ^2}}{\sqrt{a^2+b^2}}\frac{\sin
A\theta (t) }{1-\varepsilon \cos A\theta (t)}.  \tag{2.2.35}
\end{equation}
The differential equation (2.2.26) yields
\begin{equation}
\frac{dy}{dt}=\frac{ar^2\sin ^2A\theta }{a^2+b^2}\left( 1-\varepsilon \cos
A\theta \right) =\frac a{a^2+b^2}\frac{\sin ^2A\theta }{1-\varepsilon \cos
A\theta }.  \tag{2.2.36}
\end{equation}
From (2.2.31) and (2.2.36) it results
\begin{equation}
y(t)=\frac{aA\left( 1-\varepsilon ^2\right) ^{3/2}}{\left( a^2+b^2\right)
^{3/2}}\int \frac{\sin ^2A\theta }{\left( 1-\varepsilon \cos A\theta \right)
^3}d\theta (t).  \tag{2.2.37}
\end{equation}
In the same way we obtain
\begin{equation}
z(t)=\frac{bA\left( 1-\varepsilon ^2\right) ^{3/2}}{\left( a^2+b^2\right)
^{3/2}}\int \frac{\sin ^2A\theta }{\left( 1-\varepsilon \cos A\theta \right)
^3}d\theta (t).  \tag{2.2.38}
\end{equation}
In the particular case of $\varepsilon =0$ we obtain a distributional system
with quadratic cost
\[
F=\sqrt{\left( u^1\right) ^2+\left( u^2\right) ^2+\left( u^3\right) ^2},
\]
with the solution
\[
x(t)=\pm \frac{\sin \alpha t}\alpha ,
\]
\[
y(t)=\frac{at}{2\alpha ^2}-\frac{a\sin 2\alpha t}{4\alpha ^{3/2}},
\]
\[
z(t)=\frac{bt}{2\alpha ^2}-\frac{a\sin 2\alpha t}{4\alpha ^{3/2}}.
\]
where $\alpha =\sqrt{a^2+b^2}$, which are the geodesics in the framework of
the so called sub-Riemannian geometry \cite{Bel}.

\newpage\

\subsubsection{\textbf{Distributional systems with degenerate cost}}

Let us consider in the two dimensional space $R^2$ the drift less control
affine system \cite{Po22}
\begin{equation}
\dot x(t)=u^1X_1+u^2X_2,  \tag{2.2.39}
\end{equation}
with
\[
X_1=\left(
\begin{array}{c}
1 \\
0
\end{array}
\right) ,\quad X_2=\left(
\begin{array}{c}
0 \\
x
\end{array}
\right) ,
\]
and minimizing the cost
\begin{equation}
\min_{u(.)}\int_IF(u(t))dt,  \tag{2.2.40}
\end{equation}
where $F=\sqrt{(u^1)^2+(u^2)^2}+u^1$ is the positive homogeneous cost
(Kropina metric).

I have to remark that in the case $u^1\leq 0$, $u^2=0$ it results $F=0$ that
is we obtain a degenerate cost (metric). The distribution $D$ is generated
by the vectors $X_1,X_2$ and we have
\[
rankD=\left\{
\begin{array}{c}
2\quad if\quad x\not =0, \\
1\quad if\quad x=0.
\end{array}
\right.
\]
In the canonical base of $R^2$ given by $\left\{ \frac \partial {\partial
x},\frac \partial {\partial y}\right\} $ we can write
\[
X_1=\frac \partial {\partial x},\quad X_2=x\frac \partial {\partial y}
\]
and the Lie brackets are given by
\[
\left[ X_1,X_2\right] =\frac \partial {\partial y}=X_3\notin D.
\]
It results that the distribution is nonholonomic, but is bracket generating,
because the vectors fields $\left\{ X_1,X_2,X_3=\left[ X_1,X_2\right]
\right\} $ generate the entire space $R^2$. From (2.2.39) we obtain
\begin{eqnarray*}
\frac{dx}{dt} &=&u^1=s^1. \\
\frac{dy}{dt} &=&u^2x=s^2.
\end{eqnarray*}
The cost function can be written in the following form ($x\not =0$)
\[
F=\sqrt{(u^1)^2+(u^2)^2}+u^1=\sqrt{(s^1)^2+\frac{(s^2)^2}{x^2}}+s^1=\sqrt{
g_{ij}s^is^j}+\sum_{i=1}^3b^is^i
\]
(Einstein's summation, $i,j=1,2$) where $b^1=1$, $b^2=0$ and
\[
g_{ij}=\left(
\begin{array}{cc}
1 & 0 \\
1 & 1/x^2
\end{array}
\right) .
\]
The Lagrangian function has the form $L=\frac 12F^2$ and using \cite{Mi5} we
obtain the Hamiltonian in the form
\begin{equation}
H=\frac 12\left( \frac{g^{ij}p_ip_j}{2b^ip_i}\right) ^2,  \tag{2.2.41}
\end{equation}
where
\[
g^{ij}=\left(
\begin{array}{cc}
1 & 0 \\
1 & x^2
\end{array}
\right) ,
\]
is the inverse of the matrix $g_{ij}$. In these conditions we obtain that
\[
H=\frac 12\left( \frac{p_1^2+x^2p_2^2}{2p_1}\right) ^2,
\]
or, in the equivalent form
\[
H=\frac{p_1^2}8+\frac{p_2^2}{8p_1^2}x^4+\frac{p_2^2}4x^2.
\]
The Hamilton's equations (2.1.1) lead to the following differential
equations
\begin{equation}
\frac{dx}{dt}=\frac{\partial H}{dp_1}=\frac{p_1}4-\frac{p_2^4x^4}{4p_1^3}
\tag{2.2.42}
\end{equation}
\begin{equation}
\frac{dy}{dt}=\frac{\partial H}{dp_2}=\frac{p_2^3x^4}{2p_1^2}+\frac{p_2x^2}2
\tag{2.2.43}
\end{equation}
\begin{equation}
\frac{dp_1}{dt}=-\frac{\partial
H}{dx}=-\frac{p_2^4x^3}{2p_1^2}-\frac{p_2^2x} 2  \tag{2.2.44}
\end{equation}
\[
\frac{dp_2}{dt}=-\frac{\partial H}{dy}=0\Rightarrow p_2=a=const.
\]
and it results
\[
\frac{dx}{dt}=\frac{p_1^4-a^4x^4}{4p_1^3},
\]
\[
\frac{dy}{dt}=\frac{a^3x^4+ap_1^2x^2}{2p_1^2},
\]
\[
\frac{dp_1}{dt}=-\frac{a^2x}2-\frac{a^4x^3}{2p_1^2}.
\]
In these conditions the expression of the Hamiltonian leads to the following
change of variables:
\begin{equation}
x(t)=\frac{r(t)\sin A\theta (t)}a,\quad p_1(t)=r(t)\cos A\theta (t)
\tag{2.2.45}
\end{equation}
We obtain the equations
\[
\frac{dx}{dt}=\frac{r\left( \cos ^4A\theta -\sin ^4A\theta \right) }{4\cos
^3A\theta },
\]
and
\[
\frac{dp_1}{dt}=-\frac{ar\sin A\theta }2-\frac{ar\sin ^3A\theta }{2\cos
^2A\theta }.
\]
But on the other hand
\[
\frac{dx}{dt}=\frac{\stackrel{\cdot }{r}}a\sin A\theta
+\frac{rA\stackrel{ \cdot }{\theta }}a\cos A\theta
\]
\[
\frac{dp_1}{dt}=\stackrel{\cdot }{r}\cos A\theta -rA\stackrel{\cdot }{\theta
}\sin A\theta
\]
and it follows
\begin{equation}
\frac 1a\left( \frac{\stackrel{\cdot }{r}}r\sin A\theta
+A\stackrel{\cdot }{ \theta }\cos A\theta \right) =\frac{\cos
^4A\theta -\sin ^4A\theta }{4\cos ^3A\theta },  \tag{2.2.46}
\end{equation}
\begin{equation}
\frac 1a\left( \frac{\stackrel{\cdot }{r}}r\cos A\theta
-A\stackrel{\cdot }{ \theta }\sin A\theta \right) =-\frac{\sin
A\theta }2-\frac{\sin ^3A\theta }{ 2\cos ^2A\theta }.
\tag{2.2.47}
\end{equation}
The equation (2.2.46) multiplied by $\cos A\theta ,$ minus equation (2.2.47)
multiplied by $\sin A\theta $ leads to the equation
\begin{equation}
\frac Aa\frac{d\theta }{dt}=\frac 1{4\cos ^3A\theta },  \tag{2.2.48}
\end{equation}
which yields
\[
t=\frac{4A}a\int \cos ^2A\theta d\theta ,
\]
and it follows
\[
t=\frac 1a\left( \sin 2A\theta +2A\theta \right) .
\]
Moreover, the equation (2.2.46) multiplied by $\sin A\theta ,$ plus equation
(2.2.47) multiplied by $\cos A\theta $ leads to the equation
\begin{equation}
\frac 1a\frac{\stackrel{\cdot }{r}}r=-\frac{\sin A\theta }{4\cos ^3A\theta }.
\tag{2.2.49}
\end{equation}
Using the equations (2.2.48) and (2.2.49) we obtain
\[
\frac{dr}r=-A\frac{\sin A\theta }{\cos A\theta }d\theta ,
\]
which leads to the following result
\[
\ln \left| r\right| +\ln c_1=\ln \left| \cos A\theta \right| ,\quad c_1\in
R^{+},
\]
and we get
\begin{equation}
r(t)=\frac 1{c_1}\cos A\theta (t).  \tag{2.2.50}
\end{equation}
Using the change of variables (2.2.45) the Hamiltonian become
\[
H=\frac{r^2\cos ^2A\theta }8+\frac{r^2\sin ^4A\theta }{8\cos
^2A\theta }+ \frac{r^2\sin ^2A\theta }4=\frac{r^2}{8\cos ^2A\theta
},
\]
and from (2.2.50) we get
\[
H=\frac 1{8c_1^2}.
\]
Considering the integral curves parameterized by arclength, that
corresponds to fix the level $1/2$ of the Hamiltonian, we have
$c_1=\pm 1/2$ and it results
\begin{equation}
r=\pm 2\cos A\theta .  \tag{2.2.51}
\end{equation}
which together with (2.2.45) lead to the result
\begin{equation}
x(t)=\pm \frac{\sin 2A\theta (t)}a  \tag{2.2.52}
\end{equation}
The differential equation (2.2.43) yields by direct computation to
\begin{equation}
\frac{dy}{dt}=\frac{2\sin ^2A\theta }a,  \tag{2.2.53}
\end{equation}
and from (2.2.48) and (2.2.53) it follows
\[
dy=\frac{2A}{a^2}\sin ^2(2A\theta )d\theta ,
\]
which yields
\begin{equation}
y(t)=\frac{A\theta (t)}{a^2}-\frac{\sin 4A\theta (t)}{4a^2}.  \tag{2.2.54}
\end{equation}
Finally, the solution is
\begin{eqnarray*}
x(t) &=&\pm \frac{\sin 2A\theta (t)}a, \\
y(t) &=&\frac{A\theta (t)}{a^2}-\frac{\sin 4A\theta (t)}{4a^2}.
\end{eqnarray*}

\newpage\

\subsection{\textbf{Sub-Riemannian geometry }}

\ If $M$ is a smooth $n$-dimensional manifold then a sub-Riemannian
structure on $M$ is a pair $(D,g),$ where $D$ is a distribution of rank $m$
and $g$ is a Riemannian metric on $D$. A sub-Riemannian manifold $(M,D,g)$
is a smooth manifold $M$ equipped with a sub-Riemannian structure \cite
{Str1, Mo, Bel}.\\A piecewise smooth curve $c:I\subset R\rightarrow M$ is
called \textit{\ horizontal} if its tangent vectors are in $D$, i.e. $\dot
c(t)\in D_{c(t)}\subset TM,$ for almost every $t\in I.$ In sub-Riemannian
geometry the \textit{length} of a horizontal curve $c$ is defined by
\begin{equation}
L(c)=\int_I\sqrt{g(\dot c(t))}dt,  \tag{2.3.1}
\end{equation}
where $g$ is a Riemannian metric on $D$. The \textit{distance}
from $a$ to $ b $ is
\[
d(a,b)=\inf (L(c)),
\]
where the infimum is taken over all horizontal curves connecting $a$ to $b$.
The distance is assumed to be infinite if there is no horizontal curve that
connects these two points.

If locally the distribution $D$ of rank $m$ is generated by $X_i,$
$i= \overline{1,m}$, a sub-Riemannian structure on $M$ is locally
given by a control system
\begin{equation}
\dot x=\sum_{i=1}^mu_i(t)X_i(x),  \tag{2.3.2}
\end{equation}
of constant rank $m$, with the controls $u(.).$ The
\textit{controlled paths} are obtained by integrating the above
system. If $D$ is assumed to be \textit{bracket generating,} i.e.
sections of $D$ and iterated brackets span the entire tangent
space $TM$, by a well-known theorem of Chow \cite{Ch} the system
(2.3.2) is \textit{controllable}, that is for any two points $a$
and $ b$ there exists a horizontal curve which connects these
points ($M$\ is assumed to be connected).

The concept of the sub-Riemannian geometry can be extended to a more general
setting, by replacing the Riemannian metric with a Finslerian one. For the
theory of optimal control this extension is equivalent to the change of the
quadratic cost of a control affine system with a positive homogeneous cost.
Also, the case when the rank of $D$ is not constant may produce interesting
examples. We do not intend to develop a comprehensive study of the
sub-Riemannian geometry. In the present section we introduce two particular
sub-Finslerian geometries \cite{Hr3}: the Grushin plane and the Heisenberg
group, endowed with some special Randers metrics \cite{Bao}. The geodesics
of these geometries are obtained by using two different approaches: a direct
application of the Pontryagin Maximum Principle for the Grushin plane and
the same principle but combined with some results on Lie algebroids for the
Heisenberg group.

\begin{definition}
A sub-Finslerian structure on $M$ is a triple $(E,\sigma
,\mathcal{F})$ where
\\1) ($E,\pi ,M)$ is a vector bundle over $M$, with the projection map $\pi
. $\\2) $\sigma :E\rightarrow TM$ is a morphism of vector
bundles.\\3) $ \mathcal{F}$ is a Finsler metric on $E$, i.e.
$\mathcal{F}:E\rightarrow [0,\infty )$ and satisfies the following
properties:\\a) $\mathcal{F}$ is $ C^{\infty \text{ }}$on
$E\setminus \{0\}.$\\b) $\mathcal{F}(\lambda u)=\lambda
\mathcal{F}(u)$ for $\lambda >0$ and $u\in E_x,x\in M.$\\c) For
each $y\in E_x\backslash \{0\}$ the quadratic form
\[
g_y(u,v)=\frac 12\frac{\partial \mathcal{F}^2}{\partial s\partial
t} (y+su+tv)_{s,t=0}
\]
$u,v\in E_x,$ $x\in M$, is positive definite$.$
\end{definition}

\begin{example}
(i) $E=$ $M\times \Bbb{R}^m$, $\{X_1,...,X_m\}$ is a system of \textit{m}
vector fields on $M$, $\sigma :E\rightarrow TM$, given by
\begin{equation}
\sigma (x,u)=\sum_{i=1}^mu_i(t)X_i(x),  \tag{2.3.3}
\end{equation}
and $\mathcal{F}$ is a Minkowski norm on $\Bbb{R}^m$.\\(ii) $E=D,$ $\sigma
:D\rightarrow TM$ the inclusion and $\mathcal{F}$ a Finsler metric on $D$.
\end{example}

We can associate to any sub-Finslerian structure $(E,\sigma ,\mathcal{F})$ a
Finsler metric on $Im\sigma \subset TM$, defined as following
\begin{equation}
F(v)={\inf_u}\{\mathcal{F}(u)|\quad u\in E_x,\ \sigma (u)=v\}.  \tag{2.3.4}
\end{equation}
for each
\[
v\in (Im\sigma )_x\subset T_xM,\ x\in M.
\]

\begin{definition}
A curve $u:I\rightarrow E$ is called admissible if there is an absolute
continuous curve $c:I\rightarrow M$, called horizontal, such that $\pi
(u(t))=c(t)$ and $\sigma ((u(t))=\dot c(t)$ , $t\in I$.
\end{definition}

The length of an absolutely continuous horizontal curve $c(t)$, $t\in I$ is
\begin{equation}
length(c)=\int_I\mathcal{F}(u(t))dt=\int_IF(\dot c(t))dt.  \tag{2.3.5}
\end{equation}
We can also consider the sub-Finslerian distance
\[
d(a,b)=\inf length(c),
\]
where the infimum is taken over all horizontal curves connecting
$a$ and $b$. This distance is infinite if there is no admissible
curve joining $a$ and $ b$. However, if we assume that $Im\sigma $
is bracket generating, Chow's theorem guarantees that the
sub-Finslerian distance between points is finite.

\begin{definition}
A length minimizing geodesic or shortly a minimizer is an absolutely
continuous horizontal curve on $M$ that makes the distance between two
points.
\end{definition}

\begin{remark}
The energy of a horizontal curve is
\[
E(c)=\frac 12\int_IF^2(\dot c(t))dt,
\]
and it can easily be proved that if a curve is parameterized to a
constant speed, then it minimize the length integral if and only
if it minimize the energy integral.
\end{remark}

If we take the Lagrangians $L=\frac 12F^2$ and $\mathcal{L}=\frac
12\mathcal{ F}^2$ we have $\mathcal{L}=L\circ \sigma .$

The Fenchel-Legendre dual of $L$ is the Hamiltonian \cite{Hr3}
\[
H(p)=\sup_v\left\{ \left\langle p,v\right\rangle -L(v)\right\} =
\]
\[
\sup_v\left\{ \left\langle p,v\right\rangle +\sup_u\left\{
-\mathcal{L} (u);\sigma (u)=v\right\} \right\} =
\]
\[
{\sup_{u,v}}\left\{ \left\langle p,v\right\rangle -\mathcal{L}(u);\sigma
(u)=v\right\} =
\]
\[
{\sup_u}\left\{ \left\langle p,\sigma (u)\right\rangle
-\mathcal{L} (u)\right\} =
\]
\[
{\sup_u}\left\{ \left\langle \sigma ^{\star }(p),u\right\rangle
-\mathcal{L} (u)\right\} =\mathcal{H}(\sigma ^{\star }(p))
\]
Hence
\begin{equation}
H(p)=\mathcal{H}(\mu ),\quad \mu =\sigma ^{\star }(p)  \tag{2.3.6}
\end{equation}
$p$ $\in T_x^{*}M$, $\mu \in E_x^{*}.$ The Hamiltonian $H$ on $T^{*}M$ is
degenerate on $Ker\sigma ^{\star }$.

The Hamiltonian $H$ generates a system of differential equations which can
be written in terms of canonical coordinates $(x,p)$ in the standard form:
\begin{equation}
\stackrel{.}{x}^i=\frac{\partial H}{\partial p_i},\text{
}\stackrel{.}{p}_i=- \frac{\partial H}{\partial x^i}  \tag{2.3.7}
\end{equation}
Based on the Pontryagin Maximum Principle we can prove:

\begin{theorem}
Let $(x(t),p(t))$ be a solution of the Hamilton$\pi $ equations. Then every
sufficient short subarc of $x(t)$ is a minimizing sub-Finslerian geodesic.
This subarc is the unique minimizer joining its end points.
\end{theorem}

\begin{definition}
The projected curve $x(t)$ will be called normal geodesic or shortly
geodesic.
\end{definition}

\begin{remark}
Contrary to Finslerian geometry, not every minimizer is the projection of a
solution of (2.3.7) i.e. is a normal geodesic. Those minimizers that are not
normal geodesics are called singular geodesics \cite{Str1, Mo}.
\end{remark}

\newpage\

\subsubsection{\textbf{Grushin case}}

We consider the following distributional system with positive homogeneous
cost (Grushin Plane) \cite{Hr3}
\[
\dot x=u^1X_1+u^2X_2,\quad x=\left(
\begin{array}{c}
x^1 \\
x^2
\end{array}
\right) \in R^2,\quad X_1=\left(
\begin{array}{c}
1 \\
0
\end{array}
\right) ,\quad X_2=\left(
\begin{array}{c}
0 \\
x^1
\end{array}
\right)
\]
\[
\min_{u(.)}\int_0^T\mathcal{F}(u(t))dt,\quad \mathcal{F}(u)=\left\|
u\right\| +\left\langle b,u\right\rangle ,\ b=(\varepsilon ,0)^t,\
u=(u^1,u^2)^t,\ 0\leq \varepsilon <1.
\]
\[
x(0)=0,\quad x(T)=x_T.
\]
We are looking for the geodesics starting from the origin and
parameterized by arclength. The distribution $D=<$ $X_1,X_2>$ is
bracket generating and has not a constant rank on $\Bbb{R}^2$. If
we take the regular Lagrangian
\[
\mathcal{L}=\frac 12\mathcal{F}^2=\frac 12\left(
\sqrt{u_1^2+u_2^2} +\varepsilon u_1\right) ^2,
\]
on $E=D$ and use a result from \cite{Mi5} we obtain a regular
Hamiltonian $ \mathcal{H}$ on $E^{*}$
\[
\mathcal{H}=\frac 12\left( \sqrt{\frac{(\mu _1)^2}{(1-\varepsilon
^2)^2}+ \frac{(\mu _2)^2}{(1-\varepsilon ^2)}}-\frac{\varepsilon
\mu _1}{ 1-\varepsilon ^2}\right) ^2,
\]
and from (2.3.6) we get
\[
\left(
\begin{array}{c}
\mu _1 \\
\mu _2
\end{array}
\right) =\left(
\begin{array}{cc}
1 & 0 \\
0 & x^1
\end{array}
\right) \left(
\begin{array}{c}
p_1 \\
p_2
\end{array}
\right) ,
\]
and the corresponding Hamiltonian $H$ on $T^{*}M$ is
\[
H=\frac 12\left( \sqrt{\frac{(p_1)^2}{(1-\varepsilon ^2)^2}+\frac{
(p_2)^2(x^1)^2}{(1-\varepsilon ^2)}}-\frac{\varepsilon
p_1}{1-\varepsilon ^2} \right) ^2,
\]
From the Hamilton equations (2.3.7) we obtain
\[
\dot x^1=\frac{(1+\varepsilon )^2p_1}{\left( 1-\varepsilon
^2\right) ^2} -\frac \varepsilon {1-\varepsilon
^2}\sqrt{\frac{(p_1)^2}{(1-\varepsilon
^2)^2}+\frac{a^2(x^1)^2}{1-\varepsilon ^2}}-\frac{\varepsilon
p_1^2}{\left( 1-\varepsilon ^2\right) ^3}\frac
1{\sqrt{\frac{(p_1)^2}{(1-\varepsilon ^2)^2}
+\frac{a^2(x^1)^2}{1-\varepsilon ^2}}}
\]
\begin{equation}
\dot x^2=\frac{(x^1)^2a}{1-\varepsilon ^2}-\frac{\varepsilon
a(x^1)^2p_1}{ (1-\varepsilon ^2)^2}\frac
1{\sqrt{\frac{(p_1)^2}{(1-\varepsilon ^2)^2}+
\frac{a^2(x^1)^2}{1-\varepsilon ^2}}}  \tag{2.3.8}
\end{equation}
\[
\stackrel{\cdot }{p_1}=-\frac{x^1a^2}{1-\varepsilon ^2}+\frac{\varepsilon
p_1x^1a^2}{(1-\varepsilon ^2)^2}\frac 1{\sqrt{\frac{(p_1)^2}{(1-\varepsilon
^2)^2}+\frac{a^2(x^1)^2}{1-\varepsilon ^2}}},\qquad p_2=a=ct.
\]
We make the following change of variables
\[
\left\{
\begin{array}{l}
x^1=\frac{\sqrt{1-\varepsilon ^2}}ar(t)\sin A\theta (t), \\
\\
p_1=(1-\varepsilon ^2)r(t)\cos A\theta (t),
\end{array}
\right.
\]
and from (2.3.8) we get
\[
\frac{A(1-\varepsilon ^2)(\sqrt{1-\varepsilon ^2}}a\frac{d\theta
}{dt} =(1-\varepsilon \cos A\theta )^2,
\]
\[
\frac{(1-\varepsilon ^2)(\sqrt{1-\varepsilon ^2}}a\frac{dr}{dt}=\varepsilon
r\sin A\theta (\varepsilon \cos A\theta -1).
\]
Hence
\[
r=\frac 1{c(1-\varepsilon \cos A\theta )},\qquad c\in R
\]
and
\[
t=\frac{A(1-\varepsilon ^2)\sqrt{1-\varepsilon ^2}}a\int
\frac{d\theta }{ \left( 1-\varepsilon \cos A\theta \right) ^2}.
\]
But the geodesics are parameterized by arclength, that corresponds
to fix the level $1/2$ of the Hamiltonian and we have
\[
H=\frac{r^2}2\left( 1-\varepsilon \cos A\theta \right) ^2=\frac 1{2c^2},
\]
so $c=\pm 1$ and therefore
\[
r=\pm \frac 1{1-\varepsilon \cos A\theta }
\]
and finally we obtain
\begin{equation}
x^1=\pm \sqrt{1-\varepsilon ^2}\frac{\sin A\theta }{a(1-\varepsilon \cos
A\theta )},\quad \quad  \tag{2.3.9}
\end{equation}
\[
x^2=\frac{A(1-\varepsilon ^2)\sqrt{1-\varepsilon ^2}}a\int \frac{\sin
^2A\theta }{\left( 1-\varepsilon \cos A\theta \right) ^3}d\theta .
\]

\begin{remark}
1.The above geodesics are the only minimizers of this sub-Finslerian
geometry.

2. For $\varepsilon =0$ we obtain the geodesics of the Grushin plane endowed
with the standard Euclidean metric, i.e. a sub-Riemannian geometry
(distributional systems with quadratic cost)
\begin{equation}
x_1=\pm \frac{\sin at}a,\quad x_2=\frac t{2a}-\frac{\sin 2at}{4a^2}.
\tag{2.3.10}
\end{equation}
(\cite{Bel, Fa}).
\end{remark}

\newpage\ \quad

\subsubsection{\textbf{Heisenberg case}}

We consider the following distributional system with positive homogeneous
cost (Heisenberg group) \cite{Hr3}
\[
\dot x=u^1X_1+u^2X_2,\quad x=\left(
\begin{array}{c}
x^1 \\
x^2 \\
x^3
\end{array}
\right) \in \Bbb{R}^3,\ X_1=\left(
\begin{array}{c}
1 \\
0 \\
-\frac{x^2}2
\end{array}
\right) ,\ X_2=\left(
\begin{array}{c}
0 \\
1 \\
\frac{x^1}2
\end{array}
\right) ,
\]
\[
\min_{u(.)}\int_0^T\mathcal{F}(u(t))dt,\quad \mathcal{F}(u)=\left\|
u\right\| +\left\langle b,u\right\rangle ,\ b=(\varepsilon ,0)^t,\
u=(u^1,u^2)^t,\ 0\leq \varepsilon <1.
\]
We are looking for the geodesics starting from the origin and
parameterized by arclength. Here
\[
X_1=\frac \partial {\partial x^1}-\frac{x^2}2\frac \partial {\partial x^3},
\]
\[
X_2=\frac \partial {\partial x^2}+\frac{x^1}2\frac \partial {\partial x^3},
\]
\[
X_3=[X_1,X_2]=\frac \partial {\partial x^3}
\]
and $\langle X_1,X_2,X_3\rangle \equiv T\Bbb{R}^3,$ hence the
distribution $ D=<$ $X_1,X_2>$ of constant rank is strong bracket
generating \cite{Str1}.

\begin{remark}
We can try to work directly on the cotangent bundle by computing the
Hamiltonian $H(x,p)$ $=\mathcal{H}(x,\mu )$ , $\mu =\sigma ^{*}(p).$ Since
\begin{equation}
\left(
\begin{array}{c}
\mu _1 \\
\mu _2
\end{array}
\right) =\left(
\begin{array}{ccc}
1 & 0 & -x^2/2 \\
0 & 1 & x^1/2
\end{array}
\right) \left(
\begin{array}{c}
p_1 \\
p_2 \\
p_3
\end{array}
\right) ,  \tag{2.3.11}
\end{equation}
we obtain
\begin{equation}
H=\frac 12\left( \sqrt{\frac{\left( p_1-\frac{p_3x^2}2\right) ^2}{
(1-\varepsilon ^2)^2}+\frac{\left( p_2+\frac{p_3x^1}2\right) ^2}{
1-\varepsilon ^2}}-\frac{\varepsilon
(p_1-\frac{p_3x^2}2)}{1-\varepsilon ^2} \right) ^2.  \tag{2.3.12}
\end{equation}
Unfortunately, with this Hamiltonian (2.3.7) is a very complicated system of
implicit differential equations.
\end{remark}

We will use a different approach. Let us take $M=\Bbb{R}^3$ and
$E=TM$ with the basis $\{X_{1,}X_2,X_3\}$. $E$ is a Lie algebroid
over $M$ with at least one structural function nonzero. The anchor
$\sigma :E\rightarrow TM$ is just the identity and the matrix of
$\sigma $ with respect to the basis of $ E $ and $TM$ basis is
\begin{equation}
\sigma _\alpha ^i=\left(
\begin{array}{ccc}
1 & 0 & 0 \\
0 & 1 & 0 \\
-\frac{x^2}2 & \frac{x^1}2 & 1
\end{array}
\right)  \tag{2.3.13}
\end{equation}
Now the above control system can be written
\[
\dot x=u^1X_1+u^2X_2+0X_3,\quad x=\left(
\begin{array}{c}
x^1 \\
x^2 \\
x^3
\end{array}
\right) \in \Bbb{R}^3,
\]
\[
\min_{u(.)}\int_0^T\mathcal{F}(u(t))dt=\min_{u(.)}\int_0^T\mathcal{L}
(u(t))dt,
\]
where $\mathcal{F}=\sqrt{(u^1)^2+(u^2)^2}+\varepsilon u^1$ and $\mathcal{L}$
$=\frac 12\mathcal{F}^2$. To solve this problem we form the augmented
Lagrangian
\[
\stackrel{\sim }{\mathcal{L}}\mathcal{=L}(u)+\lambda u^3,
\]
($\lambda $ is Lagrange a multiplier) and still work via the maximum
principle but at the level of the Lie algebroid. If we set
\[
\mu _i=\frac{\partial \stackrel{\sim }{\mathcal{L}}}{\partial u^i},\quad
\mathcal{H}=\mu _iu^i-\stackrel{\sim }{\mathcal{L}},
\]
we obtain
\[
\mu _1=\frac{\partial \mathcal{L}}{\partial u^1},\ \mu _2=\frac{\partial
\mathcal{L}}{\partial u^2},\ \mu _3=\lambda ,
\]
and therefore
\[
\mathcal{H}=\mu _iu^i-\stackrel{\sim }{\mathcal{L}}=\frac{\partial
\mathcal{L }}{\partial u^1}u^1+\frac{\partial
\mathcal{L}}{\partial u^2}u^2+\lambda u^3- \mathcal{L}-\lambda
u^3=\mathcal{L},
\]
because $\mathcal{L}$ is 2-homogeneous. Again, using a result from
\cite{Mi5} , the Hamiltonian on $E^{*}$ is given by
\begin{equation}
\mathcal{H}=\frac 12\left( \sqrt{\frac{(\mu _1)^2}{(1-\varepsilon
^2)^2}+ \frac{(\mu _2)^2}{(1-\varepsilon ^2)}}-\frac{\varepsilon
\mu _1}{ 1-\varepsilon ^2}\right) ^2,\quad \mu _3=\lambda .
\tag{2.3.14}
\end{equation}
From the structure equations of Lie algebroids (1.2.8) and the
relation $ [X_\alpha ,X_\beta ]=L_{\alpha \beta }^\gamma X_\gamma
$ we obtain the non-zero structural functions
\[
L_{12}^3=1,\quad L_{21}^3=-1.
\]
Now, using the Hamilton equations on Lie algebroids (2.1.3) we get the
following systems of differential equations
\begin{equation}
\dot x^1=\frac{\partial \mathcal{H}}{\partial \mu _1},\quad \dot
x^2=\frac{
\partial \mathcal{H}}{\partial \mu _2},\quad \quad \quad \quad \quad \quad
\tag{2.3.15}
\end{equation}
\[
\dot x^3=-\frac{x^2}2\frac{\partial \mathcal{H}}{\partial \mu
_1}+\frac{x^1}2 \frac{\partial \mathcal{H}}{\partial \mu
_2}=\frac{x^1\dot x^2}2-\frac{\dot x^1x^2}2
\]
and
\[
\ \dot \mu _1=-L_{12}^3\mu _3\frac{\partial \mathcal{H}}{\partial
\mu _2} =-\lambda \frac{\partial \mathcal{H}}{\partial \mu _2}
\]
\begin{equation}
\ \dot \mu _2=-L_{21}^3\mu _3\frac{\partial \mathcal{H}}{\partial
\mu _1} =\lambda \frac{\partial \mathcal{H}}{\partial \mu _1}
\tag{2.3.16}
\end{equation}
\[
\ \dot \mu _3=0\Rightarrow \lambda =ct.
\]
where
\begin{eqnarray*}
&&\ \ \ \ \ \ \ \ \ \ \ \left. \frac{\partial
\mathcal{H}}{\partial \mu _1}= \frac{\left( 1+\varepsilon
^2\right) \mu _1}{(1-\varepsilon ^2)^2}-\frac{ \varepsilon
\sqrt{\frac{(\mu _1)^2}{(1-\varepsilon ^2)^2}+\frac{(\mu _2)^2}{
1-\varepsilon ^2}}}{1-\varepsilon ^2}-\frac{\varepsilon \mu _1^2}{
(1-\varepsilon ^2)^3\sqrt{\frac{(\mu _1)^2}{(1-\varepsilon
^2)^2}+\frac{(\mu
_2)^2}{1-\varepsilon ^2}}}\right. \\
&&\ \ \ \ \ \ \ \ \ \ \ \left. \frac{\partial
\mathcal{H}}{\partial \mu _2}= \frac{\mu _2}{1-\varepsilon
^2}-\frac{\varepsilon \mu _1\mu _2}{ (1-\varepsilon
^2)^2\sqrt{\frac{(\mu _1)^2}{(1-\varepsilon ^2)^2}+\frac{(\mu
_2)^2}{1-\varepsilon ^2}}}.\right.
\end{eqnarray*}
We may use to the following transformations
\begin{equation}
\begin{array}{c}
\mu _1(t)=(1-\varepsilon ^2)r(t)(a\cos A\theta (t)-b\sin A\theta (t)) \\
\\
\mu _2(t)=\sqrt{1-\varepsilon ^2}r(t)(a\sin A\theta (t)+b\cos A\theta (t))
\end{array}
\tag{2.3.17}
\end{equation}
such that $a^2+b^2=1$ and we get
\[
\sqrt{\frac{(\mu _1)^2}{(1-\varepsilon ^2)^2}+\frac{(\mu _2)^2}{
(1-\varepsilon ^2)}}=\left| r\right| .
\]
Also, from (2.3.16) we obtain
\[
t=\frac{(1-\varepsilon ^2)\sqrt{1-\varepsilon ^2}A}\lambda \int
\frac{ d\theta }{\left( 1-\varepsilon \left( a\cos A\theta -b\sin
A\theta \right) \right) ^2},
\]
and
\[
r=\frac 1{c\left( 1-\varepsilon \left( a\cos A\theta -b\sin A\theta \right)
\right) }.
\]
But the geodesics are parameterized by arclength, that corresponds
to fix the level $1/2$ of the Hamiltonian and we have
\[
\mathcal{H}=\frac{r^2}2\left( 1-\varepsilon \left( a\cos A\theta -b\sin
A\theta \right) \right) ^2=\frac 1{2c^2}.
\]
so $c=\pm 1$ and therefore
\[
r=\pm \frac 1{1-\varepsilon \left( a\cos A\theta -b\sin A\theta \right) }.
\]
From (2.3.17) we obtain
\begin{eqnarray*}
\mu _1(t) &=&\pm \frac{(1-\varepsilon ^2)(a\cos A\theta -b\sin
A\theta )}{
1-\varepsilon \left( a\cos A\theta -b\sin A\theta \right) }, \\
&& \\
\mu _2(t) &=&\frac{\sqrt{1-\varepsilon ^2}(a\sin A\theta +b\cos
A\theta )}{ 1-\varepsilon \left( a\cos A\theta -b\sin A\theta
\right) }.
\end{eqnarray*}
But
\[
\ \dot \mu _2=\lambda \dot x^1,
\]
so
\[
x^1(\theta )=\frac{\mu _2}\lambda -a_1.
\]
Since we are looking for geodesics starting from the origin, we
have $a_1= \frac{\sqrt{1-\varepsilon ^2}b}{\lambda (1-\varepsilon
a)}$ and therefore
\[
x^1(\theta )=\frac{\sqrt{1-\varepsilon ^2}(a\sin A\theta +b\cos
A\theta )}{ \lambda \left( 1-\varepsilon \left( a\cos A\theta
-b\sin A\theta \right) \right) }-\frac{\sqrt{1-\varepsilon
^2}b}{\lambda (1-\varepsilon a)},
\]
and from
\[
\ \dot \mu _1=-\lambda \dot x^2,
\]
we obtain
\[
x^2(\theta )=\frac{(1-\varepsilon ^2)(b\sin A\theta -a\cos A\theta
)}{ \lambda \left( 1-\varepsilon \left( a\cos A\theta -b\sin
A\theta \right) \right) }+\frac{(1-\varepsilon ^2)a}{\lambda
(1-\varepsilon a)}.
\]
Finally, from
\[
\dot x^3=\frac{x^1\dot x^2}2-\frac{\dot x^1x^2}2,
\]
by straightforward computation we get
\[
x^3(\theta )=\frac{(1-\varepsilon ^2)\sqrt{1-\varepsilon ^2}A}{2\lambda
^2(1-\varepsilon a)}\int \frac{1-\cos A\theta }{\left( 1-\varepsilon \left(
a\cos A\theta -b\sin A\theta \right) \right) ^2}d\theta .
\]

\begin{remark}
For $\varepsilon =0$ we obtain the sub-Riemannian case (distributional
systems with quadratic cost) with the solutions
\[
x(t)=\frac{a\sin \lambda t-b(1-\cos \lambda \theta )}\lambda
\]
\[
y(\theta )=\frac{b\sin \lambda t+a(1-\cos \lambda t)}\lambda
\]
\[
z(\theta )=\frac{\lambda t-\sin \lambda t}{2\lambda ^2},\quad a^2+b^2=1
\]
(see also \cite{Bel, Fa}).
\end{remark}

\end{document}